%% file: Article.tex
\documentclass[hidelinks,onefignum,onetabnum]{siamart220329}

\input{ex_shared}

\ifpdf
\hypersetup{
  pdftitle={Stopping Criteria for the Conjugate Gradient Algorithm in High-Order Finite Element Methods},
  pdfauthor={}
}
\fi

\colorlet{shadecolor}{blue!20}
\begin{document}

\maketitle

\begin{abstract}
We \added{consider}\deleted{introducethree new} stopping criteria that balance algebraic and discretization errors for the conjugate gradient algorithm applied to high-order finite element discretizations of Poisson problems.
\deleted{The current state of the art stopping criteria compare a posteriori estimates of discretization error against estimates of the algebraic error. Firstly, we propose a new error indicator derived from a recovery-based error estimator that is less computationally expensive and more reliable.
Secondly,}\added{Firstly, }we introduce a new stopping criterion that suggests stopping when the norm of the \deleted{linear residual}\added{linear system residual} is less than a small fraction of an error indicator derived directly from the residual. This indicator shares the same mesh size and polynomial degree scaling  as the norm of the residual, resulting in a robust criterion regardless of the mesh size,  the polynomial degree, and the shape regularity of the mesh.
\added{Secondly}\deleted{Thirdly}, 
for solving Poisson problems with highly variable piecewise constant coefficients, we introduce a subdomain-based criterion that recommends stopping when the norm of the \deleted{linear residual}\added{linear system residual} restricted to each subdomain is smaller than the corresponding indicator also restricted to that subdomain. 
\added{Reliability and efficiency theorems for the first criterion are established.}
Numerical experiments, including tests with \deleted{anisotropic meshes and }highly variable piecewise constant coefficients \added{and a GPU-accelerated three-dimensional elliptic solver}, demonstrate that the proposed criteria efficiently avoid both premature termination and over-solving.
\end{abstract}



\section{Introduction}

\input{1-introduction.tex}

\input{2-1-formulation}

\input{2-2-APosterioriError}

\input{3-etaRF.tex}
\input{3-3}

\input{4-NumericalResults}

\input{5-conclusion}

\input{6-appendix}

\bibliographystyle{siamplain}
\bibliography{references}

\end{document}

%% file: ex_shared.tex
\usepackage{lipsum}
\usepackage{amsfonts}
\usepackage{graphicx}
\usepackage{epstopdf}
\usepackage{algorithmic}
\usepackage{amsmath,amsfonts,amssymb}
\usepackage{multirow}
\usepackage{makecell}
\usepackage{pgfplots}
\usepackage{appendix}
\usepackage[a4paper,vmargin=3cm]{geometry}
\usepackage{float}
\usepackage{tikz}
\usepackage{subcaption}
\usepackage{caption}
\usepackage[T1]{fontenc}
\usepackage[utf8]{inputenc}
\usepackage{lmodern}
\usepackage{babel}
\usepackage{enumerate}
\usepackage[shortlabels]{enumitem}
\usetikzlibrary{matrix}
\usepgfplotslibrary{groupplots}
\pgfplotsset{compat=newest}
\usetikzlibrary{shapes.geometric, arrows}
\usepackage{xcolor}         
\usepackage{cleveref}       
\usepackage{microtype}

\makeatletter
\newcommand{\deleted}[1]{%
}
\makeatother
\newcommand{\added}[1]{\textcolor{black}{#1}} 
\usepackage{filecontents}   

\def\lamminm{\lambda_{{\min}}\left(\BM\right)}
\def\lammaxm{\lambda_{{\max}}\left(\BM\right)}

\def\lamminml{\lambda_{{\min}}\left(\BM_{\ell}\right)}
\def\lammaxml{\lambda_{{\max}}\left(\BM_{\ell}\right)}
\def\lammina{\lambda_{{\min}}\left(\BA\right)}
\def\lammaxa{\lambda_{{\max}}\left(\BA\right)}

\def\hatk{\widehat{K}}

\def\underC{\underline{C}}
\def\overC{\overline{C}}

\tikzstyle{startstop} = [rectangle, minimum width=3cm, minimum height=1cm,text centered, draw=black]

\tikzstyle{process} = [rectangle, minimum width=3cm, minimum height=1cm, text centered, text width=3cm, draw=black]

\tikzstyle{decision} = [rectangle, rounded corners, minimum width=3cm, minimum height=1cm, text centered, text width=3cm, draw=black]
\tikzstyle{arrow} = [thick,->,>=stealth]

\ifpdf
  \DeclareGraphicsExtensions{.eps,.pdf,.png,.jpg}
\else
  \DeclareGraphicsExtensions{.eps}
\fi

\theoremstyle{definition}
\newtheorem{exmp}{Example}[subsection]

\definecolor{myindigo}{rgb}{0.2, 0.13, 0.53}
\definecolor{mycyan}{rgb}{0.52, 0.8, 0.93}
\definecolor{myteal}{rgb}{0.26, 0.67, 0.6}
\definecolor{mygreen}{rgb}{0.067, 0.4667, 0.2}
\definecolor{myolive}{rgb}{0.6, 0.6, 0.2}
\definecolor{mysand}{rgb}{0.86, 0.8, 0.46}
\definecolor{myrose}{rgb}{0.8, 0.4, 0.46}
\definecolor{mywine}{rgb}{0.52, 0.13, 0.33}
\definecolor{mypurple}{rgb}{0.67, 0.26, 0.6}

\newsiamremark{remark}{Remark}
\newsiamremark{assumption}{Assumption}
\newsiamremark{hypothesis}{Hypothesis}
\crefname{hypothesis}{Hypothesis}{Hypotheses}
\newsiamthm{claim}{Claim}

\headers{CG Stopping Criteria for High-Order Finite Element Methods}{Guo, de Sturler, and Warburton}

\title{Stopping Criteria for the Conjugate Gradient Algorithm in High-Order Finite Element Methods
\thanks{Submitted to the editor May 12 2023.
\funding{
This work was supported under NSF DMS 2208470.
The third author was supported in part by the Exascale Computing Project, a collaborative effort of two U.S. Department of Energy organizations (Office of Science and the National Nuclear Security Administration) responsible for the planning and preparation of a capable exascale ecosystem, including software, applications, hardware, advanced system engineering, and early testbed platforms, in support of the nation’s exascale computing imperative. He was also supported in part by the John K. Costain Faculty Chair in Science at Virginia Tech. 
The first author acknowledges the generous support from the John K. Costain Graduate Fellowship. 
}}}

\author{Yichen Guo\thanks{Department of Mathematics,  Virginia Tech, 
Blacksburg, VA 24061
	 (\email{ycguo@vt.edu}, \email{sturler@vt.edu}, \email{tcew@vt.edu})} 
  \and Eric de Sturler\footnotemark[2] 
  \and Tim Warburton\footnotemark[2]
}

\usepackage{amsopn}

\def\C{\mathbb{C}}

\def\CNN{\C^{N\times N}\kern-2pt}

  \def\BA{{\bf A}} 
\def\Bb{{\bf b}}   
   
\def\Bd{{\bf d}}   
\def\Be{{\bf e}}   
  \def\BF{{\bf F}}

  \def\BL{{\bf L}} 
  \def\BM{{\bf M}} 
\def\Bn{{\bf n}}   
   
\def\Bp{{\bf p}}  \def\BP{{\bf P}} \def\CP{{\cal P}}
   
\def\Br{{\bf r}}  \def\BR{{\bf R}}

\def\Bu{{\bf u}}   
\def\Bv{{\bf v}}  \def\BV{{\bf V}} 
\def\Bw{{\bf w}}   
\def\Bx{{\bf x}}   
\def\By{{\bf y}}   
\def\Bz{{\bf z}}

\def\BAs{\BA{\kern-1pt}}      
\def\BPs{\BP{\kern-1.2pt}}    
\def\BVs{\BV{\kern-1.2pt}}    
\def\CPs{\CP{\kern-1.2pt}}    

\def\Vhn{\mathcal{V}_{h,N}}

\def\diam{\text{diam}}

\usepackage{framed}


%% file: 1-introduction.tex
Solving linear elliptic partial differential equations (PDEs) involves two main steps: discretization of the PDE and solving the resulting algebraic linear system. 
As a result, two primary sources of error 
emerge: discretization error and algebraic error, which results from the iterative solution of the linear system.
The efficient termination of iterative solvers achieves a balance between discretization error and algebraic error.
Ideally, a stopping criterion for the iterative solver suggests stopping the iteration when the algebraic error is \added{marginally lower than}\deleted{dominated by} the discretization error.
A desirable stopping criterion should be reliable in the sense of maintaining the overall accuracy of the finite element solution, and efficient in the sense of terminating the iterative solver as early as possible.
Moreover, the criterion should be inexpensive to compute and the computation should be memory efficient.
In this paper, we consider the Poisson problem discretized with high-order finite element methods (FEM), and we  solve the linear system using the Conjugate Gradient algorithm (CG).
For this study, we assume the finite element space has already been fixed, and we aim to iterate until the error due to the given discretization \deleted{dominates}\added{is slightly greater than} the error due to the linear system.

The design of stopping criteria in finite element frameworks has been explored in numerous papers \cite{arioli2004stopping, arioli2005stopping, picasso2009stopping, jiranek2010posteriori, arioli2013interplay, ern2013adaptive, arioli2013stopping, papevz2020sharp}. 
One commonly adopted criterion in these works involves assessing when the ratio of estimated algebraic error to the estimated total error falls below a threshold.
Algebraic error estimation was discussed early on by Hestenes and Stiefel in \cite{hestenes1952methods}, and further developed in \cite{golub1997matrices, meurant2013computing, meurant2021accurate, meurant2006lanczos, golub2009matrices}. 
It is common for this type of algebraic error estimation to rely on computing the difference between the computed solutions at two different iterations with a heuristically chosen gap between the iterations. 
If the iterative convergence rate is slow, then a larger delay (i.e. a greater number of additional iterations) may be necessary. 
As indicated by numerical experiment 4.1.2 presented in \cite{arioli2004stopping}, it can be  challenging to determine a reliable delay parameter for the Poisson problem with a highly variable  coefficient, due to the potential need to use a large delay to compensate for slow iterative convergence.

\added{The discretization error, which is another source of the total error, can be estimated using a posteriori error estimators.}
\deleted{
Monitoring the total error in finite element methods is also a key ingredient when devising stopping criteria.
In the context of our study, we refer to an error indicator as an ``error estimator'' if it is equivalent to the exact total error. Specifically, an error estimator should provide both upper and lower bounds to the total error up to a constant independent of the polynomial degree and the mesh size. }
Babu\v{s}ka and Rheinboldt \cite{babuvvska1978error} proposed a residual-based a posteriori error estimator for low-order FEM  on one-dimensional domains in the late 1970s. Subsequently, Melenk and Wohlmuth \cite{melenk2001residual} extended the estimator for $hp$\text{-}FEM by generalizing the Cl\'{e}ment interpolation operator to the $hp$-finite element discretization. However, their estimator is an upper bound \added{on}\deleted{for} the \added{discretization} error only up to an unknown constant that depends on the shape regularity of the triangulation, which may 
lead to a significant overestimation of \added{errors.}\deleted{the error.}
Methods for estimating the unknown constant, as developed in \cite{carstensen1999fully, veeser2009explicit}, require solving local eigenvalue problems, or obtaining trace inequalities and Poincar\'e-type inequalities with explicit constants.
On the other hand, flux recovery error estimation techniques \cite{bastian2003superconvergence, bank2003asymptotically, zhang2005new, cai2017improved} introduced in \cite{zienkiewicz1987simple, zienkiewicz1992superconvergent} reconstruct an approximation to the flux and compare the reconstructed flux with the numerical flux. The efficiency of this approach is robust with respect to the polynomial degree; however, solving this requires a significant amount of computation and memory for high-order finite element approximation.
\added{
All previously mentioned estimators require the exact solution to the linear system, which is unavailable. In \cite{becker2009convergence, arioli2013stopping}, authors estimate the discretization error using estimators based on the approximate solution. This error estimate is not equivalent to the exact discretization error due to the absence of the Galerkin orthogonality assumption. A discussion on removing this assumption is provided in \cite{papevz2018residual}.}
\deleted{
Our goal is to design a reliable and efficient stopping criterion that is also  robust with respect to the mesh size, the polynomial degree, the shape regularity of the mesh, and the diffusion coefficient. }
\deleted{To design such a stopping criterion, we propose three main innovations.
Firstly, as for stopping criteria based on a comparison of the algebraic error estimator and the a posteriori error estimator 
\mbox{ \cite{arioli2004stopping, picasso2009stopping, ern2013adaptive, arioli2013stopping}, }
we propose an alternative error indicator that serves as a lower bound \added{on}\deleted{for} the recovery-based error estimator, using the singular value decomposition of elliptic lifting matrices arising from local problems in flux recovery. The indicator is computationally less expensive to evaluate while iterating and more reliable than the recovery-based error estimator for high-order elements.}

\added{
Our goal is to design a reliable and efficient stopping criterion that is  robust with respect to the mesh size, the polynomial degree, the shape regularity of the mesh, and the diffusion coefficient. 
To design such a criterion, we propose two main innovations.
Firstly,}\deleted{Secondly,} in contrast to criteria comparing error estimates, we propose a simplified stopping criterion that depends on the norm of the \deleted{linear residual}\added{linear system residual} and an error indicator for the Poisson problem with a constant diffusion coefficient. 
We decompose the \deleted{linear residual}\added{linear system residual} into a component corresponding to the strong residual tested against the basis functions and a second component corresponding to the jumps in the normal gradient at element interfaces also tested against the basis functions. We then apply the triangle inequality to derive an error indicator that is directly comparable to the norm of the \deleted{linear residual}\added{linear system residual}. 
This indicator tends to stagnate when the discretization error \deleted{dominates}\added{is above} the algebraic error, as it depends on the strong residual and jumps in the normal gradient.
Therefore, the divergence of this indicator from the norm of the \deleted{linear residual}\added{linear system residual} can be an effective proxy for identifying when the discretization error \deleted{dominates}{is greater than} the algebraic error.
This observation motivates a criterion for terminating the iterative method when the ratio of the norm of the \deleted{linear residual}\added{linear system residual} to the new indicator falls below a specific tolerance.
The proposed indicator is a natural upper bound \added{on}\deleted{for} the norm of the residual without any unknown constants to be estimated.
It has the same intrinsic mesh size and polynomial degree scaling as the  norm of the \deleted{linear residual}\added{linear system residual}, which coincides with the scaling of the energy norm of the error in two dimensions. 
Moreover, compared with criteria based on error estimation, the proposed criterion does not require estimating the algebraic error since it relies on the \deleted{linear residual}\added{linear system residual}.
Furthermore, separate computation of the component corresponding to jumps in the normal gradient  is unnecessary, as it can be obtained directly from the difference between the \deleted{linear residual}\added{linear system residual} and the component corresponding to the strong residual. By contrast, both strong element residuals and jump residuals are computed in residual a posteriori estimators \cite{babuvvska1978error, melenk2001residual}.

\deleted{Finally}\added{Secondly}, it is important to note that the diffusion coefficient scaling in the norm of the \deleted{linear residual}\added{linear system residual} and the new indicator is different from the scaling in the total error for problems with highly variable diffusion coefficient.
This implies that contributions from subdomains with small coefficients may be dominated by contributions from 
subdomains with large coefficients. 
Thus, when solving the Poisson equation with highly variable coefficients, the separation of the indicator and the norm of the \deleted{linear residual}\added{linear system residual} may occur at a different iteration than the point at which the discretization error \deleted{dominates over}\added{is greater than} the algebraic error.
To address this issue, we propose a subdomain-based criterion that only 
recommends stopping when the norm of the \deleted{linear residual}\added{linear system residual} restricted to each subdomain is relatively small compared to the indicator restricted to that subdomain. This approach ensures that the iteration achieves sufficient accuracy in all subdomains and provides a reliable stopping criterion for problems with highly variable piecewise constant coefficients.

The paper is organized as follows.
In \cref{sec:formulation}, we review stopping criteria based on a comparison of estimates of the algebraic error and a posteriori estimates of discretization error for high-order finite element methods. \deleted{Additionally, we propose an error indicator.}
In \cref{rnfn}, we introduce a new stopping criterion that  compares the norm of the residual \added{to}\deleted{and} an indicator, and a subdomain-based stopping criterion for problems with highly variable coefficients. 
\added{Furthuremore, we establish the reliability and efficiency theorems for the Poisson equation  with  constant coefficient.}
\added{
In \cref{num}, we provide numerical results including  Poisson problems with highly variable piecewise constant coefficients and a GPU-accelerated three-dimensional elliptic solver to demonstrate the effectiveness of the proposed stopping criteria.}
\deleted{
In \mbox{\cref{num}}, we provide numerical results to demonstrate the effectiveness of the stopping criteria.  
Numerical experiments include tests with anisotropic meshes, solutions with singularities, and Poisson problems with highly variable piecewise constant coefficients.}
We end with conclusions in \cref{sec:conclusions}.

Throughout this paper, we will use standard notation from Sobolev space theory. 
\added{For a bounded domain $D \subset \mathbb{R}^d$, $\left(\cdot, \cdot\right)_D$ and $\|\cdot\|_{D}$ denote the inner product and the associated norm on $L^2(D)$. Without a subscript, we use $\left(\cdot, \cdot\right)$ and $\|\cdot\|$ to represent the inner produt and norm on $L^2(\Omega)$.}
\deleted{For a bounded domain \mbox{$\Omega \subset \mathbb{R}^d$, $\left(\cdot, \cdot\right)$} and \mbox{$\|\cdot\|_{\Omega}$} denote the inner product and the associated norm on \mbox{$L^2(\Omega)$}. }
For a vector $\Bx\in \mathbb{R}^n$, $\|\Bx\|$ denotes the $l^2$ norm of $\Bx$.

%% file: 2-1-formulation.tex
\section{Formulation}
\label{sec:formulation}
We consider the Poisson problem
\begin{equation}
\label{eq:poisson}
    - \mathbf{\nabla} \cdot\left(\kappa(x) \nabla u(x) \right) = f(x) 
\end{equation}
on a bounded domain \added{{$\Omega\subset \mathbb{R}^2$}}\deleted{$\Omega\subset \mathbb{R}^d$}, with boundary conditions  
\[
\kappa(x)\frac{\partial u}{\partial n}=g \text{ on } \Gamma_N, \quad u=0 \text{ on } \Gamma_D,
\]
where $\Gamma_N \cap \Gamma_D=\varnothing, \overline{\Gamma}_N \cup \overline{\Gamma}_D=\partial \Omega$, $f\in L^2(\Omega)$, and $g\in L^2(\Gamma_N)$ describes the Neumann boundary condition. We assume there exists a constant \deleted{$\alpha$}\added{$\underline{\kappa}$} such that \added{$0<  \underline{\kappa} \leq \kappa(x)\in L^2(\Omega)$} \deleted{$0<  \alpha\leq \kappa(x)\in L^2(\Omega)$}. 

We define $H_{0, \Gamma_D}^1(\Omega):=\left\{v \in H^1(\Omega) : v|_{\Gamma_D}=0\right\}$. The weak formulation of the Poisson equation \eqref{eq:poisson} is: find $u\in H_{0, \Gamma_D}^1(\Omega)$, such that
\begin{equation}
\label{eq:weak}
   a(u,v) = \ell(v), \quad \forall v\in H_{0, \Gamma_D}^1(\Omega),
\end{equation}
where
\[
\begin{aligned}
 a\left(u ,v \right) & := \int_{\Omega} \kappa(x) \nabla u \cdot \nabla v \,
 dx, \quad u, v \in H_{0, \Gamma_D}^1(\Omega), \\
  l(v) & :=  \int_{\Omega} f v \,dx +  \int_{\Gamma_N }gv \,ds, \quad v \in H_{0, \Gamma_D}^1(\Omega).
\end{aligned}
\]
\deleted{
We use the notation $\|\cdot\|_E$ to denote the energy norm
$
\|v\|_E = \sqrt{a(v,v)}.$}

Given a family of regular affine triangulations $\mathcal{T}_h = \left\{ K \right\}$ of $\Omega$ with elements $K$. 
\added{
We define 
$$h_K = \diam(K),\quad h= \max\limits_K h_K,$$
and
\[
\rho_K = \sup\left\{\diam(B): B \text{ is a ball contained in } K\right\}.
\]
We assume the triangulation is quasi-uniform, i.e. there exist constants $\sigma_1, \sigma_2 >0$ independent of $h$ such that for all elements $K$
\begin{equation}
    \frac{h}{h_K}<\sigma_1, \quad \frac{h_K}{\rho_K}\leq \sigma_2.
\label{eq:quasi uniform}
\end{equation}
}
\added{
We denote the reference element by $\widehat{K}$ which can be either the reference square
\[
\hatk = (-1,1)^2
\]
or the reference triangle
\[
\widehat{K}=\{(x, y) \mid -1\leq x,y \leq 1,x+y \leq0\}.
\]
Each element $K$ is the image of the reference element under an affine map $F_K:\widehat{K}\to K$ with $J_K = \nabla F_K$.
We 
}\deleted{we}define the finite element space \added{$\Vhn$}\deleted{\mbox{$S_N(\mathcal{T}_h, \Omega)$}} of piecewise polynomials of degree $N$ 
\[
\added{\Vhn :=\left\{ v_h\in H_{0,\Gamma_D}^1(\Omega) : v_h|_{K} \in \Bar{\mathbb{P}}_N(K), K\in \mathcal{T}_h \right\},}
\]
\added{
where {$\Bar{\mathbb{P}}_N(K)$ = $\mathbb{P}_N(K)$}, the polynomials space on \mbox{$K$} of total degree no more than $N$,  for triangle elements and {$\Bar{\mathbb{P}}_N(K)$ = $\mathbb{Q}_N(K)$}, the polynomial space on \mbox{$K$} of degree in each variable no more than $N$, for quadrilateral elements. We denote by \mbox{$N_s$} the dimension of \mbox{$\Vhn$},  and by \mbox{$\phi_n$} basis functions of \mbox{$\Vhn$}, where \mbox{$n = 1, \ldots, N_s$}. In this work  \mbox{$\phi_n$} denotes the Lagrange interpolating basis function associated with the $n$-th node. We use Warp \& Blend nodes for the triangle \mbox{\cite{warburton2006explicit}} and Gauss-Legendre-Lobatto nodes for the quadrilateral.
}
\deleted{
We denote by  \mbox{$\mathbb{Q}_N(K)$} the linear space of polynomials of total degree \mbox{$\leq N$ } on \mbox{$K$}, by \mbox{$N_s$} the dimension of \mbox{$\Vhn$},  and by \mbox{$\phi_n$} basis functions of \mbox{$\Vhn$}, where \mbox{$n = 1, \ldots, N_s$}. In this work  \mbox{$\phi_n$} denotes the Lagrange interpolating basis function associated with the n-th node.Warp \& Blend nodes for the triangle \mbox{\cite{warburton2006explicit}}.
}
We refer to $\mathcal{E}$ as the set of all \deleted{$(d-1)$-dimensional} element edges \deleted{(faces in $\mathbb{R}^3$)} of $\mathcal{T}_h$. 
\added{
Furthermore, we define $\mathcal{E}_{\text{bd}}^N$ and  $\mathcal{E}_{\text{bd}}^D$ as the set of element edges that lie on $\Gamma_N$ and $\Gamma_D$, respectively. We then decompose $\mathcal{E}$ into $\mathcal{E}_{\text{bd}}^N$, $\mathcal{E}_{\text{bd}}^D$, and the interior set $\mathcal{E}_{\text{int}} = \mathcal{E} \setminus \left(\mathcal{E}_{\text{bd}}^N \cup \mathcal{E}_{\text{bd}}^D\right).$}
\deleted{
Furthermore, we define $\mathcal{E}_{\text{bd}}$ as the set of element edges that lie on $\Gamma_N$, and decompose $\mathcal{E}$ into $\mathcal{E}_{\text{bd}}$ and $\mathcal{E}_{\text{int}} = \mathcal{E} \setminus \mathcal{E}_{\text{bd}}.$}
The finite element approximation to \eqref{eq:poisson} is: find $u_h \in \Vhn$ such that 
\begin{equation}
\label{eq:fem}
   a \left( u_h, v \right) = l (v), \quad \forall v\in \Vhn.
\end{equation}
Equations \cref{eq:weak} and \cref{eq:fem} give rise to the Galerkin orthogonality condition
\begin{equation}
    \label{eq:galerkin orth}
    a(u-u_h, v) = 0 \quad \forall  v\in  \Vhn.
\end{equation}

The approximation problem \cref{eq:fem} is equivalent to the linear system:
\begin{equation}
\label{eq:linear eq}
    \BA \Bx = \Bb,
\end{equation}
where $\BA \in \mathbb{R}^{N_s\times N_s}$ and  $\Bb \in \mathbb{R}^{N_s}$ are defined as follows,
\[
\BA_{ij} = a(\phi_j, \phi_i) , \quad \Bb_i = \ell(\phi_i).
\]
The matrix $\BA$ is symmetric and positive definite. 
\deleted{We define the  $\BA$-norm of $\Bx$ as $\|\Bx\|_{\BA} = \left(\Bx^T \BA\Bx \right)^{1/2}$.}\added{We define the  $\BA$-norm of $\By$ as $\|\By\|_{\BA} = \left(\By^T \BA\By \right)^{1/2}$.}
We assume that $\Bx_k\in \mathbb{R}^{N_s}$ is an approximate solution to \cref{eq:linear eq} obtained by an iterative method at the $k$-step , which in turn provides an approximate finite element solution 
$
u_h^k = \sum_{i=1}^{N_s}x_i^k \phi_i.
$
\added{
We define the residual  as 
\begin{equation}
    \Br_k = \Bb-\BA\Bx_k.
    \label{eq:res}
\end{equation}
}
The total error, the discretization error, and the algebraic error are  denoted by 
$$e:= u - u_h^k, \quad e_{\text{dis}} := u - u_h, \quad e_{\text{alg}} := u_h - u_h^k,$$  respectively. From the relation 
$a(e_{\text{alg}}, e_{\text{alg}}) = \left(\Bx-\Bx_k\right)^T\BA \left(\Bx-\Bx_k\right)  $, we obtain
\begin{equation}
    \|e_{\text{alg}}\|_E = \|\Bx-\Bx_k\|_{\BA}.
    \label{cg_alg}
\end{equation}
\added{
Here we use the notation $\|\cdot\|_E$ to denote the energy norm
\[
\|v\|_E = \sqrt{a(v,v)}.
\]
}
The Galerkin orthogonality condition \cref{eq:galerkin orth} implies 
\[
\|e\|_E^2  = \|e_{\text{dis}}\|_E^2 + \|e_{\text{alg}}\|_E^2. 
\]
As the iteration proceeds, the algebraic error gradually approaches zero, leading the total error to converge to the discretization error. Ideally, the iteration is terminated when the discretization error is dominant in the total error, i.e., 
\begin{equation}
    \|e_{\text{alg}}\|_E \leq \tau \|e\|_E,
    \label{ideal sc}
\end{equation}
for a chosen tolerance $\tau$, where $0<\tau<\added{1/\sqrt{2}} \deleted{1}$.
Since the total error and the algebraic error are unknown in general, we use error estimators $\eta_{\text{alg}}$ and $\eta_{\text{total}}$ to estimate the energy norm of the algebraic error, $\|e_{\text{alg}}\|_E$, and the total error, $\|e\|_E$, respectively. 
\added{
Consequently, it motivates the following stopping condition
\begin{equation}
\eta_{\text{alg}} \leq \tau \eta_{\text{total}}.
    \label{stopping criterion}
\end{equation}
}
\added{
A good stopping criterion should meet the following conditions:
}
\begin{enumerate}
    \item \added{Reliability: It should not terminate the iteration too early, ensuring that when the stopping criterion is met, the optimal stopping condition \cref{ideal sc} is also satisfied. }
    \item \added{Efficiency: It should not continue the iteration longer than necessary. Once the optimal stopping condition \cref{ideal sc} is achieved, the stopping criterion should also indicate that the iteration can stop.}
    \item \added{Independence from $h$ and $N$: The performance of the criterion should be $h$ and $N$ independent.}
\end{enumerate}

We review \added{the} estimation of the algebraic error for the conjugate gradient algorithm in \cref{CG} and several estimators  for high-order finite element discretization error in \cref{survey of a posteriori}.

\subsection{Error estimation for the conjugate gradient algorithm}
\label{CG}
The conjugate gradient algorithm was introduced by  Hestenes and  Stiefel \cite{hestenes1952methods} in 1952, and they also proposed a method to estimate the error. In \cite{strakovs2002error}, Strako\v{s} and Tich\'{y} showed that the estimation proposed in \cite{hestenes1952methods} is numerically stable. 
For the sake of completeness,
we briefly discuss the conjugate gradient algorithm and the error estimator proposed in \cite{hestenes1952methods}. We use the error estimator of CG as the algebraic error estimator $\eta_{\text{alg}}$ because of the equivalence of the $\BA$-norm of CG error and the energy norm of the algebraic error \cref{cg_alg}.  A comprehensive summary of CG is given in \cite{meurant1999numerical}. 

The conjugate gradient algorithm is as follows. Given $\Bx_0$, $\Br_0 = \Bb - \BA\Bx_0$, $\Bp_0=\Br_0$. For $k=1,2,\dotsc,$
\begin{equation*}
   \begin{aligned}
    & \gamma_{k-1} = \frac{\|\Br_{k-1}\|^2}{\|\Bp_{k-1}\|_{\BA}^2}, &\quad 
    & \Bx_{k} =  \Bx_{k-1} + \gamma_{k-1} \Bp_{k-1}, \quad
    &\Br_{k} = \Br_{k-1} - \gamma_{k-1}\BA\Bp_{k-1}, \quad
    \\ 
    &\beta_k = \frac{\|\Br_k\|^2}{\|\Br_{k-1}\|^2},& \quad
    & \Bp_k = \Br_{k} + \beta_{k} \Bp_{k-1}.  \quad
    & \\
\end{aligned} 
\label{cg}
\end{equation*}
The algorithm computes directions $\Bp_i$ that are $\BA$-orthogonal, i.e.  
$\Bp_i^T\BA\Bp_j = 0$, $  i \neq j$.
The approximate solution at the $k$-th step is
\[
 \Bx_{k} = \Bx_0+\sum_{i=0}^{k-1} \gamma_i\Bp_i.
 \]
To illustrate the idea of the error estimation, we assume that the CG algorithm can be run for $N_s$ steps and the exact solution $\Bx$ satisfies
\begin{equation*}
   \Bx = \Bx_0+\sum_{i=0}^{N_s}  \gamma_i\Bp_i.
\end{equation*} 
The $\BA$-norm of the CG error is
\begin{equation*}
    \|\Bx - \Bx_k\|_{\BA} = \left( \sum_{i=k}^{N_s} \gamma_i^2 \|\Bp_i\|_{\BA}^2 \right)^{1/2}.
\end{equation*}
If the \added{\textit{delay parameter}} $d$ satisfies $\|\Bx - \Bx_{k+d}\|_{\BA} \ll \|\Bx - \Bx_k\|_{\BA}$, then as 
$$ \|\Bx - \Bx_k \|_{\BA}^2 = \|\Bx_{k+d} - \Bx_k\|_{\BA}^2 + \|\Bx - \Bx_{k+d}\|_{\BA}^2,$$ 
Hestenes and Stiefel \cite{hestenes1952methods} estimate the $\BA$-norm of the CG error $ \|\Bx - \Bx_{k}\|_{\BA}$ by 
\begin{equation}
\label{eq:alg err}
    \eta_{\text{alg}}(u_h^k) := \|\Bx_{k+d} - \Bx_k\|_{\BA}.
\end{equation}
Hence,  \added{$d$ additional}\deleted{additional $d$} iterations are required to compute the estimator at the $k$-th step.

It is challenging to choose $d$ in advance, since the parameter depends on the convergence rate of CG. To achieve the same accuracy, the slower CG converges, the larger $d$ has to be. If $ \alpha \|\Bx - \Bx_{k}\|_{\BA}^2 =  \|\Bx - \Bx_{k+d}\|_{\BA}^2$ with $\alpha\in (0,1)$, the \added{\textit{effectivity}} of \cref{eq:alg err} is
\[
 \frac{\eta_{\text{alg}}(u_h^k)}{\|\Bx - \Bx_{k}\|_{\BA}} =  \left(1 - \alpha\right)^{1/2}.
\]
We demonstrate in \cref{num} that, with $d=10$, $\eta_{\text{alg}}$ is a good estimator if the algebraic error decreases fast, while it is unsatisfactory for some problems where the error  remains almost constant for a number of iterations. An increase in $d$ improves the accuracy of the estimator; however, it also leads to an increased number of additional iterations, which is undesirable.

%% file: 2-2-APosterioriError.tex
\subsection{Survey of A Posteriori error estimators}
\label{survey of a posteriori}

In this subsection, we review  error estimators based on the residual and flux reconstruction.
To simplify notation, \added{for all $w_h\in \Vhn$, }we define the element residual, \added{$r_E(w_h): \Omega\to \mathbb{R}$}\deleted{$r_E: \Omega\to \mathbb{R}$}, and the edge residual, \added{$r_J(w_h): \mathcal{E} \to \mathbb{R}$} \deleted{$r_J: \mathcal{E} \to \mathbb{R}$}, by
\added{
\begin{align}
\label{rerj}
  &   r_E(w_h)\big|_{K}  = f + \mathbf{\nabla} \cdot\left(\kappa(x) \nabla w_h \right) \text{ in } K ,\\ 
   &   \left.r_J(w_h)\right|_\ell  = 
   \begin{cases} 
   - \left[ \left(\kappa(x) \nabla w_h\right) \cdot \mathbf{n}_\ell\right] & \text { if } \ell \in \mathcal{E}_{\text{int}},  \\ 
   g - \left(\kappa(x) \nabla w_h \right)\cdot \mathbf{n}_\ell  & \text { if } \ell \in \mathcal{E}_{\text{bd}}^N, \\
   0 & \text { if } \ell \in \mathcal{E}_{\text{bd}}^D,
   \end{cases} 
\end{align}
}
where we denote the jump of the normal component of $\Bu$ across the edge $\ell$ by $\left[ \Bu \cdot \Bn_{\ell} \right]$, and $\Bn_{\ell}$ is the \added{unit} outward normal vector.

\subsubsection{Residual estimate}
\label{residual est}
The first error estimator for lower-order FEM  was proposed by Babuška and Rheinboldt \cite{babuvvska1978error}, and it has become a widely-used estimator in the literature,
\added{
\begin{equation}
\label{eq:etaR}
\eta^2(u_h) = \sum_{K\in\mathcal{T}_h}  h_K^2\| r_E(u_h)\|_{
    K}^2+ \sum_{\ell\in \mathcal{E}} h_\ell \left\| r_J(u_h)\right\|_\ell^2 .
\end{equation}
}
Here $h_K$ is the diameter of $K$ and $h_l$ is the length of the edge $\ell$. It is proved that the estimator is an upper bound \added{on} \deleted{for} the exact \added{discretization} error up to a constant \added{$C_{\text{BR}}$}\deleted{$C$},
\added{
\begin{equation}
\begin{aligned}
\| u - u_h\|_E \leq C_{\text{BR}} \eta(u_h) , 
\end{aligned}
\label{eq:etaconst}
\end{equation}
}
where \added{$C_{\text{BR}}$}\deleted{$C$} is independent of $h_K$. However, the constant \added{$C_{\text{BR}}$}\deleted{$C$} depends on the shape regularity of the mesh, polynomial degree $N$, and the diffusion coefficient $\kappa(x)$.

Based on estimator \eqref{eq:etaR}, Melenk and Wohlmuth developed a residual-based error 
estimator for $hp$-FEM in \cite{melenk2001residual} and proved that the 
estimator provides an upper bound \added{on}\deleted{for} the exact error up to a constant \added{$C_{\text{MW}}$},
\added{
\begin{equation}
\label{eq:etaRN}
\eta^2(u_h) = \sum_{K\in\mathcal{T}_h} \frac{h_K^2}{N^2}\| r_E(u_h)\|_{
    K}^2+ \sum_{\ell\in \mathcal{E}} \frac{h_\ell}{N} \left\| r_J(u_h)\right\|_\ell^2.
\end{equation}
}
The constant \added{$C_{\text{MW}}$}\deleted{$C$} shown in the upper bound (similar to \eqref{eq:etaconst}) is independent of $h_K$ and $N$, but depends on the shape regularity of the mesh and the diffusion coefficient $\kappa(x)$. 

In \cite{oden1989toward, bernardi2000adaptive, petzoldt2001regularity}, estimator \eqref{eq:etaR} is extended to an estimator explicitly depending on $\kappa(x)$ for linear FEM,
\added{
\begin{equation}
\label{eq:etaRNkappa}
\eta^2(u_h) = \sum_{K\in\mathcal{T}_h}  \frac{h_K^2}{\kappa_K}\| r_E(u_h)\|_{
    K}^2+  \sum_{\ell\in \mathcal{E} } \frac{h_\ell}{\kappa_{\ell}} \left\| r_J(u_h)\right\|_\ell^2 .
\end{equation}
}
Here $\kappa_K = \max\limits_{x\in K} \kappa(x)$ and 
\deleted{$\kappa_{\ell} =\max\limits_{\ell \in \partial K} \kappa_K$.}
\added{$\kappa_{\ell} =\max\left\{\kappa_K \, | \, K\in \mathcal{T}_h \text{ and }\ell \in \partial K \right\}.$}
Assuming $\kappa(x)$ is quasimono-\\tonically distributed, i.e. $ \kappa(x)$ has at most one local maximum around each node, \eqref{eq:etaRNkappa} is an upper bound \added{on} \deleted{for} the exact \deleted{total}\added{discretization} error up to a constant \added{$C_{\kappa}$} depending only on the shape regularity of the mesh for linear element approximation \cite{petzoldt2001regularity}. 
 If this condition does not hold, the constant \added{$C_{\kappa}$} depends on the bound $\frac{\max_{x\in\Omega} \kappa(x)}{\min_{x\in\Omega} \kappa(x)}$.

We combine the $h$, $N$, and $\kappa(x)$ scaling in \eqref{eq:etaRN} and \eqref{eq:etaRNkappa} to obtain a heuristic indicator with explicit dependence on these parameters as follows,
\added{
\begin{equation}
\label{eq:etaRR}
\eta_{\text{R}}(u_h) = \left(\sum_{K\in\mathcal{T}_h} \eta_{\text{R},K}^2\right)^{1/2} , \quad \eta_{\text{R},K}^2(u_h)=
 \begin{cases} \frac{h_K^2}{\kappa_K N^2}\| r_E(u_h) \|_{K}^2  + \sum\limits_{\ell\in \mathcal{E}_{\text{int}}\cap\partial K} \frac{h_\ell}{2\kappa_{\ell} N} \left\| r_J(u_h)\right\|_{\ell} ^2, \\ 
 \frac{h_K^2}{\kappa_K N^2}\| r_E(u_h) \|_{K}^2  + \sum\limits_{\ell\in \mathcal{E}_{\text{bd}}^N\cap\partial K} \frac{h_\ell}{\kappa_{\ell} N} \left\| r_J(u_h)\right\|_{\ell} ^2.
 \end{cases}
\end{equation}
}
\added{Since $u_h$ is unknown, we use $ \eta_{\text{R}}$  evaluated at $u_h^k$ to estimate $\|e\|_E$. }\deleted{We use $\eta_{\text{alg}}$ and $ \eta_{\text{R}}$ to approximate  $\|e_{\text{alg}}\|_E$ and $\|e\|_E$ in \eqref{ideal sc}, respectively} \added{Note that as the Galerkin orthogonality condition is not satisfied, $\eta_{\text{R}}(u_h^k)$ is not theoretically equivalent to the discretization error or the total error. Further discussion on using $\eta_{\text{R}}(u_h^k)$ to estimate the total error can be found in \cite{arioli2013stopping, papevz2018residual}.}
\deleted{and they are}\added{Applying the estimate $ \eta_{\text{R}}(u_h^k)$ and $\eta_{\text{alg}}(u_h^k)$ to \cref{stopping criterion}, we obtain} \deleted{included in} the following stopping criterion:
\added{
\begin{equation*}
   \eta_{\text{alg}}(u_h^k) \leq \tau \eta_{\text{R}}(u_h^k).
    \label{etaRsc}
\end{equation*}
}
\deleted{\textbf{Modified residual estimate}}

\subsubsection{Flux recovery-based estimator} 
Recovery-based a posteriori error estimators have been  studied extensively, see  \cite{zienkiewicz1987simple, zienkiewicz1992superconvergent,
carstensen2002each, 
ern2015polynomial,
cai2017improved} for examples. In this work, we use an \added{accuracy-enhancing} \deleted{accuracy enhancing} projection to reconstruct the numerical flux and compare it with the original numerical flux $\kappa(x)\nabla u_h^k$.
\added{To keep the implementation simple, we reconstruct the numerical flux using the Brezzi-Douglas-Marini (BDM) space for triangular elements and the Raviart-Thomas space for quadrilateral elements, solving the local problem element-wise by adopting methods from \cite{fortin1991mixed, cai2012robust}. In the following, we demonstrate the flux recovery-based estimator using triangular elements as an example.}
\deleted{In the following, we reconstruct the numerical flux using Brezzi-Douglas-Marini (BDM) elements and solve the local problem element-wise, adopting methods presented in \mbox{\cite{bastian2003superconvergence, cai2012robust}}.} 
For the edge $\ell\in\mathcal{E}$\added{, let $K^+$ and $K^-$ be two elements sharing the edge $\ell$ such that} \deleted{and }$\ell = \partial K^+ \cap \partial K^-$\deleted{,}\added{. For all $w(x)\in L^2(\Omega)$,} we define the weighted average for $w$ on $\ell$ 
\[
\left\{ w \right\}^{\kappa}_{\ell} = \frac{\kappa^{-}}{\kappa^- + \kappa^+} w^+ + 
\frac{\kappa^{+}}{\kappa^- + \kappa^+}   w^-,
\]
where \deleted{$\kappa^{-}$ and $\kappa^{+}$ are the restriction of $\kappa(x)$ on $K^-$ and $K^+$ , respectively.}\added{$\kappa^{-}$ and $w^{-}$, as well as $\kappa^{+}$ and $w^{+}$, are the restrictions of $\kappa(x)$ and $w(x)$ to $\ell$ on 
$K^-$ and $K^+$ , respectively.} Similarly, we denote the weighted jump for $w$ on $\ell$ by 
\[
\left[ w \right]^{\kappa}_{\ell} = 
\frac{\kappa^{-}}{\kappa^- + \kappa^+} \left(w^+ - w^- \right).
\]
The reconstruction is  as \added{follows. Fix an element $K$, then find} \deleted{follows, fix an element $K$, find} $ \boldsymbol{\sigma_K} \in  \left(\mathbb{P}_{N}(K)\right)^2$ satisfying
\begin{equation}
\label{eq:BDM}
\begin{aligned}
 \int_{K} \boldsymbol{\sigma_K} \cdot \nabla w \added{\,dx}& = \int_{K} \kappa(x) \nabla u^k_h \cdot \nabla w \added{\,dx}, \quad \forall w \in \mathbb{P}_{N-1}(K), \\
 \int_{K}\boldsymbol{\sigma_K} \cdot \mathbf{S}(\psi) \added{\,dx}& = \int_{K} \kappa(x) \nabla u^k_h \cdot \mathbf{S}(\psi) \added{\,dx},   \quad \forall \psi \in M_{N+1}(K),\\
 \int_{\ell_i}\left(\boldsymbol{\sigma_K} \cdot \mathbf{n} \right) z \added{\,ds}&  =  \int_{\ell_i}  \left\{ \kappa(x) \nabla u^k_h  \cdot \mathbf{n}_{\ell_{i}} \right\}^{\kappa}_{\ell_i}  z \added{\,ds},
 \quad \forall z \in \mathbb{P}_{N}\left(\ell_i\right),\ell_i \in \partial K, \quad i=1,2,3.
\end{aligned}
\end{equation}
Here, $\mathbf{S}(\psi)=\left(\partial \psi / \partial x_2,-\partial \psi / \partial x_1\right)$. Let $M_N(K)$ be the space of polynomials $\phi \in \mathbb{P}_{N}(K)$ vanishing on the boundary of $K$, 
\[
M_N(K)=\left\{\phi \in\mathbb{P}_{N}(K):\left.\phi\right|_{\partial K}=0\right\}.
\]
Let $\boldsymbol{\rho_K} = \boldsymbol{\sigma_K} - \kappa(x)\nabla u_h^k(x)$. From \eqref{eq:BDM}, \added{$\boldsymbol{\rho_K}$} \deleted{it} satisfies
\begin{equation}
\label{eq:BDMrho}
\begin{aligned}
 \int_{K} \boldsymbol{\rho_K} \cdot \nabla w \added{\,dx}& = 0, \quad \forall w \in \mathbb{P}_{N-1}(K), \\
 \int_{K}\boldsymbol{\rho_K} \cdot \mathbf{S}(\psi) \added{\,dx}& = 0,   \quad \forall \psi \in M_{N+1}(K),\\
 \int_{\ell_i}\left(\boldsymbol{\rho_K} \cdot \mathbf{n} \right) z\added{\,ds} &  =  \int_{\ell_i} \left[ \kappa(x) \nabla u^k_h \cdot \mathbf{n}_{\ell_{i}} \right]_{\ell_i}^{\kappa}  z \added{\,ds},
 \quad \forall z \in \mathbb{P}_{N}\left(\ell_i\right),\ell_i \in \partial K, \quad i=1,2,3.
\end{aligned}
\end{equation}

Now we define the error estimator based on recovery, 
\begin{equation}
    \eta_{\text{FR}}(u_h^k) := \left(\sum_{K\in\mathcal{T}_h}\left\| {\kappa(x)}^{-1/2}{\boldsymbol{\rho_K}} \right\|^2\right)^{1/2}.
    \label{eta_bdm}
\end{equation}
By solving equation \eqref{eq:BDMrho}, the jump of the normal component of the numerical flux on the edge $\ell$, $\left[ \kappa(x) \nabla u^k_h \cdot \mathbf{n}_{\ell_{i}} \right]_{\ell_i}^{\kappa}$, is lifted to the elementwise function, $\boldsymbol{\rho_K}$. We represent the discretized form of the lifting operator by 
\begin{equation}
    \By_K = \BL_K  \Bd_K,
\label{lifting}
\end{equation}
where $\By_K$ corresponds to $\boldsymbol{\rho_K}$, $\Bd_K$ is associated with $\left[ \kappa(x) \nabla u^k_h \cdot \mathbf{n}_{\ell_{i}} \right]_{\ell_i}^{\kappa}$, and $\BL_K$ is the lifting matrix. 
In each step of the iteration, 
$\boldsymbol{\rho_K}$ is obtained by matrix-vector multiplication \eqref{lifting} for all elements. 
However, storing the matrices $\left\{\BL_K\right\}$ for all elements can take a significant amount of memory. This demands allocation for $\mathcal{O}(N^{2d-1})$ entries for each element. In double precision and for $N=8$, it requires 9.5KB per triangle or 232KB per tetrahedron to store the lifting matrix\deleted{.}\added{, which becomes impractical for large-scale computations on GPUs.}
\added{
Applying $\eta_{\text{FR}}(u_h^k)$ and $\eta_{\text{alg}}(u_h^k)$ to \cref{stopping criterion}, we derive
the following stopping criterion 
\begin{equation*}
\eta_{\text{alg}}(u_h^k) \leq \tau \eta_{\text{FR}}(u_h^k)
    \label{etaBDMsc}
\end{equation*}
}
\deleted{
To reduce the computational cost and storage, we instead introduce a lower bound \added{on} \deleted{for} $\eta_{\text{FR}}$.
For the element \mbox{$K\in\mathcal{T}_h$},}
\deleted{where \mbox{$\mu_K^2$} is the smallest eigenvalue of \mbox{$\BL_K^{ T} \BM_K \BL_K$}, 
and \mbox{$\BM_K$} is the mass matrix with elements defined by \mbox{$\left(\BM_K\right)_{ij} = \int_K \kappa(x)^{-1} \phi_j \phi_i dx$.}
We define the new error indicator,}
\deleted{
The error indicator \mbox{$\eta_{\text{\underline{BDM}}}$} is a lower bound \added{on} \deleted{for} \mbox{$\eta_{\text{FR}}$}. Before the first iteration, the linear systems \mbox{\eqref{eq:BDMrho}} are constructed, the eigenvalues \mbox{$\mu_K$} are computed and stored, and matrices are discarded.
Subsequently, only \mbox{$\mu_K$} are used to compute \mbox{\eqref{eq:BDMLB}} in the iterations.
Thus, it is less expensive to compute the error indicator \eqref{eq:BDMLB} than \eqref{eta_bdm}.
The stopping criteria for \mbox{$\eta_{\text{{BDM}}}$} and \mbox{$\eta_{\text{\underline{BDM}}}$} are }
\deleted{
However, the second criterion is stricter than the first one due to \mbox{$\eta_{\text{\underline{BDM}}}$} being a lower bound \added{on} \deleted{for} \mbox{$\eta_{\text{FR}}$}.}

%% file: 3-etaRF.tex
\section{Stopping Criteria derived from the residual}
\label{rnfn}

In this section, we describe a stopping criterion derived directly from the \deleted{linear residual} \added{linear system residual} and generalize it to the Poisson problem 
with highly variable piecewise constant coefficient.

\subsection{Globally constant coefficient}

The $n$-th component of the \deleted{linear residual}\added{linear system residual} $\Br_k = \Bb - \BA \Bx_k$ is
\[
\begin{aligned}
\left(\Br_k\right)_n & = \Bb_n - (\BA \Bx_k)_n \\
& = \left(\phi_n, f \right) - \left(\phi_n, g \right)_{\partial \Omega} -  \sum_{K\in \mathcal{T}_h}  \left(\kappa(x)\nabla \phi_n, \nabla u^k_h\right)_K.
\end{aligned}
\]
Integrating the last term by parts, we obtain
\[
\begin{aligned}
    \left(\Br_k\right)_n = & \sum_{K\in \mathcal{T}_h} \left(\phi_n, r_E(u_h^k)\right)_K  - \sum_{\ell\in \mathcal{E}} \left( \phi_n, r_J(u_h^k)\right)_{\ell} \\
    = & \left(\BR_k\right)_n  + \left(\BF_k\right)_n,
\end{aligned}
\]
where $\BR_k, \BF_k\in \mathbb{R}^{N_s}$, $\left(\BR_k\right)_n = \sum_{K\in \mathcal{T}_h}\left(\phi_n, r_E(u_h^k)\right)_K$,  and $\left(\BF_k\right)_n = -\sum_{\ell\in \mathcal{E}}\left( \phi_n, r_J(u_h^k)\right)_{\ell}$.
We introduce the indicator $\eta_{\text{RF}}$ 
\begin{equation}
    \eta_{\text{RF}}(u_h^k) := \|\BR_k\|  + \|\BF_k\|,
    \label{eq:RF}
\end{equation}
with the associated stopping criterion:
\begin{equation}
    \|\Br_k\| \leq \tau \, \eta_{\text{RF}}(u_h^k).
    \label{eq:criterion}
\end{equation}
The indicator $\eta_{\text{RF}}$ is an upper bound \added{on} \deleted{for} the norm of the residual without any unknown constants involved. Ideally, $\eta_{\text{RF}}$ should closely track $ \|\Br_k\|$ until the total error converges, and the separation between  $ \|\Br_k\|$ and $\eta_{\text{RF}}$ should indicate the deviation of  the total error from the algebraic error. 
Furthermore, it is only necessary to compute $\BR_k$; $\BF_k = \Br_k - \BR_k$ can be directly calculated once $\BR_k$ has been determined.
Additionally,  on a uniform mesh, the $h$ scaling of $\eta_{\text{RF}}$ is\deleted{$h^{d-2}$, which is} consistent with the scaling of the norm of the residual\deleted{and $\eta_{\text{MR}}$}. However, it differs from the $h$ scaling of the total error except \added{when $\Omega$ is a two-dimensional domain.}\deleted{$d=2$}. 

\subsection{Highly variable piecewise constant coefficient}

We observe that there exists  a $\kappa(x)$ scaling difference between $\eta_{\text{RF}}$ and the total error.
When $\kappa(x)$ is highly variable, the difference in $\kappa(x)$ scaling may impact the effectivity of the stopping criterion \eqref{eq:criterion}. 
If a good preconditioner is available, leading to a rapid decrease in the CG error, we employ a weighted $l^2$ norm in \eqref{eq:criterion} as an alternative to the standard $l^2$ norm. This approach ensures that $\eta_{\text{RF}}$ shares the same $\kappa(x)$ scaling as the total error.
We define a weight vector  $\Bw \in \mathbb{R}^{N_s}$ \added{with its 
$n$-th component given by}
\[
\Bw_n = \min\limits_{x\in \omega_n} \kappa(x)^{-1},
\]
where $\omega_n = \operatorname{supp}\left(\phi_n\right)$. The weight $\Bw$ is similar to the $\kappa(x)$ scaling in \eqref{eq:etaRR}.
Let 
$$\|\By\|_{\Bw} = \left(\By^T \operatorname{diag}(\Bw)\By\right)^{1/2}$$ 
be the weighted $l^2$ norm of $\By$, where $\operatorname{diag}(\Bw)$ is a diagonal matrix with diagonal elements given by  $\Bw$.
We define the indicator $\eta_{\text{RF}}^{\Bw}$ 
\begin{equation}
    \eta_{\text{RF}}^{\Bw}(u_h^k) := \| \BR_k\|_{\Bw}  + \| \BF_k\|_{\Bw},
    \label{eq:heter RF}
\end{equation}
and we substitute \eqref{eq:criterion} with: 
\begin{equation}
    \| \Br_k\|_{\Bw} \leq \tau \, \eta_{\text{RF}}^{\Bw}(u_h^k).
    \label{eq:heter criterion}
\end{equation}

In cases where the good preconditioner is not available and  CG converges slowly, 
we partition the domain $\Omega$ into several subdomains $\Omega_p,$ $p=1,\ldots, P$ based on the value of $\kappa(x)$, and compare $ \|\Br_k\|_{\Bw}$  with $\eta_{\text{RF}}^{\Bw}$ restricted to these subdomains. The details of the partition are shown in \cref{subsubsetion:4-3-1}, and the computations of subdomain indicators are presented in \cref{appen}.
We propose a subdomain-based stopping criterion
as an alternative to the criterion \cref{eq:heter criterion} 
\begin{equation}
    \|\Br_k^p\|_{\Bw} \leq \tau \, {\eta^{\Bw,p}_{\text{RF}}}(u_h^k), \quad \forall p=1,\ldots, P.
    \label{eq:patchy}
\end{equation}
Here $\eta_{\text{RF}}^{\Bw,p} = \|\BR_k^p\|_{\Bw} + \|\BF_k^p\|_{\Bw}$, and $ \BR_k^p$, $\BF_k^p$, $\Br_k^p$ are vectors obtained by restricting $ \BR_k$, $\BF_k$, $\Br_k$ to  subdomain $\Omega_p$.

In contrast to criterion \cref{eq:heter criterion}, the subdomain-based criterion \cref{eq:patchy} leads to termination when the weighted norm of the local \deleted{linear residual} \added{linear system residual} is dominated by the local indicator  in all subdomains. Although the global total error may have converged, additional iterations may be required for the local errors to converge. As a result, the subdomain-based criterion \cref{eq:patchy} may recommend stopping the iteration slightly later than the criterion \cref{eq:heter criterion} suggests. 

%% file: 3-3.tex
\subsection{\added{Reliability and efficiency of the criterion}}
\label{sec:proof}
\added{
In this subsection, we prove the reliability and efficiency of the criterion \cref{eq:criterion}. 
We establish   the reliability theorem  to ensure that the criterion does not terminate the iteration prematurely, ensuring that once the stopping condition is met, the optimal stopping condition \cref{ideal sc} is also satisfied.
We further establish   the efficiency theorem  to ensure that the criterion avoids unnecessary iterations. Once the optimal stopping condition \cref{ideal sc} is achieved, the criterion also indicates that the iteration can stop.
}

\added{
All supporting lemmas and inequalities labelled (A, $*$) can be found in \cref{sec:6-1}.
Our analysis is based on the following assumption.
\begin{assumption}
    The triangulation $\mathcal{T}_h$ is quasi-uniform \cref{eq:quasi uniform} with quadrilateral elements. $\kappa(x)=1$.
    \label{assump}
\end{assumption}
First, we state the reliability theorem. 
\begin{theorem}
    Let $\Omega\subset \mathbb{R}^2$ be a bounded domain and let $0<\varepsilon, 0<\tau<1/2$. 
    Suppose $u\in H_{0, \Gamma_D}^1(\Omega)$ is the solution to the weak form given in \cref{eq:weak}, and $f\in L^2(\Omega)$ is the forcing function. Let  $u_h\in \Vhn$ satisfy the finite element approximation in \cref{eq:fem}. Furthermore, for all $w_h \in \Vhn$,
    the indicator $\eta_{\text{RF}}(w_h)$ is defined in \cref{eq:RF} and 
    and the residual $\Br(w_h) \in \mathbb{R}^{N_s}$ has its $n$-th component given by $\Br_n(w_h) = a(w_h, \phi_n) - l(\phi_n)$.
    Under \cref{assump}, for any $w_h \in \Vhn$, 
    if the following stopping criterion \cref{eq:criterion} is satisfied:
    \[
    \|\Br(w_h)\|_{l^2} \leq  \tau  \eta_{\text{RF}}(w_h),
    \]
    then there exists a constant $C(\varepsilon)>0$ depending on $\varepsilon$, but independent of the polynomial degree $N$ and mesh size $h$, such that
    \begin{equation}
    \|\nabla(u_h - w_h)\|  \leq  \tau C(\varepsilon) \frac{N^{2+\varepsilon}}{h}  
 \left( \|\nabla(u - w_h)\|  + \sum\limits_{K\in\mathcal{T}_h} \frac{h_K}{N^{3/2-3\varepsilon}} \|f_{h,K} - f\|_K \right),
 \label{eq:etaRFtoerror}
 \end{equation}
 where $f_{h,K}$ is the $L^2(K)$ projection of $f$ on the polynomial space of degree $N-1$.
\label{thm:etaRFtoerror}
\end{theorem}
\begin{proof}
   Combining \cref{etaRFetaR2,res and error,etaRerror} with the stopping criterion \cref{stopping criterion}, we can derive
\[
\begin{aligned}
    \lammina \|\nabla(u_h - w_h)\|^2  \leq  \|\Br(w_h)\|_{l^2}^2 
     \leq \tau^2   \eta_{\text{RF}}^2(w_h)
     \leq C\tau^2\lammaxm \frac{N^2}{h^2} \eta_R^2(w_h)\\
    \leq \tau^2 \lammaxm \frac{N^2}{h^2} C(\varepsilon)\left\{ N^{2+2\varepsilon}  \|\nabla(u - w_h)\|^2 + \sum\limits_{K\in\mathcal{T}_h} \frac{h_K^2}{N^{1-4\varepsilon}} \|f_{h,K} - f\|_K^2 \right\}, 
\end{aligned}
\]
where $C, C(\varepsilon)>0$ are independent of $N$ and $h$. Thus,
\[
\begin{aligned}
    \|\nabla(u_h - w_h)\| & \leq \tau \sqrt{C(\varepsilon)}  \sqrt{\frac{\lammaxm}{\lammina }} \frac{N}{h}  \left\{ N^{2+2\varepsilon}  \|\nabla(u - w_h)\|^2  + \sum\limits_{K\in\mathcal{T}_h}\frac{h_K^2}{N^{1-4\varepsilon}} \|f_{h,K} - f\|_K^2 \right\}^{1/2} \\
 & \leq  \tau\sqrt{C(\varepsilon)} \frac{\overC_{\BM}}{\underC_{\BA}} \frac{N^{2+\varepsilon}}{h}  
 \left( \|\nabla(u - w_h)\|  + \sum\limits_{K\in\mathcal{T}_h}\frac{h_K}{N^{3/2-3\varepsilon}} \|f_{h,K} - f\|_{K} \right),   
\end{aligned}
\]   
where the last inequality is from \cref{lem:lambda}.
\end{proof}
The previous theorem ensures that applying \cref{stopping criterion} results in an algebraic error bounded by the total error multiplied by a constant that depends on $\tau$. Therefore, by appropriately selecting the parameter $\tau$, we can achieve the desired accuracy of the algebraic error. 
Conversely, to demonstrate that the stopping criterion \cref{stopping criterion} is met given a small algebraic error relative to the total error, we present the efficiency theorem.
\begin{theorem}
Let $0<\nu<1$.
With the same notations and assumption as in \cref{thm:etaRFtoerror},
if the following condition holds: 
\[
\|\nabla(u_h - w_h)\| \leq \nu  \|\nabla(u - w_h)\|,
\]
then there exists a constant $C>0$, independent of $N$ and $h$, such that
\begin{equation}
\|\Br(w_h)\|_{l^2} \leq C \frac{\nu}{\sqrt{1-\nu^2}}  \left( {N}^{3/2} \eta_{\text{RF}}(u_h) + \sum\limits_{K\in\mathcal{T}_h} \frac{h_K}{\sqrt{N}} \|f_{h,K} - f\|_K\right).
\label{eq:errortoetaRF}
\end{equation}
\label{thm:errortoetaRF}
\end{theorem}
\begin{proof}
According to \cref{eq:galerkin orth},
the desirable stopping condition \cref{ideal sc}
 is equivalent to 
\[
 \|\nabla(u_h - w_h)\|^2 \leq  \frac{\nu^2}{1-\nu^2}  \|\nabla(u - u_h)\|^2.
\]
Combining \cref{lem:lambda}, \cref{etaRbound}, \cref{etaRFetaR2} with the equivalent desirable stopping condition above, 
we establish the following result, 
\[
\begin{aligned}
    \|\Br(w_h)\|_{l^2}^2 & \leq \lammaxa \|\nabla(u_h - w_h)\|^2  
     \leq   \lammaxa  \frac{\nu^2}{1-\nu^2} \|\nabla(u - u_h)\|^2  \\
    & \leq \frac{\lammaxa}{C_1}  \frac{\nu^2}{1-\nu^2} \left(\eta_{R}^2(u_h) + \sum\limits_{K\in\mathcal{T}_h}\frac{h_K^2}{N^2}\|f_{h,N} - f\|_K^2\right)  \\
    & \leq \frac{\lammaxa}{C_1}  \frac{\nu^2}{1-\nu^2} \left(\frac{1}{\lamminm} \frac{C_2 h^2}{N^2}\eta_{\text{RF}}^2(u_h) + \sum\limits_{K\in\mathcal{T}_h}\frac{h_K^2}{N^2}\|f_{h,N} - f\|_K^2\right) \\
    & =  \frac{\nu^2}{C_1(1-\nu^2)} \left( \frac{\lammaxa}{\lamminm} \frac{C_2 h^2}{N^2} \eta_{\text{RF}}^2(u_h) + \lammaxa \sum\limits_{K\in\mathcal{T}_h}\frac{h_K^2}{N^2}\|f_{h,N} - f\|_K^2\right)   \\
    & \leq  \frac{\nu^2}{C(1-\nu^2)} \left( N^3 \eta_{\text{RF}}^2(u_h) +  \sum\limits_{K\in\mathcal{T}_h}\frac{h_K^2}{N}\|f_{h,N} - f\|_K^2\right).  
\end{aligned}
\]
\end{proof}
Suppose $w_h = u_h^k$ is obtained from CG iteration. The indicator $\eta_{\text{RF}}(u_h^k)$ converges to $\eta_{\text{RF}}(u_h)$ as $k$ increases. Thus, there exists an integer $k_0>0$, such that for all $k\geq k_0$, $\eta_{\text{RF}}(u_h) \leq 2 \eta_{\text{RF}}(u_h^k)$. Using \cref{thm:errortoetaRF}, we can bound $\|\Br(u_h^k)\|$ by $\eta_{\text{RF}}(u_h^k)$ for $k\geq k_0$.
}

\added{
Note that analysis in \cref{thm:etaRFtoerror} and \cref{thm:errortoetaRF} is not sharp. In next section, numerical experiments demonstrate that the performance of the stopping criterion is independent of $N$ and $h$.
}

%% file: 4-NumericalResults.tex
\section{Numerical experiments}
\label{num}
In this section, we test the effectiveness and robustness of stopping criteria with respect to the polynomial degree $N$, the shape regularity of the mesh, the diffusion coefficient $\kappa(x)$, and the singularity of the solution.
We consider four examples. 
In \cref{4-1}, we apply criteria to the Poisson problem with a constant diffusion coefficient and \deleted{a geometrically regular mesh}\added{shape regular quadrilateral mesh \cite{braess2001finite}}, demonstrating the validity of the criteria.
\deleted{In \mbox{\cref{4-2}}, we use a diamond-shaped mesh, an anisotropic mesh presented in \mbox{\cite{babuvska1992quality}}, to demonstrate the sensitivity of the criteria to the shape regularity of the mesh.}
In \cref{4-3}, we demonstrate the performance of criteria for problems with highly variable piecewise constant coefficients and 
\added{shape regular triangle mesh. Its solution has}\deleted{a solution with} singularities caused by jumps in the  coefficient and the reentrant corner of the L-shape domain, similar to those presented in \cite[section 4.1]{arioli2004stopping}, \cite[example 7.5]{carstensen2010estimator}, and \cite[section 7.6]{papevz2018estimating}.
In \cref{4-4}, we consider the same problem as described in \cref{4-3-2} and solve the linear system using the preconditioned recycling CG, to show that deflation using the recycle subspace is beneficial in achieving efficient termination of the iteration process.
\added{In \cref{4-5}, we implement the stopping criterion in a GPU-accelerated PDE solver and collect runtime of applying various criteria to show the effectiveness of the proposed criterion.}
We summarize the tested stopping criteria as follows:
\begin{enumerate}[start=1,label={\bfseries (C\arabic*)}]
    \item $\eta_{\text{alg}} \leq \tau \eta_{\text{R}}$, where $\eta_{\text{R}}$  \eqref{eq:etaRR} is the most commonly used a posteriori error estimator; \label{c1}
    \item $\eta_{\text{alg}} \leq \tau \eta_{\text{FC}}$, where $\eta_{\text{FC}}$ \eqref{eta_bdm} is the error estimator based on flux reconstruction; \label{c3}
    \item $\| \Br_k\|_{\Bw} \leq \tau \eta_{\text{RF}}^{\Bw}$, where $\eta_{\text{RF}}^{\Bw}$ \eqref{eq:heter RF} is an upper bound for $\| \Br_k\|_{\Bw}$; \label{c5}
    \item $\|\Br_k^p\|_{\Bw} \leq \tau \eta_{\text{RF}}^{\Bw,p}$, for all $p = 1, \cdots, P$, where $\eta_{\text{RF}}^{\Bw,p}$ is the subdomain indicator; \label{c6}
    \item $\|\Br_k\| \leq \text{TOL} \|\Br_0\|$, where TOL is a preset relative tolerance. \label{c7}
\end{enumerate}

In criteria \added{\ref{c1}, \ref{c3}} \deleted{\ref{c1}-\ref{c3}}, we compare 
a posteriori error estimates  $\eta_{\text{R}}$\deleted{, $\eta_{\text{MR}}$} and $\eta_{\text{FC}}$ to the estimate of the algebraic error  $\eta_{\text{alg}}$ \eqref{eq:alg err}. 
\deleted{Moreover, although $\eta_{\text{\underline{BDM}}}$ is a lower bound for the estimate $\eta_{{{\text{BDM}}}}$, we also compare it to $\eta_{\text{alg}}$ in \ref{c4}. }
Conversely, in criterion \ref{c5}, the error indicator $\eta_{\text{RF}}^{\Bw}$, derived from the \deleted{linear residual}\added{linear system residual}, shares greater similarity with the weighted norm of the \deleted{linear residual}\added{linear system residual}, and as a result, it is comparable to the weighted norm of the residual rather than the estimate of the algebraic error.
Criterion \ref{c6} is the subdomain-based criterion for problems with highly variable piecewise constant coefficients. 
Lastly, criterion \ref{c7} is an often used criterion  based on the relative residual norm. 

We define the quality ratio of a criterion as 
\begin{equation}
    \text{quality ratio}: = \frac{\|u-u^{k^*}_h\|_E}{\|u-u_h\|_E},
    \label{ev}
\end{equation}
where $u^{k^*}_h$  is the first solution that satisfies the stopping condition during the iterative process. We note that the quality ratio is always greater than one. 
If the quality ratio is much greater than one, it implies a premature termination.
It is important to note that the quality ratio, which measures the reliability of a stopping criterion, should not be confused with the effectivity index, a common term used in many a posteriori error estimate papers, which indicates the efficiency of an error estimator.

In the following subsections, experiments are performed in Matlab R2019b \added{and \texttt{libParanumal}, a collection of GPU-accelerated flow solvers \cite{ChalmersKarakusAustinSwirydowiczWarburton2020}}. When no additional details are provided, we apply the preconditioned conjugate gradient algorithm in 
\cite{arioli2004stopping} with a zero initial guess to solve the linear systems. We use the incomplete Cholesky decomposition preconditioner with  empirically selected drop tolerance of $10^{-4}$ and  diagonal shift of $0.1$.
We choose the delay parameter $d=10$ in the algebraic error estimator \eqref{eq:alg err}. 
In \deleted {the subsequent examples}\added{\cref{4-1}, \cref{4-3}, and \cref{4-4}}, we compare the approximate solution from CG to the linear system \eqref{eq:linear eq} to the solution obtained using the backslash command in MATLAB. 
In the tables presented below, we collect the numbers of iterations and quality ratios when applying stopping criteria. For criteria relying on $\eta_{\text{alg}}$, the iterations attributed to the delay in the computation of $\eta_{\text{alg}}$ are not included in the iteration count. However, in practice, all criteria that depend on $\eta_{\text{alg}}$ require $d$ additional iterations.
In the following figures, all error estimators and indicators are denoted by markers, while all exact errors and the norm of the \deleted{linear residual}\added{linear system residual} are represented without markers.

\input{4-1.tex}

{\deleted{\textbf{4.2 Test problem 2: highly anisotropic mesh}}}

\input{4-3.tex}

\input{4-4.tex}

\input{4-5.tex}

\input{4-6ResultsSummary}

%% file: 4-1.tex
\subsection{{Test problem 1: i}sotropic mesh}
\label{4-1}

We consider the Poisson problem \eqref{eq:poisson} on $\Omega = [0,1]^2$ with the homogeneous \deleted{Neumann}\added{Dirichlet} boundary condition, $\kappa(x)=1$, and choose the right-hand side function $f$ such that the solution to the continuous problem is given by
$$u(x, y) = (1-x^2)^2(1-y^2)^2 e^{x+y}.$$ 
We  discretize the problem on a \deleted{shape regular} mesh with \deleted{128 isosceles right triangle} \added{144 quadrilateral} elements,  using piecewise polynomials with degree $N=4,6,8$. Since $\kappa(x)=1$, the weight vector $\Bw_n = 1,  n=1,\ldots,N_s$. The weighted $l^2$ norm is the same as the  $l^2$ norm, and   $\eta_{\text{RF}}$ is the same as $\eta_{\text{RF}}^{\Bw}$.

\cref{3-1} shows the energy norm of the error and the error estimates in the iteration process with $N=6$.
We observe that $\eta_{\text{alg}}$ tracks the $\BA$-norm error accurately as CG converges fast.  
Indicators $\eta_{\text{R}}$\deleted{, $\eta_{\text{MR}}$,} and $\eta_{\text{RF}}$  slightly overestimate the total error by a factor less than 10.
The estimator \added{$\eta_{\text{FC}}$}\deleted{$\eta_{\text{BDM}}$} provides a very tight estimate for the total error. \deleted{As a lower bound for \added{$\eta_{\text{FC}}$}\deleted{$\eta_{\text{BDM}}$}, the indicator $\eta_{\text{\underline{BDM}}}$ gives a lower bound for the total error (although not always). }

\begin{figure}[h!]
\begin{tikzpicture}
\begin{groupplot}[group style={
                      group name=myplot,
                      group size= 2 by 1,
                       horizontal sep = 50pt,
                        vertical sep= 2cm
                    },height=6cm,width=7.1cm]
 \nextgroupplot[  
    xlabel={iteration},
    ymode=log,  
    ymin=1e-9,
    ymax=1e1 ]
 \input{Data/test1rf}   
\nextgroupplot[ 
    xlabel={iteration},
    ymode=log,
    ymin=1e-9,
    ymax=1e1,
]
 \input{Data/test1eta}   
    \end{groupplot}
\path (myplot c1r1.south west|-current bounding box.south)--
      coordinate(legendpos)
      (myplot c2r1.south east|-current bounding box.south);
\matrix[
    matrix of nodes,
    anchor=south,
    draw,
    inner sep=0.2em,
    draw
  ]at([yshift=-9ex]legendpos)
  {
    \ref{plot:totalerr}& total error&[5pt]
    \ref{plot:exactalg}& $\|\Bx_k - \Bx\|_{\BA}$ &[5pt]
    \ref{plot:etarf}& $\eta_{\text{RF}}$ & [5pt]
    \ref{plot:res}& $ \| \Br_k\|$  & [5pt]\\
    \ref{plot:etaalg}& $\eta_{\text{alg}}$ & [5pt] 
    \ref{plot:etar}& $\eta_{\text{R}}$ & [5pt]
    \ref{plot:bdm}& \added{$\eta_{\text{FC}}$}\deleted{$\eta_{\text{BDM}}$} & [5pt] \\ }; 
\end{tikzpicture}
\caption{Convergence history for test problem 1 (isotropic mesh) with $N=6$.  Left: the total error, the $\BA$\text{-}norm error $\|\Bx_k-\Bx\|_{\BA}$, the norm of the \deleted{linear residual} \added{linear system residual} $\|\Br_k\|$ and $\eta_{\text{RF}}$. Right: the total error, the $\BA$\text{-}norm error $\|\Bx_k-\Bx\|_{\BA}$ and its estimator $\eta_{\text{alg}}$ (delay parameter $d=10$), and the error indicators $\eta_{\text{R}} {, \eta_{\text{MR}},}  {\text{ and }} \eta_{\text{FC}} {\text{ and }\eta_{\text{\underline{BDM}}}}$. }
\label{3-1}
\end{figure}

In criteria $\textbf{(C1)}$-$\textbf{(C5)}$, the parameter $\tau$ plays a crucial role in determining when to stop the iteration. A small $\tau$ may result in early termination, while a large $\tau$ could cause unnecessary iterations. To select a reasonable $\tau$, we plot the quality ratio of stopping criteria in \cref{3-1-tau}, varying $\tau$ from $1/30$ to $1/3$. 
We find that $\tau=1/20$ is an appropriate choice, as all quality ratios remain below  1.1. In subsequent examples, we set $\tau=1/20$.

\begin{figure}[htbp]
\begin{center}
\begin{tikzpicture}
\begin{groupplot}[group style={
                      group name=myplot,
                      group size= 1 by 1,
                       horizontal sep = 50pt,
                        vertical sep= 2cm
                    },height=5cm,width=10cm]
 \nextgroupplot[  
    xlabel={$1/\tau$},
    ylabel={quality ratio},
    ymin=1,
    ymax=1.5 ]
 \input{Data/test1tau}   
    \end{groupplot}
\path (myplot c1r1.north east) -- 
      coordinate(legendpos)
      (myplot c1r1.south east); 
\matrix[
    matrix of nodes,
    anchor=north east,
    fill=white,draw,
    inner sep=0.2em,
=    column 1/.style={nodes={align=center}},
    column 2/.style={nodes={anchor=base west}},
    shift={(-2ex,10ex)} 
  ]
  at (legendpos)
  {
    \ref{plot:etar} &$\textbf{(C1)} \eta_{\text{R}}$ \\
    \ref{plot:bdm} & $\textbf{(C2)} \eta_{\text{FC}}$ \\
    \ref{plot:etarf} & $\textbf{(C3)} \eta_{\text{RF}}$ \\
};
\end{tikzpicture}
\caption{Sensitivity of the stopping criteria quality ratios with respect to $\tau$ for test problem 1. }
\label{3-1-tau}
\end{center}
\end{figure}

In \cref{T1},  we present the number of iterations and the corresponding 
quality ratios \cref{ev}  for  $N=4,6,8$. 
We note first that \deleted{the performances of}  $\eta_{\text{R}}$ \deleted{, $\eta_{\text{MR}}$,} and $\eta_{\text{RF}}$ \deleted{are very similar due to the \added{shape-}regularity of the mesh. 
They all} achieve approximately the same level of accuracy with roughly the same number of iterations. 
\deleted{Criteria  $\eta_{\text{BDM}}$ and $\eta_{\text{\underline{BDM}}}$ also provide} \added{Criterion \added{$\eta_{\text{FC}}$}\deleted{$\eta_{\text{BDM}}$} also provides} a favorable termination. 
For the empirical criterion based on the relative residual norm, 
\added{almost more than half of the iterations are unnecessary.}
\deleted{more than 10  additional iterations are 
required when $N=4$. 
Moreover, the criterion that requires relative residual norm to be less than $10^{-6}$ results in premature termination for $N=8$,  as the total error has not yet sufficiently converged at the stopping point. }
Overall, the first \deleted{five}\added{three} criteria provide reliable and efficient alternatives for this problem. 

\begin{table}[h!]
\centering
\caption{Numbers of iterations (iter) and quality ratios (qual. \eqref{ev}) resulting from applying stopping criteria to the solution of test problem 1. }
\begin{tabular}{ l|c|c|c|c|c|c} 
\hline
\multirow{ 2}{*}{Criterion} & \multicolumn{2}{c|}{$N=4$}  & \multicolumn{2}{c|}{$N=6$} & \multicolumn{2}{c}{$N=8$}  \\ \cline{2-7}
\multirow{ 2}{*}{} & iter & qual.  & iter & qual.  & iter & qual.   \\ \hline 
$\eta_{\text{alg}} \leq \tau \eta_{\text{R}}$ & 21 & 1.00 &35 & 1.02 &52 & 1.02 \\  
$\eta_{\text{alg}} \leq \tau \eta_{\text{FC}}$ &24 & 1.00 &41 & 1.00 &59 & 1.00 \\  
$\|\Br_k\| \leq \tau \eta_{\text{RF}}$ &19 & 1.08 &33 & 1.07 &49 & 1.08 \\  
$\|\Br_k\|\leq 10^{-8}\|\Br_0\|$ &43 & 1.00 &67 & 1.00 &93 & 1.00 \\  
\hline
\end{tabular}
\label{T1}
\end{table}

%% file: Data/test1rf.tex
\addplot[ 
line width=2.5pt, color=gray, solid ]
coordinates { 
(1, 2.39942)
(2, 2.13278)
(3, 1.88453)
(4, 1.67555)
(5, 1.42952)
(6, 1.25535)
(7, 1.10014)
(8, 0.966092)
(9, 0.87457)
(10, 0.764609)
(11, 0.665718)
(12, 0.56734)
(13, 0.461182)
(14, 0.398422)
(15, 0.337207)
(16, 0.281733)
(17, 0.231257)
(18, 0.168413)
(19, 0.134542)
(20, 0.103266)
(21, 0.0735124)
(22, 0.0585249)
(23, 0.0405522)
(24, 0.0282396)
(25, 0.0208284)
(26, 0.013798)
(27, 0.010938)
(28, 0.00913238)
(29, 0.00783348)
(30, 0.00656679)
(31, 0.00516038)
(32, 0.00453362)
(33, 0.00431525)
(34, 0.00421971)
(35, 0.0041132)
(36, 0.00406009)
(37, 0.00405351)
(38, 0.00403467)
(39, 0.00402552)
(40, 0.00403026)
(41, 0.00402091)
(42, 0.00402654)
(43, 0.00402264)
(44, 0.00402384)
(45, 0.00402408)
(46, 0.0040232)
(47, 0.00402389)
(48, 0.00402348)
(49, 0.00402361)
(50, 0.00402363)
(51, 0.00402358)
(52, 0.00402362)
(53, 0.00402359)
(54, 0.00402361)
(55, 0.0040236)
(56, 0.0040236)
(57, 0.0040236)
(58, 0.00402361)
(59, 0.0040236)
(60, 0.00402361)
(61, 0.0040236)
(62, 0.0040236)
(63, 0.0040236)
(64, 0.0040236)
(65, 0.0040236)
(66, 0.0040236)
(67, 0.0040236)
}; \label{plot:totalerr}
\addplot[ 
line width=2.5pt, color=mygreen, dashed]
coordinates { 
(1, 2.39938)
(2, 2.13281)
(3, 1.88452)
(4, 1.67554)
(5, 1.42953)
(6, 1.25534)
(7, 1.10012)
(8, 0.966099)
(9, 0.874556)
(10, 0.764575)
(11, 0.665746)
(12, 0.567247)
(13, 0.461247)
(14, 0.398369)
(15, 0.337147)
(16, 0.281793)
(17, 0.231108)
(18, 0.168435)
(19, 0.134499)
(20, 0.103058)
(21, 0.0736042)
(22, 0.0582677)
(23, 0.0403402)
(24, 0.0280843)
(25, 0.020238)
(26, 0.0134155)
(27, 0.0100535)
(28, 0.00826978)
(29, 0.00669102)
(30, 0.00516815)
(31, 0.00333002)
(32, 0.00201605)
(33, 0.00160953)
(34, 0.00125477)
(35, 0.000809996)
(36, 0.000599238)
(37, 0.000449127)
(38, 0.00030943)
(39, 0.000201169)
(40, 0.000134574)
(41, 9.01569e-05)
(42, 5.62338e-05)
(43, 3.64938e-05)
(44, 2.56074e-05)
(45, 1.66079e-05)
(46, 9.93169e-06)
(47, 6.30005e-06)
(48, 3.5121e-06)
(49, 2.19678e-06)
(50, 1.46223e-06)
(51, 1.04715e-06)
(52, 8.89395e-07)
(53, 7.82104e-07)
(54, 7.12311e-07)
(55, 5.61168e-07)
(56, 4.34979e-07)
(57, 2.9276e-07)
(58, 2.15875e-07)
(59, 1.51793e-07)
(60, 1.07064e-07)
(61, 6.83493e-08)
(62, 4.50855e-08)
(63, 2.92928e-08)
(64, 1.91922e-08)
(65, 1.12589e-08)
(66, 7.56354e-09)
(67, 5.54618e-09)
}; \label{plot:exactalg}
\addplot[ 
line width=1.5pt, color= mycyan, densely dotted]
coordinates { 
(1, 0.681658)
(2, 0.745982)
(3, 0.674745)
(4, 0.888574)
(5, 0.668977)
(6, 0.472025)
(7, 0.537239)
(8, 0.399458)
(9, 0.362168)
(10, 0.358191)
(11, 0.235373)
(12, 0.260458)
(13, 0.242502)
(14, 0.170482)
(15, 0.163713)
(16, 0.138964)
(17, 0.130559)
(18, 0.113351)
(19, 0.0774709)
(20, 0.0683271)
(21, 0.0541147)
(22, 0.0340696)
(23, 0.0281969)
(24, 0.0195368)
(25, 0.0144457)
(26, 0.00970627)
(27, 0.00499676)
(28, 0.00392127)
(29, 0.00386322)
(30, 0.00320159)
(31, 0.00302342)
(32, 0.00126414)
(33, 0.000871586)
(34, 0.000855002)
(35, 0.000621488)
(36, 0.000353392)
(37, 0.000269922)
(38, 0.000234962)
(39, 0.000157954)
(40, 0.000106396)
(41, 6.55235e-05)
(42, 4.36402e-05)
(43, 2.7811e-05)
(44, 1.91122e-05)
(45, 1.51037e-05)
(46, 7.45748e-06)
(47, 4.99022e-06)
(48, 3.41809e-06)
(49, 1.55625e-06)
(50, 1.19917e-06)
(51, 5.50448e-07)
(52, 4.29523e-07)
(53, 3.2473e-07)
(54, 2.89552e-07)
(55, 3.49591e-07)
(56, 2.31411e-07)
(57, 2.25289e-07)
(58, 1.27421e-07)
(59, 1.08922e-07)
(60, 7.6092e-08)
(61, 5.56655e-08)
(62, 3.28661e-08)
(63, 2.26145e-08)
(64, 1.41367e-08)
(65, 1.0241e-08)
(66, 4.73157e-09)
(67, 3.53018e-09)
}; \label{plot:res}
\addplot[ 
line width=1.5pt, color= myteal,  solid, mark=+,  mark repeat= 7, mark phase = 3]
coordinates { 
(1, 0.830069)
(2, 0.950043)
(3, 0.95773)
(4, 1.40134)
(5, 1.12022)
(6, 0.837534)
(7, 1.03204)
(8, 0.793578)
(9, 0.706318)
(10, 0.65861)
(11, 0.402489)
(12, 0.451265)
(13, 0.450699)
(14, 0.321337)
(15, 0.306468)
(16, 0.250191)
(17, 0.226677)
(18, 0.201426)
(19, 0.141579)
(20, 0.133737)
(21, 0.110016)
(22, 0.0735956)
(23, 0.0596398)
(24, 0.0454681)
(25, 0.0371367)
(26, 0.0284984)
(27, 0.0249576)
(28, 0.023318)
(29, 0.0239692)
(30, 0.0230768)
(31, 0.0236484)
(32, 0.0227581)
(33, 0.0229519)
(34, 0.0227696)
(35, 0.0228867)
(36, 0.0227839)
(37, 0.0228238)
(38, 0.0228025)
(39, 0.0228021)
(40, 0.0228114)
(41, 0.0227987)
(42, 0.022809)
(43, 0.0228013)
(44, 0.022806)
(45, 0.0228034)
(46, 0.0228044)
(47, 0.0228042)
(48, 0.0228043)
(49, 0.0228041)
(50, 0.0228044)
(51, 0.0228042)
(52, 0.0228043)
(53, 0.0228042)
(54, 0.0228043)
(55, 0.0228042)
(56, 0.0228043)
(57, 0.0228042)
(58, 0.0228043)
(59, 0.0228042)
(60, 0.0228042)
(61, 0.0228042)
(62, 0.0228042)
(63, 0.0228042)
(64, 0.0228042)
(65, 0.0228042)
(66, 0.0228042)
(67, 0.0228042)
}; \label{plot:etarf}

%% file: Data/test1eta.tex
\addplot[ 
line width=2.5pt, color=gray, solid ]
coordinates { 
(1, 2.39942)
(2, 2.13278)
(3, 1.88453)
(4, 1.67555)
(5, 1.42952)
(6, 1.25535)
(7, 1.10014)
(8, 0.966092)
(9, 0.87457)
(10, 0.764609)
(11, 0.665718)
(12, 0.56734)
(13, 0.461182)
(14, 0.398422)
(15, 0.337207)
(16, 0.281733)
(17, 0.231257)
(18, 0.168413)
(19, 0.134542)
(20, 0.103266)
(21, 0.0735124)
(22, 0.0585249)
(23, 0.0405522)
(24, 0.0282396)
(25, 0.0208284)
(26, 0.013798)
(27, 0.010938)
(28, 0.00913238)
(29, 0.00783348)
(30, 0.00656679)
(31, 0.00516038)
(32, 0.00453362)
(33, 0.00431525)
(34, 0.00421971)
(35, 0.0041132)
(36, 0.00406009)
(37, 0.00405351)
(38, 0.00403467)
(39, 0.00402552)
(40, 0.00403026)
(41, 0.00402091)
(42, 0.00402654)
(43, 0.00402264)
(44, 0.00402384)
(45, 0.00402408)
(46, 0.0040232)
(47, 0.00402389)
(48, 0.00402348)
(49, 0.00402361)
(50, 0.00402363)
(51, 0.00402358)
(52, 0.00402362)
(53, 0.00402359)
(54, 0.00402361)
(55, 0.0040236)
(56, 0.0040236)
(57, 0.0040236)
(58, 0.00402361)
(59, 0.0040236)
(60, 0.00402361)
(61, 0.0040236)
(62, 0.0040236)
(63, 0.0040236)
(64, 0.0040236)
(65, 0.0040236)
(66, 0.0040236)
(67, 0.0040236)
};
\addplot[ 
line width=2.5pt, color=mygreen, dashed]
coordinates { 
(1, 2.39938)
(2, 2.13281)
(3, 1.88452)
(4, 1.67554)
(5, 1.42953)
(6, 1.25534)
(7, 1.10012)
(8, 0.966099)
(9, 0.874556)
(10, 0.764575)
(11, 0.665746)
(12, 0.567247)
(13, 0.461247)
(14, 0.398369)
(15, 0.337147)
(16, 0.281793)
(17, 0.231108)
(18, 0.168435)
(19, 0.134499)
(20, 0.103058)
(21, 0.0736042)
(22, 0.0582677)
(23, 0.0403402)
(24, 0.0280843)
(25, 0.020238)
(26, 0.0134155)
(27, 0.0100535)
(28, 0.00826978)
(29, 0.00669102)
(30, 0.00516815)
(31, 0.00333002)
(32, 0.00201605)
(33, 0.00160953)
(34, 0.00125477)
(35, 0.000809996)
(36, 0.000599238)
(37, 0.000449127)
(38, 0.00030943)
(39, 0.000201169)
(40, 0.000134574)
(41, 9.01569e-05)
(42, 5.62338e-05)
(43, 3.64938e-05)
(44, 2.56074e-05)
(45, 1.66079e-05)
(46, 9.93169e-06)
(47, 6.30005e-06)
(48, 3.5121e-06)
(49, 2.19678e-06)
(50, 1.46223e-06)
(51, 1.04715e-06)
(52, 8.89395e-07)
(53, 7.82104e-07)
(54, 7.12311e-07)
(55, 5.61168e-07)
(56, 4.34979e-07)
(57, 2.9276e-07)
(58, 2.15875e-07)
(59, 1.51793e-07)
(60, 1.07064e-07)
(61, 6.83493e-08)
(62, 4.50855e-08)
(63, 2.92928e-08)
(64, 1.91922e-08)
(65, 1.12589e-08)
(66, 7.56354e-09)
(67, 5.54618e-09)
}; 
\addplot[ 
line width=1.5pt, color= myolive,  solid, mark=square, mark repeat= 7, mark phase = 6]
coordinates { 
(1, 2.30517)
(2, 2.05599)
(3, 1.82721)
(4, 1.6275)
(5, 1.38921)
(6, 1.22331)
(7, 1.07557)
(8, 0.951303)
(9, 0.864152)
(10, 0.757598)
(11, 0.661665)
(12, 0.564246)
(13, 0.459479)
(14, 0.397378)
(15, 0.336539)
(16, 0.281474)
(17, 0.230889)
(18, 0.168232)
(19, 0.134332)
(20, 0.102929)
(21, 0.0735288)
(22, 0.0582328)
(23, 0.040308)
(24, 0.0280563)
(25, 0.0202218)
(26, 0.0134021)
(27, 0.0100435)
(28, 0.00826399)
(29, 0.006688)
(30, 0.00516639)
(31, 0.0033288)
(32, 0.00201527)
(33, 0.00160912)
(34, 0.0012545)
(35, 0.000809825)
(36, 0.000599155)
(37, 0.000449082)
(38, 0.00030941)
(39, 0.000201157)
(40, 0.000134566)
(41, 9.01508e-05)
(42, 5.62268e-05)
(43, 3.64854e-05)
(44, 2.55975e-05)
(45, 1.65984e-05)
(46, 9.92216e-06)
(47, 6.29325e-06)
(48, 3.50546e-06)
(49, 2.19153e-06)
(50, 1.4583e-06)
(51, 1.04492e-06)
(52, 8.88251e-07)
(53, 7.81555e-07)
(54, 7.12052e-07)
(55, 5.61055e-07)
(56, 4.34914e-07)
(57, 2.92708e-07)
}; \label{plot:etaalg}
\addplot[ 
line width=1.5pt, color= mypurple,  solid, mark=o,  mark repeat= 7, mark phase = 1]
coordinates { 
(1, 0.605116)
(2, 0.664672)
(3, 0.691164)
(4, 0.97723)
(5, 0.763303)
(6, 0.594384)
(7, 0.738598)
(8, 0.571415)
(9, 0.524235)
(10, 0.50004)
(11, 0.307976)
(12, 0.347801)
(13, 0.336256)
(14, 0.245642)
(15, 0.240522)
(16, 0.201174)
(17, 0.187439)
(18, 0.153253)
(19, 0.107403)
(20, 0.0954296)
(21, 0.0747834)
(22, 0.0546544)
(23, 0.0407109)
(24, 0.0330291)
(25, 0.0280172)
(26, 0.0204108)
(27, 0.0202466)
(28, 0.0177191)
(29, 0.0188847)
(30, 0.0181206)
(31, 0.0183279)
(32, 0.0179352)
(33, 0.0178957)
(34, 0.0178729)
(35, 0.0179474)
(36, 0.0178207)
(37, 0.0179001)
(38, 0.0178522)
(39, 0.0178567)
(40, 0.0178804)
(41, 0.0178496)
(42, 0.0178739)
(43, 0.0178589)
(44, 0.0178655)
(45, 0.0178641)
(46, 0.017863)
(47, 0.0178646)
(48, 0.0178634)
(49, 0.0178639)
(50, 0.017864)
(51, 0.0178637)
(52, 0.0178639)
(53, 0.0178638)
(54, 0.0178639)
(55, 0.0178638)
(56, 0.0178638)
(57, 0.0178638)
(58, 0.0178638)
(59, 0.0178638)
(60, 0.0178638)
(61, 0.0178638)
(62, 0.0178638)
(63, 0.0178638)
(64, 0.0178638)
(65, 0.0178638)
(66, 0.0178638)
(67, 0.0178638)
}; \label{plot:etar}
\addplot[ 
line width=1.5pt, color= myrose ,   solid, mark=x, mark repeat= 7, mark phase = 8]
coordinates { 
(1, 0.106656)
(2, 0.156608)
(3, 0.179917)
(4, 0.269205)
(5, 0.228775)
(6, 0.171983)
(7, 0.209442)
(8, 0.155483)
(9, 0.130679)
(10, 0.107208)
(11, 0.0579275)
(12, 0.0702387)
(13, 0.077665)
(14, 0.0565715)
(15, 0.0548441)
(16, 0.0477684)
(17, 0.0459406)
(18, 0.0400506)
(19, 0.0270266)
(20, 0.0230831)
(21, 0.0186032)
(22, 0.0111761)
(23, 0.0084058)
(24, 0.0065475)
(25, 0.00558638)
(26, 0.00427031)
(27, 0.00283119)
(28, 0.00248395)
(29, 0.00260341)
(30, 0.00225727)
(31, 0.00245581)
(32, 0.00218991)
(33, 0.00227924)
(34, 0.00218671)
(35, 0.00225196)
(36, 0.00220553)
(37, 0.00221561)
(38, 0.00221501)
(39, 0.00220869)
(40, 0.00221617)
(41, 0.00220911)
(42, 0.00221417)
(43, 0.00221069)
(44, 0.00221275)
(45, 0.00221145)
(46, 0.0022122)
(47, 0.00221177)
(48, 0.00221211)
(49, 0.00221188)
(50, 0.00221201)
(51, 0.00221193)
(52, 0.00221198)
(53, 0.00221194)
(54, 0.00221197)
(55, 0.00221194)
(56, 0.00221196)
(57, 0.00221195)
(58, 0.00221196)
(59, 0.00221195)
(60, 0.00221196)
(61, 0.00221195)
(62, 0.00221195)
(63, 0.00221196)
(64, 0.00221195)
(65, 0.00221195)
(66, 0.00221195)
(67, 0.00221195)
}; \label{plot:bdm}

%% file: Data/test1tau.tex
\addplot[ 
line width=1.5pt, color= myteal,  solid, mark=+,  mark repeat= 71, mark phase = 73]
coordinates { 
(3, 2.71845)
(3.054108e+00, 2.71845)
(3.108216e+00, 2.71845)
(3.162325e+00, 2.71845)
(3.216433e+00, 2.71845)
(3.270541e+00, 2.71845)
(3.324649e+00, 2.71845)
(3.378758e+00, 2.71845)
(3.432866e+00, 2.71845)
(3.486974e+00, 2.71845)
(3.541082e+00, 2.71845)
(3.595190e+00, 2.71845)
(3.649299e+00, 2.71845)
(3.703407e+00, 2.71845)
(3.757515e+00, 2.71845)
(3.811623e+00, 2.71845)
(3.865731e+00, 2.71845)
(3.919840e+00, 2.71845)
(3.973948e+00, 2.71845)
(4.028056e+00, 2.71845)
(4.082164e+00, 2.71845)
(4.136273e+00, 2.71845)
(4.190381e+00, 2.71845)
(4.244489e+00, 2.71845)
(4.298597e+00, 2.71845)
(4.352705e+00, 2.71845)
(4.406814e+00, 2.71845)
(4.460922e+00, 2.71845)
(4.515030e+00, 2.71845)
(4.569138e+00, 2.71845)
(4.623246e+00, 2.71845)
(4.677355e+00, 2.71845)
(4.731463e+00, 2.71845)
(4.785571e+00, 2.71845)
(4.839679e+00, 2.71845)
(4.893788e+00, 2.71845)
(4.947896e+00, 2.71845)
(5.002004e+00, 2.2697)
(5.056112e+00, 2.2697)
(5.110220e+00, 2.2697)
(5.164329e+00, 2.2697)
(5.218437e+00, 2.2697)
(5.272545e+00, 2.2697)
(5.326653e+00, 2.2697)
(5.380762e+00, 2.2697)
(5.434870e+00, 2.2697)
(5.488978e+00, 2.2697)
(5.543086e+00, 2.2697)
(5.597194e+00, 2.2697)
(5.651303e+00, 2.2697)
(5.705411e+00, 2.2697)
(5.759519e+00, 2.2697)
(5.813627e+00, 2.2697)
(5.867735e+00, 2.2697)
(5.921844e+00, 2.2697)
(5.975952e+00, 1.94688)
(6.030060e+00, 1.94688)
(6.084168e+00, 1.94688)
(6.138277e+00, 1.94688)
(6.192385e+00, 1.94688)
(6.246493e+00, 1.63207)
(6.300601e+00, 1.63207)
(6.354709e+00, 1.63207)
(6.408818e+00, 1.63207)
(6.462926e+00, 1.63207)
(6.517034e+00, 1.63207)
(6.571142e+00, 1.63207)
(6.625251e+00, 1.63207)
(6.679359e+00, 1.63207)
(6.733467e+00, 1.63207)
(6.787575e+00, 1.63207)
(6.841683e+00, 1.63207)
(6.895792e+00, 1.63207)
(6.949900e+00, 1.63207)
(7.004008e+00, 1.63207)
(7.058116e+00, 1.63207)
(7.112224e+00, 1.63207)
(7.166333e+00, 1.63207)
(7.220441e+00, 1.28253)
(7.274549e+00, 1.28253)
(7.328657e+00, 1.28253)
(7.382766e+00, 1.28253)
(7.436874e+00, 1.28253)
(7.490982e+00, 1.28253)
(7.545090e+00, 1.28253)
(7.599198e+00, 1.28253)
(7.653307e+00, 1.28253)
(7.707415e+00, 1.28253)
(7.761523e+00, 1.28253)
(7.815631e+00, 1.28253)
(7.869739e+00, 1.12676)
(7.923848e+00, 1.12676)
(7.977956e+00, 1.12676)
(8.032064e+00, 1.12676)
(8.086172e+00, 1.12676)
(8.140281e+00, 1.12676)
(8.194389e+00, 1.12676)
(8.248497e+00, 1.12676)
(8.302605e+00, 1.12676)
(8.356713e+00, 1.12676)
(8.410822e+00, 1.12676)
(8.464930e+00, 1.12676)
(8.519038e+00, 1.12676)
(8.573146e+00, 1.12676)
(8.627255e+00, 1.12676)
(8.681363e+00, 1.12676)
(8.735471e+00, 1.12676)
(8.789579e+00, 1.12676)
(8.843687e+00, 1.12676)
(8.897796e+00, 1.12676)
(8.951904e+00, 1.12676)
(9.006012e+00, 1.12676)
(9.060120e+00, 1.12676)
(9.114228e+00, 1.12676)
(9.168337e+00, 1.12676)
(9.222445e+00, 1.12676)
(9.276553e+00, 1.12676)
(9.330661e+00, 1.12676)
(9.384770e+00, 1.12676)
(9.438878e+00, 1.12676)
(9.492986e+00, 1.12676)
(9.547094e+00, 1.12676)
(9.601202e+00, 1.12676)
(9.655311e+00, 1.12676)
(9.709419e+00, 1.12676)
(9.763527e+00, 1.12676)
(9.817635e+00, 1.12676)
(9.871743e+00, 1.12676)
(9.925852e+00, 1.12676)
(9.979960e+00, 1.12676)
(1.003407e+01, 1.12676)
(1.008818e+01, 1.12676)
(1.014228e+01, 1.12676)
(1.019639e+01, 1.12676)
(1.025050e+01, 1.12676)
(1.030461e+01, 1.12676)
(1.035872e+01, 1.12676)
(1.041283e+01, 1.12676)
(1.046693e+01, 1.12676)
(1.052104e+01, 1.12676)
(1.057515e+01, 1.12676)
(1.062926e+01, 1.12676)
(1.068337e+01, 1.12676)
(1.073747e+01, 1.12676)
(1.079158e+01, 1.12676)
(1.084569e+01, 1.12676)
(1.089980e+01, 1.12676)
(1.095391e+01, 1.12676)
(1.100802e+01, 1.12676)
(1.106212e+01, 1.12676)
(1.111623e+01, 1.12676)
(1.117034e+01, 1.12676)
(1.122445e+01, 1.12676)
(1.127856e+01, 1.12676)
(1.133267e+01, 1.12676)
(1.138677e+01, 1.12676)
(1.144088e+01, 1.12676)
(1.149499e+01, 1.12676)
(1.154910e+01, 1.12676)
(1.160321e+01, 1.12676)
(1.165731e+01, 1.12676)
(1.171142e+01, 1.12676)
(1.176553e+01, 1.12676)
(1.181964e+01, 1.12676)
(1.187375e+01, 1.12676)
(1.192786e+01, 1.12676)
(1.198196e+01, 1.12676)
(1.203607e+01, 1.12676)
(1.209018e+01, 1.12676)
(1.214429e+01, 1.12676)
(1.219840e+01, 1.12676)
(1.225251e+01, 1.12676)
(1.230661e+01, 1.12676)
(1.236072e+01, 1.12676)
(1.241483e+01, 1.12676)
(1.246894e+01, 1.12676)
(1.252305e+01, 1.12676)
(1.257715e+01, 1.12676)
(1.263126e+01, 1.12676)
(1.268537e+01, 1.12676)
(1.273948e+01, 1.12676)
(1.279359e+01, 1.12676)
(1.284770e+01, 1.12676)
(1.290180e+01, 1.12676)
(1.295591e+01, 1.12676)
(1.301002e+01, 1.12676)
(1.306413e+01, 1.12676)
(1.311824e+01, 1.12676)
(1.317234e+01, 1.12676)
(1.322645e+01, 1.12676)
(1.328056e+01, 1.12676)
(1.333467e+01, 1.12676)
(1.338878e+01, 1.12676)
(1.344289e+01, 1.12676)
(1.349699e+01, 1.12676)
(1.355110e+01, 1.12676)
(1.360521e+01, 1.12676)
(1.365932e+01, 1.12676)
(1.371343e+01, 1.12676)
(1.376754e+01, 1.12676)
(1.382164e+01, 1.12676)
(1.387575e+01, 1.12676)
(1.392986e+01, 1.12676)
(1.398397e+01, 1.12676)
(1.403808e+01, 1.12676)
(1.409218e+01, 1.12676)
(1.414629e+01, 1.12676)
(1.420040e+01, 1.12676)
(1.425451e+01, 1.12676)
(1.430862e+01, 1.12676)
(1.436273e+01, 1.12676)
(1.441683e+01, 1.12676)
(1.447094e+01, 1.12676)
(1.452505e+01, 1.12676)
(1.457916e+01, 1.12676)
(1.463327e+01, 1.12676)
(1.468737e+01, 1.12676)
(1.474148e+01, 1.12676)
(1.479559e+01, 1.12676)
(1.484970e+01, 1.12676)
(1.490381e+01, 1.12676)
(1.495792e+01, 1.12676)
(1.501202e+01, 1.12676)
(1.506613e+01, 1.12676)
(1.512024e+01, 1.12676)
(1.517435e+01, 1.12676)
(1.522846e+01, 1.12676)
(1.528257e+01, 1.12676)
(1.533667e+01, 1.12676)
(1.539078e+01, 1.12676)
(1.544489e+01, 1.12676)
(1.549900e+01, 1.12676)
(1.555311e+01, 1.12676)
(1.560721e+01, 1.12676)
(1.566132e+01, 1.12676)
(1.571543e+01, 1.12676)
(1.576954e+01, 1.12676)
(1.582365e+01, 1.12676)
(1.587776e+01, 1.12676)
(1.593186e+01, 1.12676)
(1.598597e+01, 1.12676)
(1.604008e+01, 1.12676)
(1.609419e+01, 1.12676)
(1.614830e+01, 1.12676)
(1.620240e+01, 1.12676)
(1.625651e+01, 1.12676)
(1.631062e+01, 1.12676)
(1.636473e+01, 1.12676)
(1.641884e+01, 1.12676)
(1.647295e+01, 1.12676)
(1.652705e+01, 1.12676)
(1.658116e+01, 1.12676)
(1.663527e+01, 1.12676)
(1.668938e+01, 1.12676)
(1.674349e+01, 1.12676)
(1.679760e+01, 1.12676)
(1.685170e+01, 1.12676)
(1.690581e+01, 1.12676)
(1.695992e+01, 1.12676)
(1.701403e+01, 1.12676)
(1.706814e+01, 1.12676)
(1.712224e+01, 1.12676)
(1.717635e+01, 1.12676)
(1.723046e+01, 1.12676)
(1.728457e+01, 1.12676)
(1.733868e+01, 1.12676)
(1.739279e+01, 1.12676)
(1.744689e+01, 1.12676)
(1.750100e+01, 1.12676)
(1.755511e+01, 1.12676)
(1.760922e+01, 1.12676)
(1.766333e+01, 1.12676)
(1.771743e+01, 1.12676)
(1.777154e+01, 1.12676)
(1.782565e+01, 1.12676)
(1.787976e+01, 1.12676)
(1.793387e+01, 1.12676)
(1.798798e+01, 1.12676)
(1.804208e+01, 1.07248)
(1.809619e+01, 1.07248)
(1.815030e+01, 1.07248)
(1.820441e+01, 1.07248)
(1.825852e+01, 1.07248)
(1.831263e+01, 1.07248)
(1.836673e+01, 1.07248)
(1.842084e+01, 1.07248)
(1.847495e+01, 1.07248)
(1.852906e+01, 1.07248)
(1.858317e+01, 1.07248)
(1.863727e+01, 1.07248)
(1.869138e+01, 1.07248)
(1.874549e+01, 1.07248)
(1.879960e+01, 1.07248)
(1.885371e+01, 1.07248)
(1.890782e+01, 1.07248)
(1.896192e+01, 1.07248)
(1.901603e+01, 1.07248)
(1.907014e+01, 1.07248)
(1.912425e+01, 1.07248)
(1.917836e+01, 1.07248)
(1.923246e+01, 1.07248)
(1.928657e+01, 1.07248)
(1.934068e+01, 1.07248)
(1.939479e+01, 1.07248)
(1.944890e+01, 1.07248)
(1.950301e+01, 1.07248)
(1.955711e+01, 1.07248)
(1.961122e+01, 1.07248)
(1.966533e+01, 1.07248)
(1.971944e+01, 1.07248)
(1.977355e+01, 1.07248)
(1.982766e+01, 1.07248)
(1.988176e+01, 1.07248)
(1.993587e+01, 1.07248)
(1.998998e+01, 1.07248)
(2.004409e+01, 1.07248)
(2.009820e+01, 1.07248)
(2.015230e+01, 1.07248)
(2.020641e+01, 1.07248)
(2.026052e+01, 1.07248)
(2.031463e+01, 1.07248)
(2.036874e+01, 1.07248)
(2.042285e+01, 1.07248)
(2.047695e+01, 1.07248)
(2.053106e+01, 1.07248)
(2.058517e+01, 1.07248)
(2.063928e+01, 1.07248)
(2.069339e+01, 1.07248)
(2.074749e+01, 1.07248)
(2.080160e+01, 1.07248)
(2.085571e+01, 1.07248)
(2.090982e+01, 1.07248)
(2.096393e+01, 1.07248)
(2.101804e+01, 1.07248)
(2.107214e+01, 1.07248)
(2.112625e+01, 1.07248)
(2.118036e+01, 1.07248)
(2.123447e+01, 1.07248)
(2.128858e+01, 1.07248)
(2.134269e+01, 1.07248)
(2.139679e+01, 1.07248)
(2.145090e+01, 1.07248)
(2.150501e+01, 1.07248)
(2.155912e+01, 1.07248)
(2.161323e+01, 1.07248)
(2.166733e+01, 1.07248)
(2.172144e+01, 1.07248)
(2.177555e+01, 1.07248)
(2.182966e+01, 1.07248)
(2.188377e+01, 1.07248)
(2.193788e+01, 1.07248)
(2.199198e+01, 1.07248)
(2.204609e+01, 1.07248)
(2.210020e+01, 1.07248)
(2.215431e+01, 1.07248)
(2.220842e+01, 1.07248)
(2.226253e+01, 1.07248)
(2.231663e+01, 1.07248)
(2.237074e+01, 1.07248)
(2.242485e+01, 1.07248)
(2.247896e+01, 1.07248)
(2.253307e+01, 1.07248)
(2.258717e+01, 1.07248)
(2.264128e+01, 1.07248)
(2.269539e+01, 1.07248)
(2.274950e+01, 1.07248)
(2.280361e+01, 1.07248)
(2.285772e+01, 1.07248)
(2.291182e+01, 1.07248)
(2.296593e+01, 1.07248)
(2.302004e+01, 1.07248)
(2.307415e+01, 1.07248)
(2.312826e+01, 1.07248)
(2.318236e+01, 1.07248)
(2.323647e+01, 1.07248)
(2.329058e+01, 1.07248)
(2.334469e+01, 1.07248)
(2.339880e+01, 1.07248)
(2.345291e+01, 1.07248)
(2.350701e+01, 1.07248)
(2.356112e+01, 1.07248)
(2.361523e+01, 1.07248)
(2.366934e+01, 1.07248)
(2.372345e+01, 1.07248)
(2.377756e+01, 1.07248)
(2.383166e+01, 1.07248)
(2.388577e+01, 1.07248)
(2.393988e+01, 1.07248)
(2.399399e+01, 1.07248)
(2.404810e+01, 1.07248)
(2.410220e+01, 1.07248)
(2.415631e+01, 1.07248)
(2.421042e+01, 1.07248)
(2.426453e+01, 1.07248)
(2.431864e+01, 1.07248)
(2.437275e+01, 1.07248)
(2.442685e+01, 1.07248)
(2.448096e+01, 1.07248)
(2.453507e+01, 1.07248)
(2.458918e+01, 1.07248)
(2.464329e+01, 1.07248)
(2.469739e+01, 1.07248)
(2.475150e+01, 1.07248)
(2.480561e+01, 1.07248)
(2.485972e+01, 1.07248)
(2.491383e+01, 1.07248)
(2.496794e+01, 1.07248)
(2.502204e+01, 1.07248)
(2.507615e+01, 1.07248)
(2.513026e+01, 1.07248)
(2.518437e+01, 1.07248)
(2.523848e+01, 1.07248)
(2.529259e+01, 1.07248)
(2.534669e+01, 1.07248)
(2.540080e+01, 1.07248)
(2.545491e+01, 1.07248)
(2.550902e+01, 1.07248)
(2.556313e+01, 1.07248)
(2.561723e+01, 1.07248)
(2.567134e+01, 1.07248)
(2.572545e+01, 1.07248)
(2.577956e+01, 1.07248)
(2.583367e+01, 1.07248)
(2.588778e+01, 1.07248)
(2.594188e+01, 1.07248)
(2.599599e+01, 1.07248)
(2.605010e+01, 1.07248)
(2.610421e+01, 1.07248)
(2.615832e+01, 1.07248)
(2.621242e+01, 1.07248)
(2.626653e+01, 1.07248)
(2.632064e+01, 1.07248)
(2.637475e+01, 1.04874)
(2.642886e+01, 1.04874)
(2.648297e+01, 1.04874)
(2.653707e+01, 1.04874)
(2.659118e+01, 1.04874)
(2.664529e+01, 1.02227)
(2.669940e+01, 1.02227)
(2.675351e+01, 1.02227)
(2.680762e+01, 1.02227)
(2.686172e+01, 1.02227)
(2.691583e+01, 1.02227)
(2.696994e+01, 1.02227)
(2.702405e+01, 1.02227)
(2.707816e+01, 1.02227)
(2.713226e+01, 1.02227)
(2.718637e+01, 1.02227)
(2.724048e+01, 1.02227)
(2.729459e+01, 1.02227)
(2.734870e+01, 1.02227)
(2.740281e+01, 1.02227)
(2.745691e+01, 1.02227)
(2.751102e+01, 1.02227)
(2.756513e+01, 1.02227)
(2.761924e+01, 1.02227)
(2.767335e+01, 1.02227)
(2.772745e+01, 1.02227)
(2.778156e+01, 1.02227)
(2.783567e+01, 1.02227)
(2.788978e+01, 1.02227)
(2.794389e+01, 1.02227)
(2.799800e+01, 1.02227)
(2.805210e+01, 1.02227)
(2.810621e+01, 1.02227)
(2.816032e+01, 1.02227)
(2.821443e+01, 1.02227)
(2.826854e+01, 1.02227)
(2.832265e+01, 1.02227)
(2.837675e+01, 1.02227)
(2.843086e+01, 1.02227)
(2.848497e+01, 1.02227)
(2.853908e+01, 1.02227)
(2.859319e+01, 1.02227)
(2.864729e+01, 1.02227)
(2.870140e+01, 1.02227)
(2.875551e+01, 1.02227)
(2.880962e+01, 1.02227)
(2.886373e+01, 1.02227)
(2.891784e+01, 1.02227)
(2.897194e+01, 1.02227)
(2.902605e+01, 1.02227)
(2.908016e+01, 1.02227)
(2.913427e+01, 1.02227)
(2.918838e+01, 1.02227)
(2.924248e+01, 1.02227)
(2.929659e+01, 1.02227)
(2.935070e+01, 1.02227)
(2.940481e+01, 1.02227)
(2.945892e+01, 1.02227)
(2.951303e+01, 1.02227)
(2.956713e+01, 1.02227)
(2.962124e+01, 1.02227)
(2.967535e+01, 1.02227)
(2.972946e+01, 1.02227)
(2.978357e+01, 1.02227)
(2.983768e+01, 1.02227)
(2.989178e+01, 1.02227)
(2.994589e+01, 1.02227)
(30, 1.02227)
}; 
\addplot[ 
line width=1.5pt, color= mypurple,  solid, mark=o,  mark repeat= 71, mark phase = 71]
coordinates { 
(3, 1.63207)
(3.054108e+00, 1.63207)
(3.108216e+00, 1.63207)
(3.162325e+00, 1.63207)
(3.216433e+00, 1.63207)
(3.270541e+00, 1.63207)
(3.324649e+00, 1.63207)
(3.378758e+00, 1.63207)
(3.432866e+00, 1.63207)
(3.486974e+00, 1.63207)
(3.541082e+00, 1.28253)
(3.595190e+00, 1.28253)
(3.649299e+00, 1.28253)
(3.703407e+00, 1.28253)
(3.757515e+00, 1.28253)
(3.811623e+00, 1.28253)
(3.865731e+00, 1.28253)
(3.919840e+00, 1.28253)
(3.973948e+00, 1.28253)
(4.028056e+00, 1.28253)
(4.082164e+00, 1.28253)
(4.136273e+00, 1.28253)
(4.190381e+00, 1.28253)
(4.244489e+00, 1.28253)
(4.298597e+00, 1.28253)
(4.352705e+00, 1.28253)
(4.406814e+00, 1.28253)
(4.460922e+00, 1.28253)
(4.515030e+00, 1.28253)
(4.569138e+00, 1.28253)
(4.623246e+00, 1.28253)
(4.677355e+00, 1.28253)
(4.731463e+00, 1.28253)
(4.785571e+00, 1.28253)
(4.839679e+00, 1.28253)
(4.893788e+00, 1.28253)
(4.947896e+00, 1.28253)
(5.002004e+00, 1.28253)
(5.056112e+00, 1.28253)
(5.110220e+00, 1.28253)
(5.164329e+00, 1.28253)
(5.218437e+00, 1.28253)
(5.272545e+00, 1.28253)
(5.326653e+00, 1.28253)
(5.380762e+00, 1.28253)
(5.434870e+00, 1.28253)
(5.488978e+00, 1.28253)
(5.543086e+00, 1.12676)
(5.597194e+00, 1.12676)
(5.651303e+00, 1.12676)
(5.705411e+00, 1.12676)
(5.759519e+00, 1.12676)
(5.813627e+00, 1.12676)
(5.867735e+00, 1.12676)
(5.921844e+00, 1.12676)
(5.975952e+00, 1.12676)
(6.030060e+00, 1.12676)
(6.084168e+00, 1.12676)
(6.138277e+00, 1.12676)
(6.192385e+00, 1.12676)
(6.246493e+00, 1.12676)
(6.300601e+00, 1.12676)
(6.354709e+00, 1.12676)
(6.408818e+00, 1.12676)
(6.462926e+00, 1.12676)
(6.517034e+00, 1.12676)
(6.571142e+00, 1.12676)
(6.625251e+00, 1.12676)
(6.679359e+00, 1.12676)
(6.733467e+00, 1.12676)
(6.787575e+00, 1.12676)
(6.841683e+00, 1.12676)
(6.895792e+00, 1.12676)
(6.949900e+00, 1.12676)
(7.004008e+00, 1.12676)
(7.058116e+00, 1.12676)
(7.112224e+00, 1.12676)
(7.166333e+00, 1.12676)
(7.220441e+00, 1.12676)
(7.274549e+00, 1.12676)
(7.328657e+00, 1.12676)
(7.382766e+00, 1.12676)
(7.436874e+00, 1.12676)
(7.490982e+00, 1.12676)
(7.545090e+00, 1.12676)
(7.599198e+00, 1.12676)
(7.653307e+00, 1.12676)
(7.707415e+00, 1.12676)
(7.761523e+00, 1.12676)
(7.815631e+00, 1.12676)
(7.869739e+00, 1.12676)
(7.923848e+00, 1.12676)
(7.977956e+00, 1.12676)
(8.032064e+00, 1.12676)
(8.086172e+00, 1.12676)
(8.140281e+00, 1.12676)
(8.194389e+00, 1.12676)
(8.248497e+00, 1.12676)
(8.302605e+00, 1.12676)
(8.356713e+00, 1.12676)
(8.410822e+00, 1.12676)
(8.464930e+00, 1.12676)
(8.519038e+00, 1.12676)
(8.573146e+00, 1.12676)
(8.627255e+00, 1.12676)
(8.681363e+00, 1.12676)
(8.735471e+00, 1.12676)
(8.789579e+00, 1.12676)
(8.843687e+00, 1.12676)
(8.897796e+00, 1.12676)
(8.951904e+00, 1.07248)
(9.006012e+00, 1.07248)
(9.060120e+00, 1.07248)
(9.114228e+00, 1.07248)
(9.168337e+00, 1.07248)
(9.222445e+00, 1.07248)
(9.276553e+00, 1.07248)
(9.330661e+00, 1.07248)
(9.384770e+00, 1.07248)
(9.438878e+00, 1.07248)
(9.492986e+00, 1.07248)
(9.547094e+00, 1.07248)
(9.601202e+00, 1.07248)
(9.655311e+00, 1.07248)
(9.709419e+00, 1.07248)
(9.763527e+00, 1.07248)
(9.817635e+00, 1.07248)
(9.871743e+00, 1.07248)
(9.925852e+00, 1.07248)
(9.979960e+00, 1.07248)
(1.003407e+01, 1.07248)
(1.008818e+01, 1.07248)
(1.014228e+01, 1.07248)
(1.019639e+01, 1.07248)
(1.025050e+01, 1.07248)
(1.030461e+01, 1.07248)
(1.035872e+01, 1.07248)
(1.041283e+01, 1.07248)
(1.046693e+01, 1.07248)
(1.052104e+01, 1.07248)
(1.057515e+01, 1.07248)
(1.062926e+01, 1.07248)
(1.068337e+01, 1.07248)
(1.073747e+01, 1.07248)
(1.079158e+01, 1.07248)
(1.084569e+01, 1.07248)
(1.089980e+01, 1.07248)
(1.095391e+01, 1.07248)
(1.100802e+01, 1.07248)
(1.106212e+01, 1.07248)
(1.111623e+01, 1.07248)
(1.117034e+01, 1.04874)
(1.122445e+01, 1.04874)
(1.127856e+01, 1.04874)
(1.133267e+01, 1.04874)
(1.138677e+01, 1.04874)
(1.144088e+01, 1.04874)
(1.149499e+01, 1.04874)
(1.154910e+01, 1.04874)
(1.160321e+01, 1.04874)
(1.165731e+01, 1.04874)
(1.171142e+01, 1.04874)
(1.176553e+01, 1.04874)
(1.181964e+01, 1.04874)
(1.187375e+01, 1.04874)
(1.192786e+01, 1.04874)
(1.198196e+01, 1.04874)
(1.203607e+01, 1.04874)
(1.209018e+01, 1.04874)
(1.214429e+01, 1.04874)
(1.219840e+01, 1.04874)
(1.225251e+01, 1.04874)
(1.230661e+01, 1.04874)
(1.236072e+01, 1.04874)
(1.241483e+01, 1.04874)
(1.246894e+01, 1.04874)
(1.252305e+01, 1.04874)
(1.257715e+01, 1.04874)
(1.263126e+01, 1.04874)
(1.268537e+01, 1.04874)
(1.273948e+01, 1.04874)
(1.279359e+01, 1.04874)
(1.284770e+01, 1.04874)
(1.290180e+01, 1.04874)
(1.295591e+01, 1.04874)
(1.301002e+01, 1.04874)
(1.306413e+01, 1.04874)
(1.311824e+01, 1.04874)
(1.317234e+01, 1.04874)
(1.322645e+01, 1.04874)
(1.328056e+01, 1.04874)
(1.333467e+01, 1.04874)
(1.338878e+01, 1.04874)
(1.344289e+01, 1.04874)
(1.349699e+01, 1.04874)
(1.355110e+01, 1.04874)
(1.360521e+01, 1.04874)
(1.365932e+01, 1.04874)
(1.371343e+01, 1.04874)
(1.376754e+01, 1.04874)
(1.382164e+01, 1.04874)
(1.387575e+01, 1.04874)
(1.392986e+01, 1.04874)
(1.398397e+01, 1.04874)
(1.403808e+01, 1.04874)
(1.409218e+01, 1.04874)
(1.414629e+01, 1.04874)
(1.420040e+01, 1.04874)
(1.425451e+01, 1.02227)
(1.430862e+01, 1.02227)
(1.436273e+01, 1.02227)
(1.441683e+01, 1.02227)
(1.447094e+01, 1.02227)
(1.452505e+01, 1.02227)
(1.457916e+01, 1.02227)
(1.463327e+01, 1.02227)
(1.468737e+01, 1.02227)
(1.474148e+01, 1.02227)
(1.479559e+01, 1.02227)
(1.484970e+01, 1.02227)
(1.490381e+01, 1.02227)
(1.495792e+01, 1.02227)
(1.501202e+01, 1.02227)
(1.506613e+01, 1.02227)
(1.512024e+01, 1.02227)
(1.517435e+01, 1.02227)
(1.522846e+01, 1.02227)
(1.528257e+01, 1.02227)
(1.533667e+01, 1.02227)
(1.539078e+01, 1.02227)
(1.544489e+01, 1.02227)
(1.549900e+01, 1.02227)
(1.555311e+01, 1.02227)
(1.560721e+01, 1.02227)
(1.566132e+01, 1.02227)
(1.571543e+01, 1.02227)
(1.576954e+01, 1.02227)
(1.582365e+01, 1.02227)
(1.587776e+01, 1.02227)
(1.593186e+01, 1.02227)
(1.598597e+01, 1.02227)
(1.604008e+01, 1.02227)
(1.609419e+01, 1.02227)
(1.614830e+01, 1.02227)
(1.620240e+01, 1.02227)
(1.625651e+01, 1.02227)
(1.631062e+01, 1.02227)
(1.636473e+01, 1.02227)
(1.641884e+01, 1.02227)
(1.647295e+01, 1.02227)
(1.652705e+01, 1.02227)
(1.658116e+01, 1.02227)
(1.663527e+01, 1.02227)
(1.668938e+01, 1.02227)
(1.674349e+01, 1.02227)
(1.679760e+01, 1.02227)
(1.685170e+01, 1.02227)
(1.690581e+01, 1.02227)
(1.695992e+01, 1.02227)
(1.701403e+01, 1.02227)
(1.706814e+01, 1.02227)
(1.712224e+01, 1.02227)
(1.717635e+01, 1.02227)
(1.723046e+01, 1.02227)
(1.728457e+01, 1.02227)
(1.733868e+01, 1.02227)
(1.739279e+01, 1.02227)
(1.744689e+01, 1.02227)
(1.750100e+01, 1.02227)
(1.755511e+01, 1.02227)
(1.760922e+01, 1.02227)
(1.766333e+01, 1.02227)
(1.771743e+01, 1.02227)
(1.777154e+01, 1.02227)
(1.782565e+01, 1.02227)
(1.787976e+01, 1.02227)
(1.793387e+01, 1.02227)
(1.798798e+01, 1.02227)
(1.804208e+01, 1.02227)
(1.809619e+01, 1.02227)
(1.815030e+01, 1.02227)
(1.820441e+01, 1.02227)
(1.825852e+01, 1.02227)
(1.831263e+01, 1.02227)
(1.836673e+01, 1.02227)
(1.842084e+01, 1.02227)
(1.847495e+01, 1.02227)
(1.852906e+01, 1.02227)
(1.858317e+01, 1.02227)
(1.863727e+01, 1.02227)
(1.869138e+01, 1.02227)
(1.874549e+01, 1.02227)
(1.879960e+01, 1.02227)
(1.885371e+01, 1.02227)
(1.890782e+01, 1.02227)
(1.896192e+01, 1.02227)
(1.901603e+01, 1.02227)
(1.907014e+01, 1.02227)
(1.912425e+01, 1.02227)
(1.917836e+01, 1.02227)
(1.923246e+01, 1.02227)
(1.928657e+01, 1.02227)
(1.934068e+01, 1.02227)
(1.939479e+01, 1.02227)
(1.944890e+01, 1.02227)
(1.950301e+01, 1.02227)
(1.955711e+01, 1.02227)
(1.961122e+01, 1.02227)
(1.966533e+01, 1.02227)
(1.971944e+01, 1.02227)
(1.977355e+01, 1.02227)
(1.982766e+01, 1.02227)
(1.988176e+01, 1.02227)
(1.993587e+01, 1.02227)
(1.998998e+01, 1.02227)
(2.004409e+01, 1.02227)
(2.009820e+01, 1.02227)
(2.015230e+01, 1.02227)
(2.020641e+01, 1.02227)
(2.026052e+01, 1.02227)
(2.031463e+01, 1.02227)
(2.036874e+01, 1.02227)
(2.042285e+01, 1.02227)
(2.047695e+01, 1.02227)
(2.053106e+01, 1.02227)
(2.058517e+01, 1.02227)
(2.063928e+01, 1.02227)
(2.069339e+01, 1.02227)
(2.074749e+01, 1.02227)
(2.080160e+01, 1.02227)
(2.085571e+01, 1.02227)
(2.090982e+01, 1.02227)
(2.096393e+01, 1.02227)
(2.101804e+01, 1.02227)
(2.107214e+01, 1.02227)
(2.112625e+01, 1.02227)
(2.118036e+01, 1.02227)
(2.123447e+01, 1.02227)
(2.128858e+01, 1.02227)
(2.134269e+01, 1.02227)
(2.139679e+01, 1.02227)
(2.145090e+01, 1.02227)
(2.150501e+01, 1.02227)
(2.155912e+01, 1.02227)
(2.161323e+01, 1.02227)
(2.166733e+01, 1.02227)
(2.172144e+01, 1.02227)
(2.177555e+01, 1.02227)
(2.182966e+01, 1.02227)
(2.188377e+01, 1.02227)
(2.193788e+01, 1.02227)
(2.199198e+01, 1.02227)
(2.204609e+01, 1.02227)
(2.210020e+01, 1.02227)
(2.215431e+01, 1.02227)
(2.220842e+01, 1.00907)
(2.226253e+01, 1.00907)
(2.231663e+01, 1.00907)
(2.237074e+01, 1.00907)
(2.242485e+01, 1.00907)
(2.247896e+01, 1.00907)
(2.253307e+01, 1.00907)
(2.258717e+01, 1.00907)
(2.264128e+01, 1.00907)
(2.269539e+01, 1.00907)
(2.274950e+01, 1.00907)
(2.280361e+01, 1.00907)
(2.285772e+01, 1.00907)
(2.291182e+01, 1.00907)
(2.296593e+01, 1.00907)
(2.302004e+01, 1.00907)
(2.307415e+01, 1.00907)
(2.312826e+01, 1.00907)
(2.318236e+01, 1.00907)
(2.323647e+01, 1.00907)
(2.329058e+01, 1.00907)
(2.334469e+01, 1.00907)
(2.339880e+01, 1.00907)
(2.345291e+01, 1.00907)
(2.350701e+01, 1.00907)
(2.356112e+01, 1.00907)
(2.361523e+01, 1.00907)
(2.366934e+01, 1.00907)
(2.372345e+01, 1.00907)
(2.377756e+01, 1.00907)
(2.383166e+01, 1.00907)
(2.388577e+01, 1.00907)
(2.393988e+01, 1.00907)
(2.399399e+01, 1.00907)
(2.404810e+01, 1.00907)
(2.410220e+01, 1.00907)
(2.415631e+01, 1.00907)
(2.421042e+01, 1.00907)
(2.426453e+01, 1.00907)
(2.431864e+01, 1.00907)
(2.437275e+01, 1.00907)
(2.442685e+01, 1.00907)
(2.448096e+01, 1.00907)
(2.453507e+01, 1.00907)
(2.458918e+01, 1.00907)
(2.464329e+01, 1.00907)
(2.469739e+01, 1.00907)
(2.475150e+01, 1.00907)
(2.480561e+01, 1.00907)
(2.485972e+01, 1.00907)
(2.491383e+01, 1.00907)
(2.496794e+01, 1.00907)
(2.502204e+01, 1.00907)
(2.507615e+01, 1.00907)
(2.513026e+01, 1.00907)
(2.518437e+01, 1.00907)
(2.523848e+01, 1.00907)
(2.529259e+01, 1.00907)
(2.534669e+01, 1.00907)
(2.540080e+01, 1.00907)
(2.545491e+01, 1.00907)
(2.550902e+01, 1.00907)
(2.556313e+01, 1.00907)
(2.561723e+01, 1.00907)
(2.567134e+01, 1.00907)
(2.572545e+01, 1.00907)
(2.577956e+01, 1.00907)
(2.583367e+01, 1.00907)
(2.588778e+01, 1.00907)
(2.594188e+01, 1.00907)
(2.599599e+01, 1.00907)
(2.605010e+01, 1.00907)
(2.610421e+01, 1.00907)
(2.615832e+01, 1.00907)
(2.621242e+01, 1.00907)
(2.626653e+01, 1.00907)
(2.632064e+01, 1.00907)
(2.637475e+01, 1.00907)
(2.642886e+01, 1.00907)
(2.648297e+01, 1.00907)
(2.653707e+01, 1.00907)
(2.659118e+01, 1.00907)
(2.664529e+01, 1.00907)
(2.669940e+01, 1.00907)
(2.675351e+01, 1.00907)
(2.680762e+01, 1.00907)
(2.686172e+01, 1.00907)
(2.691583e+01, 1.00907)
(2.696994e+01, 1.00907)
(2.702405e+01, 1.00907)
(2.707816e+01, 1.00907)
(2.713226e+01, 1.00907)
(2.718637e+01, 1.00907)
(2.724048e+01, 1.00907)
(2.729459e+01, 1.00907)
(2.734870e+01, 1.00907)
(2.740281e+01, 1.00907)
(2.745691e+01, 1.00907)
(2.751102e+01, 1.00907)
(2.756513e+01, 1.00907)
(2.761924e+01, 1.00907)
(2.767335e+01, 1.00907)
(2.772745e+01, 1.00907)
(2.778156e+01, 1.00907)
(2.783567e+01, 1.00907)
(2.788978e+01, 1.00907)
(2.794389e+01, 1.00907)
(2.799800e+01, 1.00907)
(2.805210e+01, 1.00907)
(2.810621e+01, 1.00907)
(2.816032e+01, 1.00907)
(2.821443e+01, 1.00907)
(2.826854e+01, 1.00907)
(2.832265e+01, 1.00907)
(2.837675e+01, 1.00907)
(2.843086e+01, 1.00907)
(2.848497e+01, 1.00907)
(2.853908e+01, 1.00907)
(2.859319e+01, 1.00907)
(2.864729e+01, 1.00907)
(2.870140e+01, 1.00907)
(2.875551e+01, 1.00907)
(2.880962e+01, 1.00907)
(2.886373e+01, 1.00907)
(2.891784e+01, 1.00907)
(2.897194e+01, 1.00907)
(2.902605e+01, 1.00907)
(2.908016e+01, 1.00907)
(2.913427e+01, 1.00907)
(2.918838e+01, 1.00907)
(2.924248e+01, 1.00907)
(2.929659e+01, 1.00907)
(2.935070e+01, 1.00907)
(2.940481e+01, 1.00907)
(2.945892e+01, 1.00907)
(2.951303e+01, 1.00907)
(2.956713e+01, 1.00907)
(2.962124e+01, 1.00907)
(2.967535e+01, 1.00907)
(2.972946e+01, 1.00907)
(2.978357e+01, 1.00743)
(2.983768e+01, 1.00743)
(2.989178e+01, 1.00743)
(2.994589e+01, 1.00743)
(30, 1.00743)
}; 
\addplot[ 
line width=1.5pt, color= myrose ,   solid, mark=x, mark repeat= 71, mark phase = 77]
coordinates { 
(3, 1.00907)
(3.054108e+00, 1.00907)
(3.108216e+00, 1.00907)
(3.162325e+00, 1.00907)
(3.216433e+00, 1.00907)
(3.270541e+00, 1.00907)
(3.324649e+00, 1.00907)
(3.378758e+00, 1.00907)
(3.432866e+00, 1.00907)
(3.486974e+00, 1.00907)
(3.541082e+00, 1.00907)
(3.595190e+00, 1.00907)
(3.649299e+00, 1.00907)
(3.703407e+00, 1.00743)
(3.757515e+00, 1.00743)
(3.811623e+00, 1.00743)
(3.865731e+00, 1.00743)
(3.919840e+00, 1.00743)
(3.973948e+00, 1.00743)
(4.028056e+00, 1.00743)
(4.082164e+00, 1.00743)
(4.136273e+00, 1.00743)
(4.190381e+00, 1.00743)
(4.244489e+00, 1.00743)
(4.298597e+00, 1.00743)
(4.352705e+00, 1.00743)
(4.406814e+00, 1.00743)
(4.460922e+00, 1.00743)
(4.515030e+00, 1.00743)
(4.569138e+00, 1.00743)
(4.623246e+00, 1.00743)
(4.677355e+00, 1.00743)
(4.731463e+00, 1.00743)
(4.785571e+00, 1.00743)
(4.839679e+00, 1.00743)
(4.893788e+00, 1.00743)
(4.947896e+00, 1.00275)
(5.002004e+00, 1.00275)
(5.056112e+00, 1.00275)
(5.110220e+00, 1.00275)
(5.164329e+00, 1.00275)
(5.218437e+00, 1.00275)
(5.272545e+00, 1.00275)
(5.326653e+00, 1.00275)
(5.380762e+00, 1.00275)
(5.434870e+00, 1.00275)
(5.488978e+00, 1.00275)
(5.543086e+00, 1.00275)
(5.597194e+00, 1.00275)
(5.651303e+00, 1.00275)
(5.705411e+00, 1.00275)
(5.759519e+00, 1.00275)
(5.813627e+00, 1.00275)
(5.867735e+00, 1.00275)
(5.921844e+00, 1.00275)
(5.975952e+00, 1.00275)
(6.030060e+00, 1.00275)
(6.084168e+00, 1.00275)
(6.138277e+00, 1.00275)
(6.192385e+00, 1.00275)
(6.246493e+00, 1.00275)
(6.300601e+00, 1.00275)
(6.354709e+00, 1.00275)
(6.408818e+00, 1.00275)
(6.462926e+00, 1.00275)
(6.517034e+00, 1.00275)
(6.571142e+00, 1.00275)
(6.625251e+00, 1.00275)
(6.679359e+00, 1.00275)
(6.733467e+00, 1.00275)
(6.787575e+00, 1.00275)
(6.841683e+00, 1.00275)
(6.895792e+00, 1.00275)
(6.949900e+00, 1.00275)
(7.004008e+00, 1.00275)
(7.058116e+00, 1.00275)
(7.112224e+00, 1.00275)
(7.166333e+00, 1.00048)
(7.220441e+00, 1.00048)
(7.274549e+00, 1.00048)
(7.328657e+00, 1.00048)
(7.382766e+00, 1.00048)
(7.436874e+00, 1.00048)
(7.490982e+00, 1.00048)
(7.545090e+00, 1.00048)
(7.599198e+00, 1.00048)
(7.653307e+00, 1.00048)
(7.707415e+00, 1.00048)
(7.761523e+00, 1.00048)
(7.815631e+00, 1.00048)
(7.869739e+00, 1.00048)
(7.923848e+00, 1.00048)
(7.977956e+00, 1.00048)
(8.032064e+00, 1.00048)
(8.086172e+00, 1.00048)
(8.140281e+00, 1.00048)
(8.194389e+00, 1.00048)
(8.248497e+00, 1.00048)
(8.302605e+00, 1.00048)
(8.356713e+00, 1.00048)
(8.410822e+00, 1.00048)
(8.464930e+00, 1.00048)
(8.519038e+00, 1.00048)
(8.573146e+00, 1.00048)
(8.627255e+00, 1.00048)
(8.681363e+00, 1.00048)
(8.735471e+00, 1.00048)
(8.789579e+00, 1.00048)
(8.843687e+00, 1.00048)
(8.897796e+00, 1.00048)
(8.951904e+00, 1.00048)
(9.006012e+00, 1.00048)
(9.060120e+00, 1.00048)
(9.114228e+00, 1.00048)
(9.168337e+00, 1.00048)
(9.222445e+00, 1.00048)
(9.276553e+00, 1.00048)
(9.330661e+00, 1.00048)
(9.384770e+00, 1.00048)
(9.438878e+00, 1.00048)
(9.492986e+00, 1.00048)
(9.547094e+00, 1.00048)
(9.601202e+00, 1.00048)
(9.655311e+00, 1.00048)
(9.709419e+00, 1.00048)
(9.763527e+00, 1.00048)
(9.817635e+00, 1.00048)
(9.871743e+00, 1.00048)
(9.925852e+00, 1.00048)
(9.979960e+00, 1.00048)
(1.003407e+01, 1.00048)
(1.008818e+01, 1.00048)
(1.014228e+01, 1.00048)
(1.019639e+01, 1.00048)
(1.025050e+01, 1.00048)
(1.030461e+01, 1.00048)
(1.035872e+01, 1.00048)
(1.041283e+01, 1.00048)
(1.046693e+01, 1.00048)
(1.052104e+01, 1.00048)
(1.057515e+01, 1.00048)
(1.062926e+01, 1.00048)
(1.068337e+01, 1.00048)
(1.073747e+01, 1.00048)
(1.079158e+01, 1.00048)
(1.084569e+01, 1.00048)
(1.089980e+01, 1.00048)
(1.095391e+01, 1.00048)
(1.100802e+01, 1.00165)
(1.106212e+01, 1.00165)
(1.111623e+01, 1.00165)
(1.117034e+01, 1.00165)
(1.122445e+01, 1.00165)
(1.127856e+01, 1.00165)
(1.133267e+01, 1.00165)
(1.138677e+01, 1.00165)
(1.144088e+01, 1.00165)
(1.149499e+01, 1.00165)
(1.154910e+01, 1.00165)
(1.160321e+01, 1.00165)
(1.165731e+01, 1.00165)
(1.171142e+01, 1.00165)
(1.176553e+01, 1.00165)
(1.181964e+01, 1.00165)
(1.187375e+01, 1.00165)
(1.192786e+01, 1.00165)
(1.198196e+01, 1.00165)
(1.203607e+01, 1.00165)
(1.209018e+01, 1.00165)
(1.214429e+01, 1.00165)
(1.219840e+01, 1.00165)
(1.225251e+01, 1.00165)
(1.230661e+01, 1.00165)
(1.236072e+01, 1.00165)
(1.241483e+01, 1.00165)
(1.246894e+01, 1.00165)
(1.252305e+01, 1.00165)
(1.257715e+01, 1.00165)
(1.263126e+01, 1.00165)
(1.268537e+01, 1.00165)
(1.273948e+01, 1.00165)
(1.279359e+01, 1.00165)
(1.284770e+01, 1.00165)
(1.290180e+01, 1.00165)
(1.295591e+01, 1.00165)
(1.301002e+01, 1.00165)
(1.306413e+01, 1.00165)
(1.311824e+01, 1.00165)
(1.317234e+01, 1.00165)
(1.322645e+01, 1.00165)
(1.328056e+01, 1.00165)
(1.333467e+01, 1.00165)
(1.338878e+01, 1.00165)
(1.344289e+01, 1.00165)
(1.349699e+01, 1.00165)
(1.355110e+01, 1.00165)
(1.360521e+01, 1.00165)
(1.365932e+01, 1.00165)
(1.371343e+01, 1.00165)
(1.376754e+01, 1.00165)
(1.382164e+01, 1.00165)
(1.387575e+01, 1.00165)
(1.392986e+01, 1.00165)
(1.398397e+01, 1.00165)
(1.403808e+01, 1.00165)
(1.409218e+01, 1.00165)
(1.414629e+01, 1.00165)
(1.420040e+01, 1.00165)
(1.425451e+01, 1.00165)
(1.430862e+01, 1.00165)
(1.436273e+01, 1.00165)
(1.441683e+01, 1.00165)
(1.447094e+01, 1.00165)
(1.452505e+01, 1.00165)
(1.457916e+01, 1.00165)
(1.463327e+01, 1.00165)
(1.468737e+01, 1.00165)
(1.474148e+01, 1.00165)
(1.479559e+01, 1.00165)
(1.484970e+01, 1.00165)
(1.490381e+01, 1.00165)
(1.495792e+01, 1.00165)
(1.501202e+01, 1.00165)
(1.506613e+01, 1.00165)
(1.512024e+01, 1.00165)
(1.517435e+01, 1.00165)
(1.522846e+01, 1.00165)
(1.528257e+01, 1.00165)
(1.533667e+01, 1.00165)
(1.539078e+01, 1.00165)
(1.544489e+01, 1.00165)
(1.549900e+01, 1.00165)
(1.555311e+01, 1.00165)
(1.560721e+01, 1.00165)
(1.566132e+01, 1.00165)
(1.571543e+01, 1.00165)
(1.576954e+01, 1.00165)
(1.582365e+01, 1.00165)
(1.587776e+01, 1.00165)
(1.593186e+01, 1.00165)
(1.598597e+01, 1.00165)
(1.604008e+01, 1.00165)
(1.609419e+01, 1.00165)
(1.614830e+01, 1.00165)
(1.620240e+01, 1.00165)
(1.625651e+01, 1.00165)
(1.631062e+01, 1.00165)
(1.636473e+01, 1.00165)
(1.641884e+01, 1.00165)
(1.647295e+01, 0.999331)
(1.652705e+01, 0.999331)
(1.658116e+01, 0.999331)
(1.663527e+01, 0.999331)
(1.668938e+01, 0.999331)
(1.674349e+01, 0.999331)
(1.679760e+01, 0.999331)
(1.685170e+01, 0.999331)
(1.690581e+01, 0.999331)
(1.695992e+01, 0.999331)
(1.701403e+01, 0.999331)
(1.706814e+01, 0.999331)
(1.712224e+01, 0.999331)
(1.717635e+01, 0.999331)
(1.723046e+01, 0.999331)
(1.728457e+01, 0.999331)
(1.733868e+01, 0.999331)
(1.739279e+01, 0.999331)
(1.744689e+01, 0.999331)
(1.750100e+01, 0.999331)
(1.755511e+01, 0.999331)
(1.760922e+01, 0.999331)
(1.766333e+01, 0.999331)
(1.771743e+01, 0.999331)
(1.777154e+01, 0.999331)
(1.782565e+01, 0.999331)
(1.787976e+01, 0.999331)
(1.793387e+01, 0.999331)
(1.798798e+01, 0.999331)
(1.804208e+01, 0.999331)
(1.809619e+01, 0.999331)
(1.815030e+01, 0.999331)
(1.820441e+01, 0.999331)
(1.825852e+01, 0.999331)
(1.831263e+01, 0.999331)
(1.836673e+01, 0.999331)
(1.842084e+01, 0.999331)
(1.847495e+01, 0.999331)
(1.852906e+01, 0.999331)
(1.858317e+01, 0.999331)
(1.863727e+01, 0.999331)
(1.869138e+01, 0.999331)
(1.874549e+01, 0.999331)
(1.879960e+01, 0.999331)
(1.885371e+01, 0.999331)
(1.890782e+01, 0.999331)
(1.896192e+01, 0.999331)
(1.901603e+01, 0.999331)
(1.907014e+01, 0.999331)
(1.912425e+01, 0.999331)
(1.917836e+01, 0.999331)
(1.923246e+01, 0.999331)
(1.928657e+01, 0.999331)
(1.934068e+01, 0.999331)
(1.939479e+01, 0.999331)
(1.944890e+01, 0.999331)
(1.950301e+01, 0.999331)
(1.955711e+01, 0.999331)
(1.961122e+01, 0.999331)
(1.966533e+01, 0.999331)
(1.971944e+01, 0.999331)
(1.977355e+01, 0.999331)
(1.982766e+01, 0.999331)
(1.988176e+01, 0.999331)
(1.993587e+01, 0.999331)
(1.998998e+01, 0.999331)
(2.004409e+01, 0.999331)
(2.009820e+01, 0.999331)
(2.015230e+01, 0.999331)
(2.020641e+01, 0.999331)
(2.026052e+01, 0.999331)
(2.031463e+01, 0.999331)
(2.036874e+01, 0.999331)
(2.042285e+01, 0.999331)
(2.047695e+01, 0.999331)
(2.053106e+01, 0.999331)
(2.058517e+01, 0.999331)
(2.063928e+01, 0.999331)
(2.069339e+01, 0.999331)
(2.074749e+01, 0.999331)
(2.080160e+01, 0.999331)
(2.085571e+01, 0.999331)
(2.090982e+01, 0.999331)
(2.096393e+01, 0.999331)
(2.101804e+01, 0.999331)
(2.107214e+01, 0.999331)
(2.112625e+01, 0.999331)
(2.118036e+01, 0.999331)
(2.123447e+01, 0.999331)
(2.128858e+01, 0.999331)
(2.134269e+01, 0.999331)
(2.139679e+01, 0.999331)
(2.145090e+01, 0.999331)
(2.150501e+01, 0.999331)
(2.155912e+01, 0.999331)
(2.161323e+01, 0.999331)
(2.166733e+01, 0.999331)
(2.172144e+01, 0.999331)
(2.177555e+01, 0.999331)
(2.182966e+01, 0.999331)
(2.188377e+01, 0.999331)
(2.193788e+01, 0.999331)
(2.199198e+01, 0.999331)
(2.204609e+01, 0.999331)
(2.210020e+01, 0.999331)
(2.215431e+01, 0.999331)
(2.220842e+01, 0.999331)
(2.226253e+01, 0.999331)
(2.231663e+01, 0.999331)
(2.237074e+01, 0.999331)
(2.242485e+01, 0.999331)
(2.247896e+01, 0.999331)
(2.253307e+01, 0.999331)
(2.258717e+01, 0.999331)
(2.264128e+01, 0.999331)
(2.269539e+01, 0.999331)
(2.274950e+01, 0.999331)
(2.280361e+01, 0.999331)
(2.285772e+01, 0.999331)
(2.291182e+01, 0.999331)
(2.296593e+01, 0.999331)
(2.302004e+01, 0.999331)
(2.307415e+01, 0.999331)
(2.312826e+01, 0.999331)
(2.318236e+01, 0.999331)
(2.323647e+01, 0.999331)
(2.329058e+01, 0.999331)
(2.334469e+01, 0.999331)
(2.339880e+01, 0.999331)
(2.345291e+01, 0.999331)
(2.350701e+01, 0.999331)
(2.356112e+01, 0.999331)
(2.361523e+01, 0.999331)
(2.366934e+01, 0.999331)
(2.372345e+01, 0.999331)
(2.377756e+01, 0.999331)
(2.383166e+01, 0.999331)
(2.388577e+01, 0.999331)
(2.393988e+01, 0.999331)
(2.399399e+01, 0.999331)
(2.404810e+01, 0.999331)
(2.410220e+01, 0.999331)
(2.415631e+01, 0.999331)
(2.421042e+01, 0.999331)
(2.426453e+01, 0.999331)
(2.431864e+01, 0.999331)
(2.437275e+01, 0.999331)
(2.442685e+01, 0.999331)
(2.448096e+01, 0.999331)
(2.453507e+01, 1.00073)
(2.458918e+01, 1.00073)
(2.464329e+01, 1.00073)
(2.469739e+01, 1.00073)
(2.475150e+01, 1.00073)
(2.480561e+01, 1.00073)
(2.485972e+01, 1.00073)
(2.491383e+01, 1.00073)
(2.496794e+01, 1.00073)
(2.502204e+01, 1.00073)
(2.507615e+01, 1.00073)
(2.513026e+01, 1.00073)
(2.518437e+01, 1.00073)
(2.523848e+01, 1.00073)
(2.529259e+01, 1.00073)
(2.534669e+01, 1.00073)
(2.540080e+01, 1.00073)
(2.545491e+01, 1.00073)
(2.550902e+01, 1.00073)
(2.556313e+01, 1.00073)
(2.561723e+01, 1.00073)
(2.567134e+01, 1.00073)
(2.572545e+01, 1.00073)
(2.577956e+01, 1.00073)
(2.583367e+01, 1.00073)
(2.588778e+01, 1.00073)
(2.594188e+01, 1.00073)
(2.599599e+01, 1.00073)
(2.605010e+01, 1.00073)
(2.610421e+01, 1.00073)
(2.615832e+01, 1.00073)
(2.621242e+01, 1.00073)
(2.626653e+01, 1.00073)
(2.632064e+01, 1.00073)
(2.637475e+01, 1.00073)
(2.642886e+01, 1.00073)
(2.648297e+01, 1.00073)
(2.653707e+01, 1.00073)
(2.659118e+01, 1.00073)
(2.664529e+01, 1.00073)
(2.669940e+01, 1.00073)
(2.675351e+01, 1.00073)
(2.680762e+01, 1.00073)
(2.686172e+01, 1.00073)
(2.691583e+01, 1.00073)
(2.696994e+01, 1.00073)
(2.702405e+01, 1.00073)
(2.707816e+01, 1.00073)
(2.713226e+01, 1.00073)
(2.718637e+01, 1.00073)
(2.724048e+01, 1.00073)
(2.729459e+01, 1.00073)
(2.734870e+01, 1.00073)
(2.740281e+01, 1.00073)
(2.745691e+01, 1.00073)
(2.751102e+01, 1.00073)
(2.756513e+01, 1.00073)
(2.761924e+01, 1.00073)
(2.767335e+01, 1.00073)
(2.772745e+01, 1.00073)
(2.778156e+01, 1.00073)
(2.783567e+01, 1.00073)
(2.788978e+01, 1.00073)
(2.794389e+01, 1.00073)
(2.799800e+01, 1.00073)
(2.805210e+01, 1.00073)
(2.810621e+01, 1.00073)
(2.816032e+01, 1.00073)
(2.821443e+01, 1.00073)
(2.826854e+01, 1.00073)
(2.832265e+01, 1.00073)
(2.837675e+01, 1.00073)
(2.843086e+01, 1.00073)
(2.848497e+01, 1.00073)
(2.853908e+01, 1.00073)
(2.859319e+01, 1.00073)
(2.864729e+01, 1.00073)
(2.870140e+01, 1.00073)
(2.875551e+01, 1.00073)
(2.880962e+01, 1.00073)
(2.886373e+01, 1.00073)
(2.891784e+01, 1.00073)
(2.897194e+01, 1.00073)
(2.902605e+01, 1.00073)
(2.908016e+01, 1.00073)
(2.913427e+01, 1.00073)
(2.918838e+01, 1.00073)
(2.924248e+01, 1.00073)
(2.929659e+01, 1.00073)
(2.935070e+01, 1.00073)
(2.940481e+01, 1.00073)
(2.945892e+01, 1.00073)
(2.951303e+01, 1.00073)
(2.956713e+01, 1.00073)
(2.962124e+01, 1.00073)
(2.967535e+01, 1.00073)
(2.972946e+01, 1.00073)
(2.978357e+01, 1.00073)
(2.983768e+01, 1.00073)
(2.989178e+01, 1.00073)
(2.994589e+01, 1.00073)
(30, 1.00073)
}; 

%% file: 4-3.tex
\subsection{{Test problem 3: h}ighly variable piecewise constant coefficients}
\label{4-3}

\begin{figure}[h]
\begin{minipage}[t]{0.45\linewidth}
\centering
\includegraphics[width=5.7cm ]{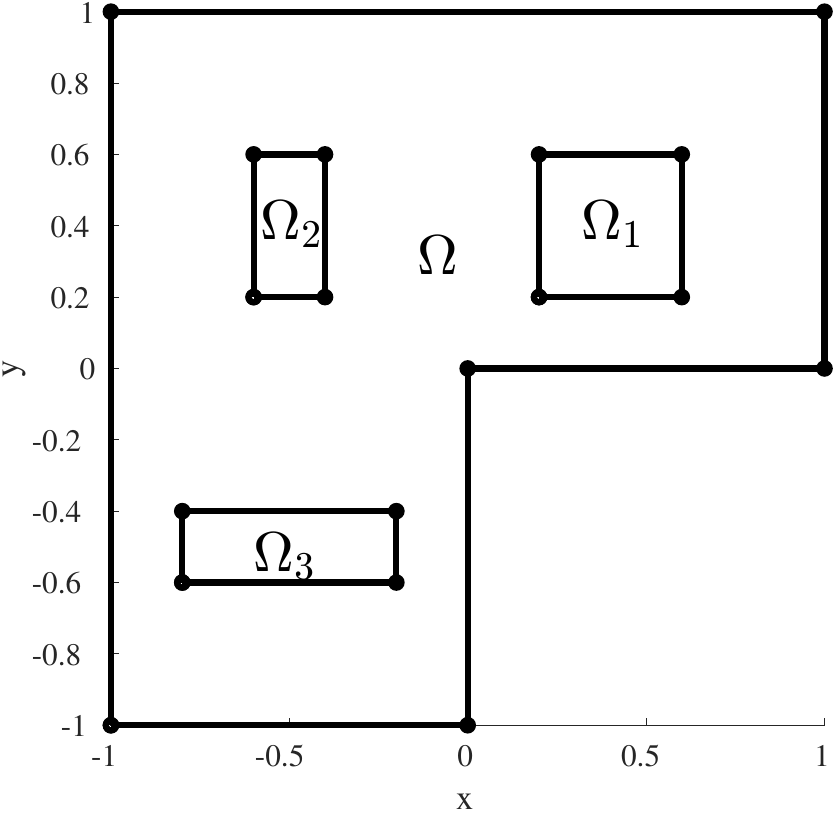}
\caption{Geometry of the domain $\Omega$ in test problem 3.}
\label{fig:D}
\end{minipage}
\hspace{0.5cm}
\begin{minipage}[t]{0.45\linewidth}
\centering
\includegraphics[width=5.7cm]{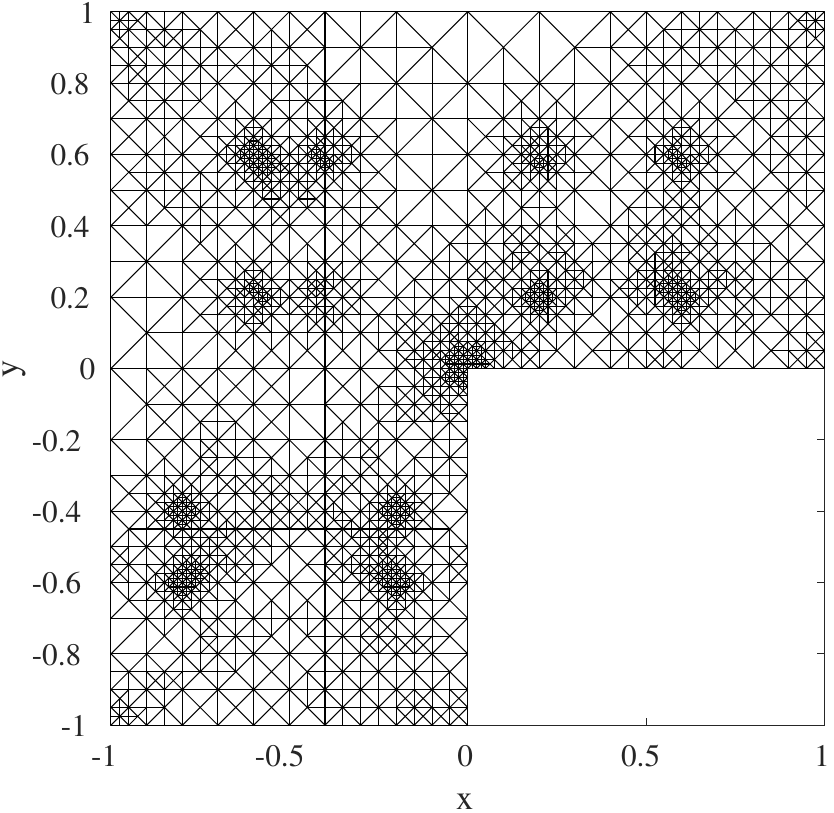}
\caption{Mesh with 3733 elements for  \cref{4-3-1} in test problem 3.}
\label{fig:mesh}
\end{minipage}
\end{figure}

In order to \deleted{demonstrate}\added{explore} the effectiveness of the stopping criteria with variable diffusion coefficient $\kappa(x)$, we consider two problems on \added{an}\deleted{a} L-shape domain with highly variable piecewise constant coefficients and the homogeneous Dirichlet boundary condition, similar to the example considered in \cite[Section 4.1]{arioli2004stopping}. As shown in \cref{fig:D}, the domain $\Omega$ is partitioned into four subdomains, and $\kappa({x})$ is constant on each subdomain. 

\begin{exmp}
\label{4-3-1}
   We choose $f_1(x) = 0.1$, and 
\[
\kappa_1({x})=\left\{\begin{aligned}{}
1, &\quad {x} \in \Omega \backslash \left\{\Omega_1\cup \Omega_2 \cup \Omega_3 \right\} \\
10^{-6}, &\quad {x} \in \Omega_1 \cup \Omega_2 \cup \Omega_3.
\end{aligned}\right.
\]
\end{exmp}

\begin{exmp}
     \label{4-3-2}
   Let $f_2(x) = 10$, and
 \[
\kappa_2({x})=\left\{\begin{aligned}{}
1, &\quad {x} \in \Omega \backslash \left\{\Omega_1\cup \Omega_2 \cup \Omega_3 \right\} \\
10^{6}, &\quad {x} \in \Omega_1 \cup \Omega_2 \cup \Omega_3.
\end{aligned}\right.
\]
\end{exmp}

We begin with a structured mesh  of $\Omega$ consisting of 150 isosceles right triangle elements, and refine the mesh adaptively, using $\eta_{\text{R},K}$ \eqref{eq:etaRR} as an error indicator. The adaptive mesh refinement strategy is to refine all elements where $\eta_{\text{R},K}$ is greater than the average $\eta_{\text{R},K}$. 
The mesh, as illustrated in \cref{fig:mesh}, consists of 3733 elements and is used in \cref{4-3-1}. The refinement is concentrated near corners of $\Omega_1$, $\Omega_2$, $\Omega_3$ and the reentrant corner of the L-shape domain. 
We assume that the solution obtained from a mesh refined six times is an accurate approximation to the exact solution of the continuous problem for \cref{4-3-1}.
\deleted{
In this example, we do not examine $\eta_{\text{MR}}$ as the mesh is regular and the performance of $\eta_{\text{MR}}$ is similar to $\eta_{\text{R}}$.}

\begin{table}[htbp]
\centering
\caption{Numbers of iterations (iter) and quality ratios (qual. \eqref{ev}) resulting from applying stopping criteria to the solution in test problem 3 with the highly variable coefficient.}
\begin{tabular}{ c|l|c|c|c|c|c|c} 
\hline
\multirow{ 2}{*}{$\kappa(x), f(x)$} & \multirow{ 2}{*}{Criterion} & \multicolumn{2}{c|}{$N=4$}  & \multicolumn{2}{c|}{$N=6$} & \multicolumn{2}{c}{$N=8$} \\ \cline{3-8}
\multirow{ 2}{*}{} & \multirow{ 2}{*}{} & iter & qual.  & iter & qual.  & iter & qual.   \\
\hline 
\multirow{6}{*}{$\kappa_1(x),f_1(x)$} & $\eta_{\text{alg}} \leq \tau \eta_{\text{R}}$ & 76 & 1.03 &139 & 1.07 &212 & 1.11  \\ 
\multirow{6}{*}{} & $\eta_{\text{alg}} \leq \tau \eta_{\text{FC}}$ &86 & 1.00 &155 & 1.00 &243 & 1.01  \\ 
\multirow{6}{*}{} & $ \|\Br_k\|_{\Bw} \leq \tau \eta_{\text{RF}}^{\Bw}$ &70 & 1.13 &131 & 1.14 &201 & 1.26  \\ 
\multirow{6}{*}{} & $\|\Br_k\|\leq 10^{-8}\|\Br_0\|$ &192 & 1.00 &334 & 1.00 &506 & 1.00  \\ 
\hline
\multirow{6}{*}{$\kappa_2(x),f_2(x)$} & $\eta_{\text{alg}} \leq \tau \eta_{\text{R}}$ & 244 & 32.95 &419 & 53.17 &631 & 93.99  \\ 
\multirow{6}{*}{} & $\eta_{\text{alg}} \leq \tau \eta_{\text{FC}}$ &252 & 32.95 &435 & 53.17 &663 & 93.99  \\ 
\multirow{6}{*}{} & $ \|\Br_k\|_{\Bw} \leq \tau \eta_{\text{RF}}^{\Bw}$ &241 & 32.95 &417 & 53.17 &632 & 93.99  \\ 
\multirow{6}{*}{} & $\|\Br_k\|\leq 10^{-8}\|\Br_0\|$ &634 & 1.00 &1104 & 1.00 &1654 & 1.00  \\ 
\hline
\end{tabular}
\label{tb:3}
\end{table}

\begin{figure}[htbp]
\begin{tikzpicture}
\begin{groupplot}[group style={
                      group name=myplot,
                      group size= 2 by 1,
                       horizontal sep = 50pt,
                        vertical sep= 2cm
            },height=6cm,width=7.1cm]
 \nextgroupplot[  
    xlabel={iteration},
    ymode=log,  
    ymin=1e-6,
    ymax=1e2 ]
 \input{Data/omega_101}   
\nextgroupplot[ 
    xlabel={iteration},
    ymode=log,
    ymin=1e-6,
    ymax=1e2,
]
 \input{Data/omega_102}   
    \end{groupplot}
\path (myplot c1r1.south west|-current bounding box.south)--
      coordinate(legendpos)
      (myplot c2r1.south east|-current bounding box.south);
\matrix[
    matrix of nodes,
    anchor=south,
    draw,
    inner sep=0.2em,
    draw
  ]at([yshift=-9ex]legendpos)
  {
    \ref{plot:totalerr}& total error&[5pt]
    \ref{plot:exactalg}& $\|\Bx_k - \Bx\|_{\BA}$ &[5pt]
    \ref{plot:etarf}& $\eta_{\text{RF}}^{\Bw}$ & [5pt]
    \ref{plot:res}& $ \| \Br_k\|_{\Bw}$  & [5pt]\\
    \ref{plot:etaalg}& $\eta_{\text{alg}}$ & [5pt] 
    \ref{plot:etar}& $\eta_{\text{R}}$ & [5pt]
    \ref{plot:bdm}& \added{$\eta_{\text{FC}}$}\deleted{$\eta_{\text{BDM}}$} \\ };
\end{tikzpicture}
\caption{Convergence history of the Poisson problem with $f_1(x)$ and a highly variable coefficient $\kappa_1(x)$ in \cref{4-3-1} and polynomial degree $N=6$. Left: the total error, the $\BA$\text{-}norm error $\|\Bx_k-\Bx\|_{\BA}$, the weighted norm of the \deleted{linear residual}\added{linear system residual} $\|\Br_k\|_{\Bw}$ and $\eta_{\text{RF}}^{\Bw}$. Right: the total error, the $\BA$\text{-}norm error $\|\Bx_k-\Bx\|_{\BA}$ and its estimator $\eta_{\text{alg}}$(delay parameter $d=10$), and error estimators $\eta_{\text{R}}$ and \added{$\eta_{\text{FC}}$}\deleted{$\eta_{\text{BDM}}$}. }
\label{iteration31_rf}
\end{figure}

\cref{iteration31_rf} exhibits the convergence history of the energy norm of the error and its error estimates in the whole domain $\Omega$ for \cref{4-3-1} with $N=6$.
The total error converges rapidly, and the norm of the residual decreases roughly monotonically. We note that the algebraic estimator $\eta_{\text{alg}}$ provides an accurate approximation to the $\BA$-norm of the algebraic error.
The left part of  \cref{iteration31_rf} demonstrates that the separation of $\eta_{\text{RF}}^{\Bw}$ and $\|\Br_k\|_{\Bw}$ is close to the separation of the total error and the algebraic error, halting the iteration at a reasonable point.
On the right part of \cref{iteration31_rf}, the estimator $\eta_{\text{R}}$ slightly overestimates the total error. The indicator \added{$\eta_{\text{FC}}$}\deleted{$\eta_{\text{BDM}}$} \deleted{and $\eta_{\underline{\text{BDM}}}$ follow} \added{follows} the total error closely. 

\cref{tb:3} displays the results for \cref{4-3-1} with $N=4,6,8$. The criterion $\eta_{\text{alg}}\leq \tau \eta_\text{R}$ exhibits satisfactory performance. 
As for results of \cref{4-1}\deleted{and \mbox{\cref{4-2}}}, applying $\eta_{\text{alg}}\leq \tau \eta_\text{BDM}$ demonstrates a favorable termination\deleted{, while the criterion $\eta_{\text{alg}}\leq \tau \eta_{\underline{\text{BDM}}}$ leads to slightly delayed termination}. 
The criterion $\|\Br_k\|_{\Bw}\leq \tau \eta_{\text{RF}}^{\Bw}$ yields a small quality ratio and requires a small number of iterations.
Similarly, although the empirical criterion $\|\Br_k\|\leq 10^{-8}\|\Br_0\|$ achieves a quality ratio 1, it requires a larger number of additional  iterations, compared with other criteria.

\begin{figure}[htbp]
\begin{tikzpicture}
\begin{groupplot}[group style={
                      group name=myplot,
                      group size= 2 by 1,
                       horizontal sep = 50pt,
                        vertical sep= 2cm
            },height=6cm,width=7.1cm]           
 \nextgroupplot[  
    xlabel={iteration},
    ymode=log,  
    ymin=1e-7,
    ymax=1e1 ]
 \input{Data/omega_121}   
\nextgroupplot[ 
    xlabel={iteration},
    ymode=log,
    ymin=1e-7,
    ymax=1e1,
]
 \input{Data/omega_122}   
    \end{groupplot}
\path (myplot c1r1.south west|-current bounding box.south)--
      coordinate(legendpos)
      (myplot c2r1.south east|-current bounding box.south);
\matrix[
    matrix of nodes,
    anchor=south,
    draw,
    inner sep=0.2em,
    draw
  ]at([yshift=-9ex]legendpos)
{
    \ref{plot:totalerr}& total error&[5pt]
    \ref{plot:exactalg}& $\|\Bx_k - \Bx\|_{\BA}$ &[5pt]
    \ref{plot:etarf}& $\eta_{\text{RF}}^{\Bw}$ & [5pt]
    \ref{plot:res}& $ \| \Br_k\|_{\Bw}$  & [5pt]\\
    \ref{plot:etaalg}& $\eta_{\text{alg}}$ & [5pt] 
    \ref{plot:etar}& $\eta_{\text{R}}$ & [5pt]
    \ref{plot:bdm}& \added{$\eta_{\text{FC}}$}\deleted{$\eta_{\text{BDM}}$} \\ };
\end{tikzpicture}
\caption{Convergence history of the Poisson problem with $f_2(x)$  and a highly variable coefficient $\kappa_2(x)$ in \cref{4-3-2} and polynomial degree $N=6$. Left: the total error, the $\BA$\text{-}norm error $\|\Bx_k-\Bx\|_{\BA}$, the weighted norm of the \deleted{linear residual}\added{linear system residual} $\|\Br_k\|_{\Bw}$ and $\eta_{\text{RF}}^{\Bw}$. Right: the total error, the $\BA$\text{-}norm error $\|\Bx_k-\Bx\|_{\BA}$ and its estimator $\eta_{\text{alg}}$ (delay parameter $d=10$), and error indicators $\eta_{\text{R}}$ and $ \eta_{\text{FC}}$.   }
\label{iteration3_eta}
\end{figure}

 \cref{4-3-2} is more challenging than \cref{4-3-1}. Analysis for similar problems in one dimensional space is provided in \cite{wang2007large}.
 We refine the mesh such that the mesh consists of 5747 elements in solving  \cref{4-3-2}.
\cref{iteration3_eta} depicts the convergence history of the energy norm of errors and their estimates. 
The energy norm of the total error and the $\BA$-norm of the algebraic error display several plateaus in the iteration process. 
Moreover, the norm of the residual is highly oscillatory when the total error is in the first three plateaus.
All error estimators and indicators follow the trend of the residual, rather than the trend of the total error since all indicators are based on local residuals and jump residuals. In particular, the right part of \cref{iteration3_eta} shows that 
with the delay parameter $d=10$, $\eta_{\text{alg}}$ does not provide an accurate approximation of the exact algebraic error. 
As also highlighted in \cite[Section 4.1]{arioli2004stopping}, a large value of $d$ is necessary to obtain an accurate algebraic error estimator. In this example, the estimator $\eta_{\text{alg}}$ with $d$ exceeding 150, may serve as an effective estimator. However, it requires an extra 150 iterations to obtain the estimator. 
Since $\eta_{\text{alg}}$ substantially underestimates the algebraic error at several phases of the iteration, its inferior performance  contributes to the failure of criteria $\eta_{\text{alg}}\leq \tau \eta_{\text{R}}$\deleted{,} \added{and  $\eta_{\text{alg}}\leq \tau \eta_{\text{FC}}$} \deleted{, and $\eta_{\text{alg}}\leq \tau \eta_{\underline{\text{BDM}}}$}. 

The results for \cref{4-3-2} are also presented in \cref{tb:3}. 
The criteria based on $\eta_\text{R}$, $\eta_\text{BDM}$, and $\eta_\text{RF}^{\Bw}$ result in early termination with the same quality ratios, as they all suggests stopping at the second plateau (approximately from 400 to 600 steps).
\deleted{The termination resulting from the criterion 
$\eta_\text{alg} \leq \tau \eta_{\underline{\text{BDM}}}$ is  slightly early, since the iteration stops at the third plateau.}
As expected, the criterion based on the relative norm of the residual requires a large number of iterations.

\begin{figure}[htbp]
    \centering
    \includegraphics[width=6cm]{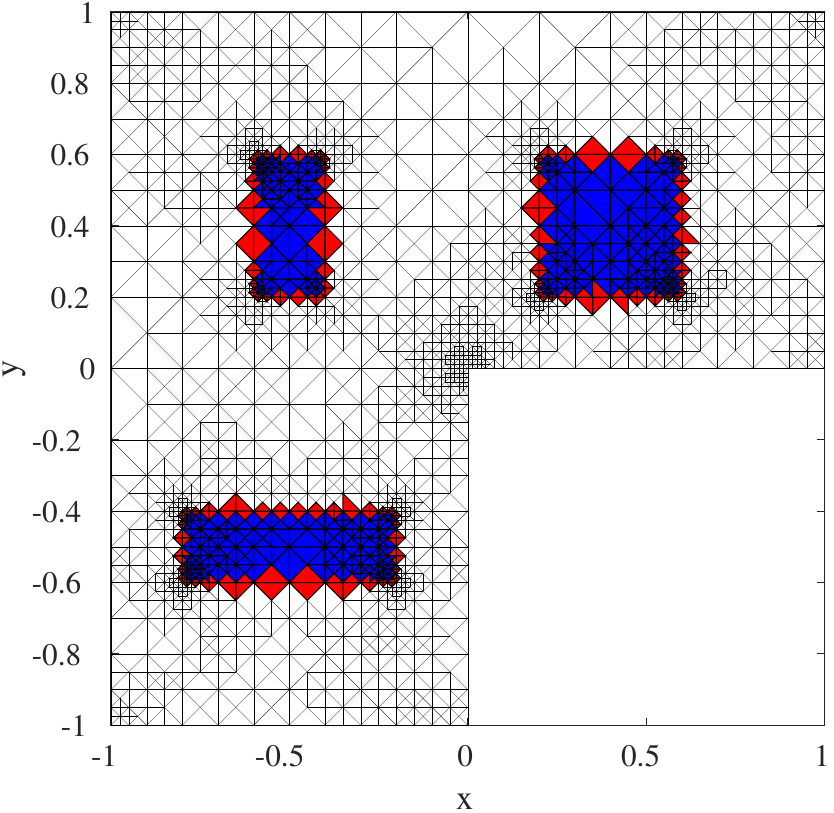}
    \caption{Partition of $\Omega$ in test problem 3: interior subdomain (blue), overlap subdomain (red) and exterior subdomain (white).  }
    \label{subdomains}
\end{figure}

\subsubsection{Subdomain-based stopping criterion}
\label{subsubsetion:4-3-1}
To address the suboptimal performance of $\eta_{\text{RF}}^{\Bw}$ for problems with a highly variable coefficient in \cref{4-3-2}, we consider the subdomain-based stopping criterion in \cref{eq:patchy}.  As illustrated in \cref{subdomains}, we partition the domain $\Omega$ into three subdomains: the interior subdomain colored in blue, the overlap subdomain colored in red and the exterior subdomain in white.
 \cref{tb:3.1} demonstrates that quality ratios of applying  the subdomain-based stopping criterion are one. 
Compared with the results from \cref{tb:3}, the subdomain-based stopping criterion results in late termination for \cref{4-3-1}.
For \cref{4-3-2}, it is the only criterion, in addition to the criterion based on relative residual norm, that leads to reasonable termination. 
Moreover, the subdomain-based stopping criterion requires fewer iterations than the relative residual norm criterion.
It strikes a balance between efficiency and reliability, thereby making it a competitive choice. 
The computation of the subdomain-based stopping criterion is presented in \cref{appen}.

\begin{table}[htbp]
\centering
\caption{Numbers of iterations (iter) and quality ratios (qual. \eqref{ev}) resulting from applying the subdomain-based stopping criterion to the solution in \cref{4-3-1} and \cref{4-3-2}.}
\begin{tabular}{ c|l|c|c|c|c|c|c} 
\hline
\multirow{ 2}{*}{$\kappa(x), f(x)$} & \multirow{ 2}{*}{Criterion} & \multicolumn{2}{c|}{$N=4$}  & \multicolumn{2}{c|}{$N=6$} & \multicolumn{2}{c}{$N=8$} \\ \cline{3-8}
\multirow{ 2}{*}{} & \multirow{ 2}{*}{} & iter & qual.  & iter & qual.  & iter & qual.   \\
\hline 
\multirow{1}{*}{$\kappa_1(x),f_1(x)$} & 
$ \|\Br_k^p\|_{\Bw} \leq \tau \eta_{\text{RF}}^{\Bw,p} , \, \forall p$  &79 & 1.02 &196 & 1.00 &399 & 1.00  \\
\multirow{1}{*}{$\kappa_2(x),f_2(x)$} & 
$ \|\Br_k^p\|_{\Bw} \leq \tau \eta_{\text{RF}}^{\Bw,p} , \, \forall p$  &551 & 1.00 &997 & 1.00 &1505 & 1.00  \\ 
\hline
\end{tabular}
\label{tb:3.1}
\end{table}

\begin{figure}[htbp]
\begin{tikzpicture}
\begin{groupplot}[group style={
                      group name=myplot,
                      group size= 2 by 2,
                       horizontal sep = 50pt,
                        vertical sep= 2cm
                       },height=6cm,width=7.1cm,
                       title style={at={(0.5,0)}, below, yshift=-3.5em}]
 \nextgroupplot[ title={(a): Whole Domain $\Omega$},
    xlabel={iteration},
    ymode=log,  
    ymin=1e-9,
    ymax=1e2 ]
 \input{Data/omega_121}   
\nextgroupplot[ 
    title={(b): Interior subdomain},
    xlabel={iteration},
    ymode=log,
    ymin=1e-9,
    ymax=1e2,
]
\input{Data/int_121}
\nextgroupplot[ 
    title={(c): Overlap subdomain},
    xlabel={iteration},
    ymode=log,
        ymin=1e-9,
    ymax=1e2,
]
\input{Data/inter_121}
\nextgroupplot[ 
    title={(d): Exterior subdomain},
    xlabel={iteration},
    ymode=log,
        ymin=1e-9,
    ymax=1e2,
]
\input{Data/ext_121}
    \end{groupplot}
\path (myplot c1r1.south west|-current bounding box.south)--
      coordinate(legendpos)
      (myplot c2r1.south east|-current bounding box.south);
\matrix[
    matrix of nodes,
    anchor=south,
    draw,
    inner sep=0.2em,
    draw
  ]at([yshift=-7ex]legendpos)
  {
    \ref{plot:totalerr}& total error&[5pt]
    \ref{plot:exactalg}& $\|\Bx_k - \Bx\|_{\BA}$ &[5pt]
    \ref{plot:etarf}& $\eta_{\text{RF}}^{\Bw,p}$ & [5pt]
    \ref{plot:res}& $ \| \Br_k^{p}\|_{\Bw}$ \\};
    
\end{tikzpicture}
\caption{
Convergence history of the Poisson problem with a highly variable coefficient, $\kappa_2(x)$, in \cref{4-3-2} and polynomial degree $N=6$. The total error, the exact $\BA$\text{-}norm error $\|\Bx_k-\Bx\|_{\BA}$, the weighted norm of the subdomain \deleted{linear residual}\added{linear system residual} $\|\Br_k^p\|_{\Bw}$ and subdomain error indicator $\eta_{\text{RF}}^{\Bw,p}$.}
\label{iteration32_rf}
\end{figure}

\cref{iteration32_rf} (b), (c), and (d) display the  convergence
history of the error indicators within subdomains for \cref{4-3-2}. In comparison to the convergence in the interior domain, $\eta_{\text{RF}}^{\Bw,p}$ in the overlap and exterior subdomains deviate from $\|\Br^p_k\|_{\Bw}$ much earlier. 
The early separation of the weighted norm of the residual $\|\Br^p_k\|_{\Bw}$ and the indicator $\eta_{\text{RF}}^{\Bw,p}$ in the exterior subdomain leading to the early termination in the whole domain, since the indicator $\eta_{\text{RF}}^{\Bw,p}$ in the exterior subdomain is dominant in the total $\eta_{\text{RF}}^{\Bw}$.
However, the local indicator $\eta_{\text{RF}}^{\Bw,p}$ in the interior subdomain diverges from the local residual $\|\Br^p\|_{\Bw}$ at around the iteration 950. The subdomain-based stopping criterion ensures that the iteration terminates when the local indicator $\eta_{\text{RF}}^{\Bw,p}$ in the interior subdomain tends to stagnate. 
In solving the Poisson problem \cref{4-3-2},  the subdomain stopping criterion is more conservative and more reliable in practice.
Similarly, the partition can be generalized to the other criteria based on error estimators.

%% file: Data/omega_101.tex
\addplot[ 
line width=2.5pt, color=gray, solid ]
coordinates { 
(1, 336.263)
(2, 288.75)
(3, 245.454)
(4, 215.484)
(5, 187.651)
(6, 163.494)
(7, 141.177)
(8, 121.069)
(9, 103.518)
(10, 88.0071)
(11, 74.4967)
(12, 62.7539)
(13, 52.7192)
(14, 44.3296)
(15, 37.2365)
(16, 31.0042)
(17, 26.0168)
(18, 22.0913)
(19, 19.0399)
(20, 16.5761)
(21, 14.7158)
(22, 13.0272)
(23, 11.2725)
(24, 9.63662)
(25, 8.20765)
(26, 7.12809)
(27, 6.32594)
(28, 5.6758)
(29, 5.16654)
(30, 4.78535)
(31, 4.4549)
(32, 4.18954)
(33, 3.99262)
(34, 3.83172)
(35, 3.70556)
(36, 3.60627)
(37, 3.52047)
(38, 3.44602)
(39, 3.38721)
(40, 3.33774)
(41, 3.29603)
(42, 3.25489)
(43, 3.21182)
(44, 3.16222)
(45, 3.10485)
(46, 3.03499)
(47, 2.9531)
(48, 2.85407)
(49, 2.73833)
(50, 2.60579)
(51, 2.46591)
(52, 2.31711)
(53, 2.15818)
(54, 1.98522)
(55, 1.81581)
(56, 1.66399)
(57, 1.53612)
(58, 1.41528)
(59, 1.29607)
(60, 1.1788)
(61, 1.0713)
(62, 0.978275)
(63, 0.896997)
(64, 0.83163)
(65, 0.787693)
(66, 0.759835)
(67, 0.738187)
(68, 0.718764)
(69, 0.700858)
(70, 0.684121)
(71, 0.669135)
(72, 0.655056)
(73, 0.63968)
(74, 0.620817)
(75, 0.595919)
(76, 0.562247)
(77, 0.520674)
(78, 0.473665)
(79, 0.428132)
(80, 0.386797)
(81, 0.352337)
(82, 0.323934)
(83, 0.29965)
(84, 0.278092)
(85, 0.259848)
(86, 0.244962)
(87, 0.231624)
(88, 0.219497)
(89, 0.208665)
(90, 0.197925)
(91, 0.186449)
(92, 0.175584)
(93, 0.165319)
(94, 0.155259)
(95, 0.145936)
(96, 0.136911)
(97, 0.128206)
(98, 0.12004)
(99, 0.111703)
(100, 0.103068)
(101, 0.0948273)
(102, 0.0869152)
(103, 0.0792339)
(104, 0.0734653)
(105, 0.0687586)
(106, 0.0642094)
(107, 0.059915)
(108, 0.0560148)
(109, 0.0526574)
(110, 0.0493724)
(111, 0.0461244)
(112, 0.0434457)
(113, 0.041297)
(114, 0.0395555)
(115, 0.0381095)
(116, 0.0369261)
(117, 0.0358864)
(118, 0.0349399)
(119, 0.0340505)
(120, 0.0331862)
(121, 0.0323241)
(122, 0.0313606)
(123, 0.0303452)
(124, 0.0292431)
(125, 0.0280451)
(126, 0.0268329)
(127, 0.0256426)
(128, 0.0246831)
(129, 0.0240088)
(130, 0.0235548)
(131, 0.0232291)
(132, 0.0229841)
(133, 0.0227724)
(134, 0.0225854)
(135, 0.0224154)
(136, 0.0222567)
(137, 0.0221024)
(138, 0.0219241)
(139, 0.0217247)
(140, 0.0215204)
(141, 0.0212956)
(142, 0.0210846)
(143, 0.0209187)
(144, 0.0207934)
(145, 0.020704)
(146, 0.0206357)
(147, 0.0205812)
(148, 0.0205391)
(149, 0.0205078)
(150, 0.0204857)
(151, 0.0204694)
(152, 0.0204568)
(153, 0.0204476)
(154, 0.0204394)
(155, 0.0204318)
(156, 0.0204247)
(157, 0.0204175)
(158, 0.0204113)
(159, 0.0204062)
(160, 0.0204014)
(161, 0.0203963)
(162, 0.0203906)
(163, 0.0203847)
(164, 0.0203793)
(165, 0.0203741)
(166, 0.0203698)
(167, 0.0203666)
(168, 0.0203638)
(169, 0.0203609)
(170, 0.0203583)
(171, 0.0203559)
(172, 0.0203538)
(173, 0.0203518)
(174, 0.02035)
(175, 0.0203485)
(176, 0.0203471)
(177, 0.0203461)
(178, 0.0203451)
(179, 0.0203443)
(180, 0.0203436)
(181, 0.0203431)
(182, 0.0203426)
(183, 0.0203422)
(184, 0.0203419)
(185, 0.0203417)
(186, 0.0203415)
(187, 0.0203413)
(188, 0.0203411)
(189, 0.0203409)
(190, 0.0203407)
(191, 0.0203405)
(192, 0.0203403)
(193, 0.0203402)
(194, 0.0203401)
(195, 0.02034)
(196, 0.02034)
(197, 0.0203399)
(198, 0.0203399)
(199, 0.0203398)
(200, 0.0203398)
(201, 0.0203398)
(202, 0.0203397)
(203, 0.0203397)
(204, 0.0203396)
(205, 0.0203396)
(206, 0.0203395)
(207, 0.0203395)
(208, 0.0203394)
(209, 0.0203394)
(210, 0.0203394)
(211, 0.0203394)
(212, 0.0203394)
(213, 0.0203394)
(214, 0.0203394)
(215, 0.0203393)
(216, 0.0203393)
(217, 0.0203393)
(218, 0.0203393)
(219, 0.0203393)
(220, 0.0203393)
(221, 0.0203393)
(222, 0.0203393)
(223, 0.0203393)
(224, 0.0203393)
(225, 0.0203393)
(226, 0.0203393)
(227, 0.0203393)
(228, 0.0203393)
(229, 0.0203393)
(230, 0.0203393)
(231, 0.0203393)
(232, 0.0203393)
(233, 0.0203393)
(234, 0.0203393)
(235, 0.0203393)
(236, 0.0203393)
(237, 0.0203393)
(238, 0.0203393)
(239, 0.0203393)
(240, 0.0203393)
(241, 0.0203393)
(242, 0.0203393)
(243, 0.0203393)
(244, 0.0203393)
(245, 0.0203393)
(246, 0.0203393)
(247, 0.0203393)
(248, 0.0203393)
(249, 0.0203393)
(250, 0.0203393)
(251, 0.0203393)
(252, 0.0203393)
(253, 0.0203393)
(254, 0.0203393)
(255, 0.0203393)
(256, 0.0203393)
(257, 0.0203393)
(258, 0.0203393)
(259, 0.0203393)
(260, 0.0203393)
(261, 0.0203393)
(262, 0.0203393)
(263, 0.0203393)
(264, 0.0203393)
(265, 0.0203393)
(266, 0.0203393)
(267, 0.0203393)
(268, 0.0203393)
(269, 0.0203393)
(270, 0.0203393)
(271, 0.0203393)
(272, 0.0203393)
(273, 0.0203393)
(274, 0.0203393)
(275, 0.0203393)
(276, 0.0203393)
(277, 0.0203393)
(278, 0.0203393)
(279, 0.0203393)
(280, 0.0203393)
(281, 0.0203393)
(282, 0.0203393)
(283, 0.0203393)
(284, 0.0203393)
(285, 0.0203393)
(286, 0.0203393)
(287, 0.0203393)
(288, 0.0203393)
(289, 0.0203393)
(290, 0.0203393)
(291, 0.0203393)
(292, 0.0203393)
(293, 0.0203393)
(294, 0.0203393)
(295, 0.0203393)
(296, 0.0203393)
(297, 0.0203393)
(298, 0.0203393)
(299, 0.0203393)
(300, 0.0203393)
(301, 0.0203393)
(302, 0.0203393)
(303, 0.0203393)
(304, 0.0203393)
(305, 0.0203393)
(306, 0.0203393)
(307, 0.0203393)
(308, 0.0203393)
(309, 0.0203393)
(310, 0.0203393)
(311, 0.0203393)
(312, 0.0203393)
(313, 0.0203393)
(314, 0.0203393)
(315, 0.0203393)
(316, 0.0203393)
(317, 0.0203393)
(318, 0.0203393)
(319, 0.0203393)
(320, 0.0203393)
(321, 0.0203393)
(322, 0.0203393)
(323, 0.0203393)
(324, 0.0203393)
(325, 0.0203393)
(326, 0.0203393)
(327, 0.0203393)
(328, 0.0203393)
(329, 0.0203393)
(330, 0.0203393)
(331, 0.0203393)
(332, 0.0203393)
(333, 0.0203393)
(334, 0.0203393)
};
\addplot[ 
line width=2.5pt, color=mygreen, dashed]
coordinates { 
(1, 336.263)
(2, 288.75)
(3, 245.454)
(4, 215.484)
(5, 187.651)
(6, 163.494)
(7, 141.177)
(8, 121.069)
(9, 103.518)
(10, 88.0071)
(11, 74.4967)
(12, 62.7538)
(13, 52.7192)
(14, 44.3296)
(15, 37.2365)
(16, 31.0042)
(17, 26.0168)
(18, 22.0913)
(19, 19.0399)
(20, 16.576)
(21, 14.7158)
(22, 13.0272)
(23, 11.2725)
(24, 9.6366)
(25, 8.20762)
(26, 7.12807)
(27, 6.3259)
(28, 5.67576)
(29, 5.1665)
(30, 4.78531)
(31, 4.45485)
(32, 4.18949)
(33, 3.99257)
(34, 3.83167)
(35, 3.7055)
(36, 3.60621)
(37, 3.52041)
(38, 3.44596)
(39, 3.38715)
(40, 3.33768)
(41, 3.29597)
(42, 3.25483)
(43, 3.21176)
(44, 3.16215)
(45, 3.10478)
(46, 3.03492)
(47, 2.95303)
(48, 2.854)
(49, 2.73826)
(50, 2.60571)
(51, 2.46583)
(52, 2.31702)
(53, 2.15809)
(54, 1.98511)
(55, 1.8157)
(56, 1.66387)
(57, 1.53599)
(58, 1.41513)
(59, 1.29591)
(60, 1.17863)
(61, 1.07111)
(62, 0.978064)
(63, 0.896767)
(64, 0.831381)
(65, 0.78743)
(66, 0.759563)
(67, 0.737907)
(68, 0.718476)
(69, 0.700562)
(70, 0.683819)
(71, 0.668826)
(72, 0.65474)
(73, 0.639357)
(74, 0.620484)
(75, 0.595571)
(76, 0.561879)
(77, 0.520277)
(78, 0.473228)
(79, 0.427649)
(80, 0.386262)
(81, 0.351749)
(82, 0.323295)
(83, 0.298959)
(84, 0.277347)
(85, 0.259051)
(86, 0.244116)
(87, 0.230729)
(88, 0.218553)
(89, 0.207672)
(90, 0.196877)
(91, 0.185336)
(92, 0.174402)
(93, 0.164063)
(94, 0.153921)
(95, 0.144512)
(96, 0.135392)
(97, 0.126582)
(98, 0.118305)
(99, 0.109836)
(100, 0.101041)
(101, 0.0926203)
(102, 0.0845018)
(103, 0.0765789)
(104, 0.0705937)
(105, 0.0656815)
(106, 0.0609029)
(107, 0.0563571)
(108, 0.0521917)
(109, 0.0485707)
(110, 0.0449883)
(111, 0.0413977)
(112, 0.0383907)
(113, 0.035941)
(114, 0.0339257)
(115, 0.0322281)
(116, 0.0308197)
(117, 0.029566)
(118, 0.0284097)
(119, 0.0273085)
(120, 0.0262229)
(121, 0.025123)
(122, 0.0238705)
(123, 0.0225198)
(124, 0.0210112)
(125, 0.0193092)
(126, 0.0175019)
(127, 0.015616)
(128, 0.0139846)
(129, 0.0127568)
(130, 0.0118804)
(131, 0.0112208)
(132, 0.0107043)
(133, 0.010242)
(134, 0.00981904)
(135, 0.00942142)
(136, 0.00903744)
(137, 0.00865051)
(138, 0.00818408)
(139, 0.00763404)
(140, 0.00703149)
(141, 0.00631021)
(142, 0.00555642)
(143, 0.00488954)
(144, 0.00432186)
(145, 0.003869)
(146, 0.00348522)
(147, 0.00314624)
(148, 0.00285804)
(149, 0.00262393)
(150, 0.00244549)
(151, 0.00230419)
(152, 0.0021898)
(153, 0.00210196)
(154, 0.00202108)
(155, 0.00194225)
(156, 0.00186576)
(157, 0.00178599)
(158, 0.00171303)
(159, 0.0016518)
(160, 0.00159164)
(161, 0.00152475)
(162, 0.00144619)
(163, 0.00136105)
(164, 0.00127639)
(165, 0.00119028)
(166, 0.00111597)
(167, 0.00105531)
(168, 0.000999176)
(169, 0.000939214)
(170, 0.000880789)
(171, 0.000824064)
(172, 0.000769879)
(173, 0.000714727)
(174, 0.000661087)
(175, 0.000612326)
(176, 0.000566638)
(177, 0.000526109)
(178, 0.00048904)
(179, 0.000453505)
(180, 0.000422322)
(181, 0.000392953)
(182, 0.000368186)
(183, 0.000348249)
(184, 0.000330539)
(185, 0.000315213)
(186, 0.000301648)
(187, 0.000288939)
(188, 0.000275567)
(189, 0.000260638)
(190, 0.000243696)
(191, 0.00022606)
(192, 0.00021087)
(193, 0.000197634)
(194, 0.000186252)
(195, 0.000176367)
(196, 0.000169192)
(197, 0.000163605)
(198, 0.000158835)
(199, 0.000154508)
(200, 0.000150369)
(201, 0.000146035)
(202, 0.000140659)
(203, 0.000134153)
(204, 0.000126416)
(205, 0.000117146)
(206, 0.000105869)
(207, 9.52425e-05)
(208, 8.67417e-05)
(209, 7.98838e-05)
(210, 7.46036e-05)
(211, 7.04933e-05)
(212, 6.72069e-05)
(213, 6.44897e-05)
(214, 6.20533e-05)
(215, 5.95154e-05)
(216, 5.68617e-05)
(217, 5.41663e-05)
(218, 5.16475e-05)
(219, 4.87882e-05)
(220, 4.57402e-05)
(221, 4.23887e-05)
(222, 3.88323e-05)
(223, 3.58717e-05)
(224, 3.38752e-05)
(225, 3.23988e-05)
(226, 3.1314e-05)
(227, 3.04129e-05)
(228, 2.95202e-05)
(229, 2.85464e-05)
(230, 2.73881e-05)
(231, 2.59788e-05)
(232, 2.44113e-05)
(233, 2.2727e-05)
(234, 2.09859e-05)
(235, 1.91799e-05)
(236, 1.72589e-05)
(237, 1.54284e-05)
(238, 1.38525e-05)
(239, 1.25019e-05)
(240, 1.1459e-05)
(241, 1.05735e-05)
(242, 9.80191e-06)
(243, 9.00641e-06)
(244, 8.18579e-06)
(245, 7.4353e-06)
(246, 6.72609e-06)
(247, 6.11074e-06)
(248, 5.65087e-06)
(249, 5.31307e-06)
(250, 5.02563e-06)
(251, 4.76187e-06)
(252, 4.53355e-06)
(253, 4.31742e-06)
(254, 4.09786e-06)
(255, 3.84574e-06)
(256, 3.55044e-06)
(257, 3.23953e-06)
(258, 2.97233e-06)
(259, 2.74978e-06)
(260, 2.55277e-06)
(261, 2.38565e-06)
(262, 2.24046e-06)
(263, 2.12608e-06)
(264, 2.0415e-06)
(265, 1.9674e-06)
(266, 1.88624e-06)
(267, 1.79294e-06)
(268, 1.69642e-06)
(269, 1.60993e-06)
(270, 1.52953e-06)
(271, 1.45458e-06)
(272, 1.37378e-06)
(273, 1.28166e-06)
(274, 1.19263e-06)
(275, 1.10782e-06)
(276, 1.03342e-06)
(277, 9.72394e-07)
(278, 9.1905e-07)
(279, 8.66932e-07)
(280, 8.14872e-07)
(281, 7.65434e-07)
(282, 7.18298e-07)
(283, 6.80358e-07)
(284, 6.49568e-07)
(285, 6.22029e-07)
(286, 5.96599e-07)
(287, 5.69595e-07)
(288, 5.41563e-07)
(289, 5.12929e-07)
(290, 4.88214e-07)
(291, 4.65592e-07)
(292, 4.44117e-07)
(293, 4.21995e-07)
(294, 3.97837e-07)
(295, 3.71309e-07)
(296, 3.46121e-07)
(297, 3.25125e-07)
(298, 3.08698e-07)
(299, 2.93119e-07)
(300, 2.78774e-07)
(301, 2.64196e-07)
(302, 2.48917e-07)
(303, 2.33892e-07)
(304, 2.18996e-07)
(305, 2.03042e-07)
(306, 1.89167e-07)
(307, 1.77349e-07)
(308, 1.66772e-07)
(309, 1.57469e-07)
(310, 1.48784e-07)
(311, 1.40776e-07)
(312, 1.33072e-07)
(313, 1.24489e-07)
(314, 1.14688e-07)
(315, 1.05224e-07)
(316, 9.67132e-08)
(317, 8.87175e-08)
(318, 8.17779e-08)
(319, 7.52886e-08)
(320, 6.91697e-08)
(321, 6.39402e-08)
(322, 5.95533e-08)
(323, 5.60605e-08)
(324, 5.30106e-08)
(325, 5.01485e-08)
(326, 4.73895e-08)
(327, 4.44793e-08)
(328, 4.14816e-08)
(329, 3.86441e-08)
(330, 3.55279e-08)
(331, 3.24484e-08)
(332, 2.98465e-08)
(333, 2.7972e-08)
(334, 2.67067e-08)
};
\addplot[ 
line width=1.5pt, color= mycyan, densely dotted]
coordinates { 
(1, 127.309)
(2, 183.038)
(3, 145.597)
(4, 128.757)
(5, 104.689)
(6, 93.1263)
(7, 87.2435)
(8, 78.9936)
(9, 66.0073)
(10, 54.8784)
(11, 44.0456)
(12, 38.7203)
(13, 32.6846)
(14, 29.9402)
(15, 25.0181)
(16, 20.0842)
(17, 16.5076)
(18, 13.015)
(19, 10.7266)
(20, 9.20608)
(21, 7.57431)
(22, 7.31997)
(23, 6.6989)
(24, 6.3438)
(25, 5.2236)
(26, 4.28826)
(27, 3.44207)
(28, 3.04358)
(29, 2.45417)
(30, 2.16607)
(31, 1.90552)
(32, 1.5551)
(33, 1.30508)
(34, 1.17543)
(35, 1.03689)
(36, 0.990588)
(37, 0.951379)
(38, 0.886955)
(39, 0.810494)
(40, 0.715338)
(41, 0.646714)
(42, 0.644404)
(43, 0.652288)
(44, 0.693497)
(45, 0.740985)
(46, 0.776043)
(47, 0.817586)
(48, 0.877236)
(49, 0.910947)
(50, 0.940034)
(51, 0.928213)
(52, 0.904845)
(53, 0.888443)
(54, 0.869874)
(55, 0.804196)
(56, 0.717744)
(57, 0.659586)
(58, 0.652606)
(59, 0.633542)
(60, 0.594201)
(61, 0.526161)
(62, 0.464799)
(63, 0.419115)
(64, 0.34925)
(65, 0.279837)
(66, 0.230379)
(67, 0.219426)
(68, 0.203623)
(69, 0.19584)
(70, 0.181815)
(71, 0.170198)
(72, 0.168317)
(73, 0.178124)
(74, 0.196014)
(75, 0.220923)
(76, 0.245026)
(77, 0.250338)
(78, 0.245641)
(79, 0.225287)
(80, 0.201803)
(81, 0.173747)
(82, 0.155401)
(83, 0.139236)
(84, 0.126545)
(85, 0.109232)
(86, 0.0980988)
(87, 0.0908075)
(88, 0.0846384)
(89, 0.0814407)
(90, 0.083159)
(91, 0.0818212)
(92, 0.0760029)
(93, 0.0711463)
(94, 0.0672931)
(95, 0.061475)
(96, 0.0585268)
(97, 0.0550996)
(98, 0.0524509)
(99, 0.0523222)
(100, 0.0491793)
(101, 0.0452148)
(102, 0.0423321)
(103, 0.0367643)
(104, 0.0301232)
(105, 0.0286096)
(106, 0.0278767)
(107, 0.026735)
(108, 0.0249457)
(109, 0.0225712)
(110, 0.0221648)
(111, 0.0201094)
(112, 0.0173663)
(113, 0.0152719)
(114, 0.0136381)
(115, 0.0124627)
(116, 0.0111461)
(117, 0.0105927)
(118, 0.00991138)
(119, 0.00962238)
(120, 0.00919585)
(121, 0.00932097)
(122, 0.00952886)
(123, 0.00970852)
(124, 0.00993277)
(125, 0.00991503)
(126, 0.00977296)
(127, 0.00893573)
(128, 0.00757339)
(129, 0.00616846)
(130, 0.00518845)
(131, 0.00446405)
(132, 0.00408359)
(133, 0.00389866)
(134, 0.00365228)
(135, 0.00355324)
(136, 0.00336105)
(137, 0.00341411)
(138, 0.00364227)
(139, 0.00367536)
(140, 0.00366044)
(141, 0.00370601)
(142, 0.00328347)
(143, 0.0028821)
(144, 0.00237078)
(145, 0.00202079)
(146, 0.00175975)
(147, 0.00155342)
(148, 0.0013648)
(149, 0.00119196)
(150, 0.00101954)
(151, 0.000911482)
(152, 0.000785498)
(153, 0.000700496)
(154, 0.000688203)
(155, 0.000658486)
(156, 0.000645022)
(157, 0.000652383)
(158, 0.000580015)
(159, 0.00053904)
(160, 0.000553046)
(161, 0.000560655)
(162, 0.000591285)
(163, 0.000558797)
(164, 0.000552363)
(165, 0.000510521)
(166, 0.000450251)
(167, 0.000402108)
(168, 0.000395996)
(169, 0.000391947)
(170, 0.000373679)
(171, 0.000358781)
(172, 0.000342009)
(173, 0.00033339)
(174, 0.000303782)
(175, 0.000286261)
(176, 0.000264717)
(177, 0.000239373)
(178, 0.000229077)
(179, 0.000210519)
(180, 0.000194875)
(181, 0.000178051)
(182, 0.000152889)
(183, 0.000138823)
(184, 0.000126543)
(185, 0.000114957)
(186, 0.000106651)
(187, 0.000102667)
(188, 0.000104124)
(189, 0.000107284)
(190, 0.000109443)
(191, 0.000102081)
(192, 8.90487e-05)
(193, 8.31295e-05)
(194, 7.73353e-05)
(195, 6.76726e-05)
(196, 5.66777e-05)
(197, 5.14112e-05)
(198, 4.7763e-05)
(199, 4.4881e-05)
(200, 4.44697e-05)
(201, 4.66969e-05)
(202, 5.08809e-05)
(203, 5.31067e-05)
(204, 5.55887e-05)
(205, 5.96873e-05)
(206, 5.96531e-05)
(207, 5.18936e-05)
(208, 4.40207e-05)
(209, 3.81531e-05)
(210, 3.26779e-05)
(211, 2.89227e-05)
(212, 2.5381e-05)
(213, 2.33729e-05)
(214, 2.2014e-05)
(215, 2.24717e-05)
(216, 2.19001e-05)
(217, 2.05928e-05)
(218, 2.0334e-05)
(219, 2.06787e-05)
(220, 2.05317e-05)
(221, 2.09829e-05)
(222, 1.9716e-05)
(223, 1.64644e-05)
(224, 1.35221e-05)
(225, 1.1673e-05)
(226, 1.00651e-05)
(227, 9.70063e-06)
(228, 9.69942e-06)
(229, 1.00702e-05)
(230, 1.06229e-05)
(231, 1.10175e-05)
(232, 1.09263e-05)
(233, 1.06238e-05)
(234, 1.02267e-05)
(235, 1.00896e-05)
(236, 9.86993e-06)
(237, 8.75282e-06)
(238, 7.81687e-06)
(239, 6.81609e-06)
(240, 5.88876e-06)
(241, 5.24455e-06)
(242, 4.86555e-06)
(243, 4.79945e-06)
(244, 4.45461e-06)
(245, 4.03367e-06)
(246, 3.73682e-06)
(247, 3.17168e-06)
(248, 2.63643e-06)
(249, 2.28501e-06)
(250, 2.16662e-06)
(251, 1.98155e-06)
(252, 1.8587e-06)
(253, 1.74797e-06)
(254, 1.76022e-06)
(255, 1.83284e-06)
(256, 1.81037e-06)
(257, 1.67994e-06)
(258, 1.42245e-06)
(259, 1.29219e-06)
(260, 1.17628e-06)
(261, 1.06014e-06)
(262, 9.63828e-07)
(263, 8.17091e-07)
(264, 7.20781e-07)
(265, 6.9673e-07)
(266, 7.1518e-07)
(267, 7.57094e-07)
(268, 7.00096e-07)
(269, 6.60955e-07)
(270, 6.25621e-07)
(271, 6.03921e-07)
(272, 6.28863e-07)
(273, 6.09065e-07)
(274, 5.75279e-07)
(275, 5.35464e-07)
(276, 4.7572e-07)
(277, 4.28017e-07)
(278, 4.02757e-07)
(279, 3.87832e-07)
(280, 3.74397e-07)
(281, 3.53082e-07)
(282, 3.187e-07)
(283, 2.70252e-07)
(284, 2.48447e-07)
(285, 2.26925e-07)
(286, 2.22727e-07)
(287, 2.23513e-07)
(288, 2.19531e-07)
(289, 2.08791e-07)
(290, 1.90301e-07)
(291, 1.83256e-07)
(292, 1.77505e-07)
(293, 1.73964e-07)
(294, 1.84735e-07)
(295, 1.80171e-07)
(296, 1.65263e-07)
(297, 1.42406e-07)
(298, 1.31118e-07)
(299, 1.26923e-07)
(300, 1.18496e-07)
(301, 1.23839e-07)
(302, 1.20589e-07)
(303, 1.13642e-07)
(304, 1.11136e-07)
(305, 1.05054e-07)
(306, 9.02687e-08)
(307, 8.37957e-08)
(308, 7.75001e-08)
(309, 7.10205e-08)
(310, 6.70024e-08)
(311, 6.36996e-08)
(312, 6.09531e-08)
(313, 6.4128e-08)
(314, 6.21456e-08)
(315, 5.49016e-08)
(316, 5.23656e-08)
(317, 4.66136e-08)
(318, 4.27379e-08)
(319, 4.05607e-08)
(320, 3.75287e-08)
(321, 3.2316e-08)
(322, 2.88006e-08)
(323, 2.52699e-08)
(324, 2.42168e-08)
(325, 2.28074e-08)
(326, 2.26586e-08)
(327, 2.24435e-08)
(328, 2.10595e-08)
(329, 2.05155e-08)
(330, 2.04106e-08)
(331, 1.8553e-08)
(332, 1.54916e-08)
(333, 1.24435e-08)
(334, 1.02685e-08)
};
\addplot[ 
line width=1.5pt, color= myteal,  solid, mark=+,  mark repeat= 33, mark phase = 5]
coordinates { 
(1, 229.626)
(2, 354.864)
(3, 290.22)
(4, 249.14)
(5, 192.22)
(6, 165.096)
(7, 151.59)
(8, 137.121)
(9, 114.081)
(10, 94.2164)
(11, 74.3617)
(12, 64.3756)
(13, 53.6049)
(14, 49.1642)
(15, 40.9613)
(16, 33.0769)
(17, 27.5003)
(18, 21.7428)
(19, 17.9513)
(20, 15.6305)
(21, 13.0324)
(22, 12.7029)
(23, 11.6229)
(24, 10.9291)
(25, 8.84038)
(26, 7.11867)
(27, 5.68711)
(28, 5.09021)
(29, 4.16233)
(30, 3.699)
(31, 3.23866)
(32, 2.60369)
(33, 2.15376)
(34, 1.93345)
(35, 1.71301)
(36, 1.64956)
(37, 1.59373)
(38, 1.48828)
(39, 1.35415)
(40, 1.18331)
(41, 1.05779)
(42, 1.04865)
(43, 1.05982)
(44, 1.12935)
(45, 1.2129)
(46, 1.26875)
(47, 1.34108)
(48, 1.44108)
(49, 1.50341)
(50, 1.55947)
(51, 1.54262)
(52, 1.4953)
(53, 1.46118)
(54, 1.43124)
(55, 1.3287)
(56, 1.19306)
(57, 1.11266)
(58, 1.11578)
(59, 1.10463)
(60, 1.04215)
(61, 0.927986)
(62, 0.817928)
(63, 0.740989)
(64, 0.627233)
(65, 0.515838)
(66, 0.434888)
(67, 0.416904)
(68, 0.389385)
(69, 0.374978)
(70, 0.351918)
(71, 0.331661)
(72, 0.329117)
(73, 0.34695)
(74, 0.379599)
(75, 0.427164)
(76, 0.474368)
(77, 0.487581)
(78, 0.482133)
(79, 0.4483)
(80, 0.408685)
(81, 0.362869)
(82, 0.332953)
(83, 0.307426)
(84, 0.287096)
(85, 0.258568)
(86, 0.239776)
(87, 0.227612)
(88, 0.217453)
(89, 0.213091)
(90, 0.217542)
(91, 0.216624)
(92, 0.207592)
(93, 0.199604)
(94, 0.193271)
(95, 0.183579)
(96, 0.17978)
(97, 0.17571)
(98, 0.173031)
(99, 0.17475)
(100, 0.17092)
(101, 0.165119)
(102, 0.160925)
(103, 0.151594)
(104, 0.140523)
(105, 0.138118)
(106, 0.137466)
(107, 0.136005)
(108, 0.133362)
(109, 0.129529)
(110, 0.129042)
(111, 0.125751)
(112, 0.121238)
(113, 0.117773)
(114, 0.115095)
(115, 0.113264)
(116, 0.111052)
(117, 0.1102)
(118, 0.109044)
(119, 0.108581)
(120, 0.107848)
(121, 0.108157)
(122, 0.108528)
(123, 0.109112)
(124, 0.109742)
(125, 0.110062)
(126, 0.110018)
(127, 0.108768)
(128, 0.106544)
(129, 0.104201)
(130, 0.102584)
(131, 0.101313)
(132, 0.100722)
(133, 0.100372)
(134, 0.0999834)
(135, 0.099795)
(136, 0.0995125)
(137, 0.0995911)
(138, 0.100044)
(139, 0.100155)
(140, 0.100221)
(141, 0.100391)
(142, 0.0997811)
(143, 0.0991919)
(144, 0.0983855)
(145, 0.0978483)
(146, 0.09744)
(147, 0.0971692)
(148, 0.096861)
(149, 0.096652)
(150, 0.0963897)
(151, 0.0962588)
(152, 0.096095)
(153, 0.0959929)
(154, 0.095995)
(155, 0.0959518)
(156, 0.0959582)
(157, 0.0959712)
(158, 0.0958946)
(159, 0.095874)
(160, 0.0958733)
(161, 0.09591)
(162, 0.0959457)
(163, 0.0959097)
(164, 0.0959296)
(165, 0.095868)
(166, 0.0958345)
(167, 0.0957771)
(168, 0.0957861)
(169, 0.095775)
(170, 0.0957711)
(171, 0.0957554)
(172, 0.0957492)
(173, 0.0957414)
(174, 0.0957246)
(175, 0.0957101)
(176, 0.095704)
(177, 0.095684)
(178, 0.095683)
(179, 0.0956745)
(180, 0.0956628)
(181, 0.0956611)
(182, 0.0956448)
(183, 0.0956451)
(184, 0.0956368)
(185, 0.0956365)
(186, 0.0956304)
(187, 0.0956333)
(188, 0.0956292)
(189, 0.095634)
(190, 0.0956292)
(191, 0.0956318)
(192, 0.0956245)
(193, 0.0956257)
(194, 0.0956234)
(195, 0.0956202)
(196, 0.0956214)
(197, 0.0956171)
(198, 0.0956201)
(199, 0.0956172)
(200, 0.0956188)
(201, 0.0956183)
(202, 0.0956193)
(203, 0.0956199)
(204, 0.0956198)
(205, 0.0956213)
(206, 0.0956212)
(207, 0.0956199)
(208, 0.0956197)
(209, 0.0956187)
(210, 0.0956186)
(211, 0.0956181)
(212, 0.0956177)
(213, 0.0956182)
(214, 0.0956172)
(215, 0.0956179)
(216, 0.0956172)
(217, 0.0956175)
(218, 0.0956169)
(219, 0.0956173)
(220, 0.0956165)
(221, 0.0956172)
(222, 0.095616)
(223, 0.0956168)
(224, 0.0956156)
(225, 0.0956163)
(226, 0.0956158)
(227, 0.0956159)
(228, 0.095616)
(229, 0.0956157)
(230, 0.0956163)
(231, 0.0956158)
(232, 0.0956163)
(233, 0.0956159)
(234, 0.0956166)
(235, 0.0956159)
(236, 0.0956168)
(237, 0.0956161)
(238, 0.0956167)
(239, 0.0956164)
(240, 0.0956167)
(241, 0.0956164)
(242, 0.0956167)
(243, 0.0956165)
(244, 0.0956166)
(245, 0.0956167)
(246, 0.0956166)
(247, 0.0956167)
(248, 0.0956166)
(249, 0.0956167)
(250, 0.0956167)
(251, 0.0956166)
(252, 0.0956166)
(253, 0.0956166)
(254, 0.0956166)
(255, 0.0956166)
(256, 0.0956166)
(257, 0.0956166)
(258, 0.0956166)
(259, 0.0956166)
(260, 0.0956165)
(261, 0.0956166)
(262, 0.0956165)
(263, 0.0956165)
(264, 0.0956165)
(265, 0.0956165)
(266, 0.0956165)
(267, 0.0956165)
(268, 0.0956165)
(269, 0.0956165)
(270, 0.0956165)
(271, 0.0956165)
(272, 0.0956165)
(273, 0.0956165)
(274, 0.0956165)
(275, 0.0956165)
(276, 0.0956165)
(277, 0.0956165)
(278, 0.0956165)
(279, 0.0956165)
(280, 0.0956165)
(281, 0.0956165)
(282, 0.0956165)
(283, 0.0956165)
(284, 0.0956165)
(285, 0.0956165)
(286, 0.0956165)
(287, 0.0956166)
(288, 0.0956165)
(289, 0.0956166)
(290, 0.0956165)
(291, 0.0956165)
(292, 0.0956165)
(293, 0.0956165)
(294, 0.0956165)
(295, 0.0956165)
(296, 0.0956165)
(297, 0.0956165)
(298, 0.0956165)
(299, 0.0956165)
(300, 0.0956165)
(301, 0.0956165)
(302, 0.0956165)
(303, 0.0956165)
(304, 0.0956165)
(305, 0.0956165)
(306, 0.0956165)
(307, 0.0956165)
(308, 0.0956165)
(309, 0.0956166)
(310, 0.0956166)
(311, 0.0956166)
(312, 0.0956166)
(313, 0.0956166)
(314, 0.0956166)
(315, 0.0956166)
(316, 0.0956166)
(317, 0.0956166)
(318, 0.0956166)
(319, 0.0956166)
(320, 0.0956166)
(321, 0.0956166)
(322, 0.0956166)
(323, 0.0956166)
(324, 0.0956166)
(325, 0.0956166)
(326, 0.0956166)
(327, 0.0956166)
(328, 0.0956166)
(329, 0.0956166)
(330, 0.0956166)
(331, 0.0956166)
(332, 0.0956166)
(333, 0.0956166)
(334, 0.0956166)
};

%% file: Data/omega_102.tex
\addplot[ 
line width=2.5pt, color=gray, solid ]
coordinates { 
(1, 336.263)
(2, 288.75)
(3, 245.454)
(4, 215.484)
(5, 187.651)
(6, 163.494)
(7, 141.177)
(8, 121.069)
(9, 103.518)
(10, 88.0071)
(11, 74.4967)
(12, 62.7539)
(13, 52.7192)
(14, 44.3296)
(15, 37.2365)
(16, 31.0042)
(17, 26.0168)
(18, 22.0913)
(19, 19.0399)
(20, 16.5761)
(21, 14.7158)
(22, 13.0272)
(23, 11.2725)
(24, 9.63662)
(25, 8.20765)
(26, 7.12809)
(27, 6.32594)
(28, 5.6758)
(29, 5.16654)
(30, 4.78535)
(31, 4.4549)
(32, 4.18954)
(33, 3.99262)
(34, 3.83172)
(35, 3.70556)
(36, 3.60627)
(37, 3.52047)
(38, 3.44602)
(39, 3.38721)
(40, 3.33774)
(41, 3.29603)
(42, 3.25489)
(43, 3.21182)
(44, 3.16222)
(45, 3.10485)
(46, 3.03499)
(47, 2.9531)
(48, 2.85407)
(49, 2.73833)
(50, 2.60579)
(51, 2.46591)
(52, 2.31711)
(53, 2.15818)
(54, 1.98522)
(55, 1.81581)
(56, 1.66399)
(57, 1.53612)
(58, 1.41528)
(59, 1.29607)
(60, 1.1788)
(61, 1.0713)
(62, 0.978275)
(63, 0.896997)
(64, 0.83163)
(65, 0.787693)
(66, 0.759835)
(67, 0.738187)
(68, 0.718764)
(69, 0.700858)
(70, 0.684121)
(71, 0.669135)
(72, 0.655056)
(73, 0.63968)
(74, 0.620817)
(75, 0.595919)
(76, 0.562247)
(77, 0.520674)
(78, 0.473665)
(79, 0.428132)
(80, 0.386797)
(81, 0.352337)
(82, 0.323934)
(83, 0.29965)
(84, 0.278092)
(85, 0.259848)
(86, 0.244962)
(87, 0.231624)
(88, 0.219497)
(89, 0.208665)
(90, 0.197925)
(91, 0.186449)
(92, 0.175584)
(93, 0.165319)
(94, 0.155259)
(95, 0.145936)
(96, 0.136911)
(97, 0.128206)
(98, 0.12004)
(99, 0.111703)
(100, 0.103068)
(101, 0.0948273)
(102, 0.0869152)
(103, 0.0792339)
(104, 0.0734653)
(105, 0.0687586)
(106, 0.0642094)
(107, 0.059915)
(108, 0.0560148)
(109, 0.0526574)
(110, 0.0493724)
(111, 0.0461244)
(112, 0.0434457)
(113, 0.041297)
(114, 0.0395555)
(115, 0.0381095)
(116, 0.0369261)
(117, 0.0358864)
(118, 0.0349399)
(119, 0.0340505)
(120, 0.0331862)
(121, 0.0323241)
(122, 0.0313606)
(123, 0.0303452)
(124, 0.0292431)
(125, 0.0280451)
(126, 0.0268329)
(127, 0.0256426)
(128, 0.0246831)
(129, 0.0240088)
(130, 0.0235548)
(131, 0.0232291)
(132, 0.0229841)
(133, 0.0227724)
(134, 0.0225854)
(135, 0.0224154)
(136, 0.0222567)
(137, 0.0221024)
(138, 0.0219241)
(139, 0.0217247)
(140, 0.0215204)
(141, 0.0212956)
(142, 0.0210846)
(143, 0.0209187)
(144, 0.0207934)
(145, 0.020704)
(146, 0.0206357)
(147, 0.0205812)
(148, 0.0205391)
(149, 0.0205078)
(150, 0.0204857)
(151, 0.0204694)
(152, 0.0204568)
(153, 0.0204476)
(154, 0.0204394)
(155, 0.0204318)
(156, 0.0204247)
(157, 0.0204175)
(158, 0.0204113)
(159, 0.0204062)
(160, 0.0204014)
(161, 0.0203963)
(162, 0.0203906)
(163, 0.0203847)
(164, 0.0203793)
(165, 0.0203741)
(166, 0.0203698)
(167, 0.0203666)
(168, 0.0203638)
(169, 0.0203609)
(170, 0.0203583)
(171, 0.0203559)
(172, 0.0203538)
(173, 0.0203518)
(174, 0.02035)
(175, 0.0203485)
(176, 0.0203471)
(177, 0.0203461)
(178, 0.0203451)
(179, 0.0203443)
(180, 0.0203436)
(181, 0.0203431)
(182, 0.0203426)
(183, 0.0203422)
(184, 0.0203419)
(185, 0.0203417)
(186, 0.0203415)
(187, 0.0203413)
(188, 0.0203411)
(189, 0.0203409)
(190, 0.0203407)
(191, 0.0203405)
(192, 0.0203403)
(193, 0.0203402)
(194, 0.0203401)
(195, 0.02034)
(196, 0.02034)
(197, 0.0203399)
(198, 0.0203399)
(199, 0.0203398)
(200, 0.0203398)
(201, 0.0203398)
(202, 0.0203397)
(203, 0.0203397)
(204, 0.0203396)
(205, 0.0203396)
(206, 0.0203395)
(207, 0.0203395)
(208, 0.0203394)
(209, 0.0203394)
(210, 0.0203394)
(211, 0.0203394)
(212, 0.0203394)
(213, 0.0203394)
(214, 0.0203394)
(215, 0.0203393)
(216, 0.0203393)
(217, 0.0203393)
(218, 0.0203393)
(219, 0.0203393)
(220, 0.0203393)
(221, 0.0203393)
(222, 0.0203393)
(223, 0.0203393)
(224, 0.0203393)
(225, 0.0203393)
(226, 0.0203393)
(227, 0.0203393)
(228, 0.0203393)
(229, 0.0203393)
(230, 0.0203393)
(231, 0.0203393)
(232, 0.0203393)
(233, 0.0203393)
(234, 0.0203393)
(235, 0.0203393)
(236, 0.0203393)
(237, 0.0203393)
(238, 0.0203393)
(239, 0.0203393)
(240, 0.0203393)
(241, 0.0203393)
(242, 0.0203393)
(243, 0.0203393)
(244, 0.0203393)
(245, 0.0203393)
(246, 0.0203393)
(247, 0.0203393)
(248, 0.0203393)
(249, 0.0203393)
(250, 0.0203393)
(251, 0.0203393)
(252, 0.0203393)
(253, 0.0203393)
(254, 0.0203393)
(255, 0.0203393)
(256, 0.0203393)
(257, 0.0203393)
(258, 0.0203393)
(259, 0.0203393)
(260, 0.0203393)
(261, 0.0203393)
(262, 0.0203393)
(263, 0.0203393)
(264, 0.0203393)
(265, 0.0203393)
(266, 0.0203393)
(267, 0.0203393)
(268, 0.0203393)
(269, 0.0203393)
(270, 0.0203393)
(271, 0.0203393)
(272, 0.0203393)
(273, 0.0203393)
(274, 0.0203393)
(275, 0.0203393)
(276, 0.0203393)
(277, 0.0203393)
(278, 0.0203393)
(279, 0.0203393)
(280, 0.0203393)
(281, 0.0203393)
(282, 0.0203393)
(283, 0.0203393)
(284, 0.0203393)
(285, 0.0203393)
(286, 0.0203393)
(287, 0.0203393)
(288, 0.0203393)
(289, 0.0203393)
(290, 0.0203393)
(291, 0.0203393)
(292, 0.0203393)
(293, 0.0203393)
(294, 0.0203393)
(295, 0.0203393)
(296, 0.0203393)
(297, 0.0203393)
(298, 0.0203393)
(299, 0.0203393)
(300, 0.0203393)
(301, 0.0203393)
(302, 0.0203393)
(303, 0.0203393)
(304, 0.0203393)
(305, 0.0203393)
(306, 0.0203393)
(307, 0.0203393)
(308, 0.0203393)
(309, 0.0203393)
(310, 0.0203393)
(311, 0.0203393)
(312, 0.0203393)
(313, 0.0203393)
(314, 0.0203393)
(315, 0.0203393)
(316, 0.0203393)
(317, 0.0203393)
(318, 0.0203393)
(319, 0.0203393)
(320, 0.0203393)
(321, 0.0203393)
(322, 0.0203393)
(323, 0.0203393)
(324, 0.0203393)
(325, 0.0203393)
(326, 0.0203393)
(327, 0.0203393)
(328, 0.0203393)
(329, 0.0203393)
(330, 0.0203393)
(331, 0.0203393)
(332, 0.0203393)
(333, 0.0203393)
(334, 0.0203393)
};
\addplot[ 
line width=2.5pt, color=mygreen, dashed]
coordinates { 
(1, 336.263)
(2, 288.75)
(3, 245.454)
(4, 215.484)
(5, 187.651)
(6, 163.494)
(7, 141.177)
(8, 121.069)
(9, 103.518)
(10, 88.0071)
(11, 74.4967)
(12, 62.7538)
(13, 52.7192)
(14, 44.3296)
(15, 37.2365)
(16, 31.0042)
(17, 26.0168)
(18, 22.0913)
(19, 19.0399)
(20, 16.576)
(21, 14.7158)
(22, 13.0272)
(23, 11.2725)
(24, 9.6366)
(25, 8.20762)
(26, 7.12807)
(27, 6.3259)
(28, 5.67576)
(29, 5.1665)
(30, 4.78531)
(31, 4.45485)
(32, 4.18949)
(33, 3.99257)
(34, 3.83167)
(35, 3.7055)
(36, 3.60621)
(37, 3.52041)
(38, 3.44596)
(39, 3.38715)
(40, 3.33768)
(41, 3.29597)
(42, 3.25483)
(43, 3.21176)
(44, 3.16215)
(45, 3.10478)
(46, 3.03492)
(47, 2.95303)
(48, 2.854)
(49, 2.73826)
(50, 2.60571)
(51, 2.46583)
(52, 2.31702)
(53, 2.15809)
(54, 1.98511)
(55, 1.8157)
(56, 1.66387)
(57, 1.53599)
(58, 1.41513)
(59, 1.29591)
(60, 1.17863)
(61, 1.07111)
(62, 0.978064)
(63, 0.896767)
(64, 0.831381)
(65, 0.78743)
(66, 0.759563)
(67, 0.737907)
(68, 0.718476)
(69, 0.700562)
(70, 0.683819)
(71, 0.668826)
(72, 0.65474)
(73, 0.639357)
(74, 0.620484)
(75, 0.595571)
(76, 0.561879)
(77, 0.520277)
(78, 0.473228)
(79, 0.427649)
(80, 0.386262)
(81, 0.351749)
(82, 0.323295)
(83, 0.298959)
(84, 0.277347)
(85, 0.259051)
(86, 0.244116)
(87, 0.230729)
(88, 0.218553)
(89, 0.207672)
(90, 0.196877)
(91, 0.185336)
(92, 0.174402)
(93, 0.164063)
(94, 0.153921)
(95, 0.144512)
(96, 0.135392)
(97, 0.126582)
(98, 0.118305)
(99, 0.109836)
(100, 0.101041)
(101, 0.0926203)
(102, 0.0845018)
(103, 0.0765789)
(104, 0.0705937)
(105, 0.0656815)
(106, 0.0609029)
(107, 0.0563571)
(108, 0.0521917)
(109, 0.0485707)
(110, 0.0449883)
(111, 0.0413977)
(112, 0.0383907)
(113, 0.035941)
(114, 0.0339257)
(115, 0.0322281)
(116, 0.0308197)
(117, 0.029566)
(118, 0.0284097)
(119, 0.0273085)
(120, 0.0262229)
(121, 0.025123)
(122, 0.0238705)
(123, 0.0225198)
(124, 0.0210112)
(125, 0.0193092)
(126, 0.0175019)
(127, 0.015616)
(128, 0.0139846)
(129, 0.0127568)
(130, 0.0118804)
(131, 0.0112208)
(132, 0.0107043)
(133, 0.010242)
(134, 0.00981904)
(135, 0.00942142)
(136, 0.00903744)
(137, 0.00865051)
(138, 0.00818408)
(139, 0.00763404)
(140, 0.00703149)
(141, 0.00631021)
(142, 0.00555642)
(143, 0.00488954)
(144, 0.00432186)
(145, 0.003869)
(146, 0.00348522)
(147, 0.00314624)
(148, 0.00285804)
(149, 0.00262393)
(150, 0.00244549)
(151, 0.00230419)
(152, 0.0021898)
(153, 0.00210196)
(154, 0.00202108)
(155, 0.00194225)
(156, 0.00186576)
(157, 0.00178599)
(158, 0.00171303)
(159, 0.0016518)
(160, 0.00159164)
(161, 0.00152475)
(162, 0.00144619)
(163, 0.00136105)
(164, 0.00127639)
(165, 0.00119028)
(166, 0.00111597)
(167, 0.00105531)
(168, 0.000999176)
(169, 0.000939214)
(170, 0.000880789)
(171, 0.000824064)
(172, 0.000769879)
(173, 0.000714727)
(174, 0.000661087)
(175, 0.000612326)
(176, 0.000566638)
(177, 0.000526109)
(178, 0.00048904)
(179, 0.000453505)
(180, 0.000422322)
(181, 0.000392953)
(182, 0.000368186)
(183, 0.000348249)
(184, 0.000330539)
(185, 0.000315213)
(186, 0.000301648)
(187, 0.000288939)
(188, 0.000275567)
(189, 0.000260638)
(190, 0.000243696)
(191, 0.00022606)
(192, 0.00021087)
(193, 0.000197634)
(194, 0.000186252)
(195, 0.000176367)
(196, 0.000169192)
(197, 0.000163605)
(198, 0.000158835)
(199, 0.000154508)
(200, 0.000150369)
(201, 0.000146035)
(202, 0.000140659)
(203, 0.000134153)
(204, 0.000126416)
(205, 0.000117146)
(206, 0.000105869)
(207, 9.52425e-05)
(208, 8.67417e-05)
(209, 7.98838e-05)
(210, 7.46036e-05)
(211, 7.04933e-05)
(212, 6.72069e-05)
(213, 6.44897e-05)
(214, 6.20533e-05)
(215, 5.95154e-05)
(216, 5.68617e-05)
(217, 5.41663e-05)
(218, 5.16475e-05)
(219, 4.87882e-05)
(220, 4.57402e-05)
(221, 4.23887e-05)
(222, 3.88323e-05)
(223, 3.58717e-05)
(224, 3.38752e-05)
(225, 3.23988e-05)
(226, 3.1314e-05)
(227, 3.04129e-05)
(228, 2.95202e-05)
(229, 2.85464e-05)
(230, 2.73881e-05)
(231, 2.59788e-05)
(232, 2.44113e-05)
(233, 2.2727e-05)
(234, 2.09859e-05)
(235, 1.91799e-05)
(236, 1.72589e-05)
(237, 1.54284e-05)
(238, 1.38525e-05)
(239, 1.25019e-05)
(240, 1.1459e-05)
(241, 1.05735e-05)
(242, 9.80191e-06)
(243, 9.00641e-06)
(244, 8.18579e-06)
(245, 7.4353e-06)
(246, 6.72609e-06)
(247, 6.11074e-06)
(248, 5.65087e-06)
(249, 5.31307e-06)
(250, 5.02563e-06)
(251, 4.76187e-06)
(252, 4.53355e-06)
(253, 4.31742e-06)
(254, 4.09786e-06)
(255, 3.84574e-06)
(256, 3.55044e-06)
(257, 3.23953e-06)
(258, 2.97233e-06)
(259, 2.74978e-06)
(260, 2.55277e-06)
(261, 2.38565e-06)
(262, 2.24046e-06)
(263, 2.12608e-06)
(264, 2.0415e-06)
(265, 1.9674e-06)
(266, 1.88624e-06)
(267, 1.79294e-06)
(268, 1.69642e-06)
(269, 1.60993e-06)
(270, 1.52953e-06)
(271, 1.45458e-06)
(272, 1.37378e-06)
(273, 1.28166e-06)
(274, 1.19263e-06)
(275, 1.10782e-06)
(276, 1.03342e-06)
(277, 9.72394e-07)
(278, 9.1905e-07)
(279, 8.66932e-07)
(280, 8.14872e-07)
(281, 7.65434e-07)
(282, 7.18298e-07)
(283, 6.80358e-07)
(284, 6.49568e-07)
(285, 6.22029e-07)
(286, 5.96599e-07)
(287, 5.69595e-07)
(288, 5.41563e-07)
(289, 5.12929e-07)
(290, 4.88214e-07)
(291, 4.65592e-07)
(292, 4.44117e-07)
(293, 4.21995e-07)
(294, 3.97837e-07)
(295, 3.71309e-07)
(296, 3.46121e-07)
(297, 3.25125e-07)
(298, 3.08698e-07)
(299, 2.93119e-07)
(300, 2.78774e-07)
(301, 2.64196e-07)
(302, 2.48917e-07)
(303, 2.33892e-07)
(304, 2.18996e-07)
(305, 2.03042e-07)
(306, 1.89167e-07)
(307, 1.77349e-07)
(308, 1.66772e-07)
(309, 1.57469e-07)
(310, 1.48784e-07)
(311, 1.40776e-07)
(312, 1.33072e-07)
(313, 1.24489e-07)
(314, 1.14688e-07)
(315, 1.05224e-07)
(316, 9.67132e-08)
(317, 8.87175e-08)
(318, 8.17779e-08)
(319, 7.52886e-08)
(320, 6.91697e-08)
(321, 6.39402e-08)
(322, 5.95533e-08)
(323, 5.60605e-08)
(324, 5.30106e-08)
(325, 5.01485e-08)
(326, 4.73895e-08)
(327, 4.44793e-08)
(328, 4.14816e-08)
(329, 3.86441e-08)
(330, 3.55279e-08)
(331, 3.24484e-08)
(332, 2.98465e-08)
(333, 2.7972e-08)
(334, 2.67067e-08)
};
\addplot[ 
line width=1.5pt, color= myolive,  solid, mark=square, mark repeat= 33, mark phase = 8]
coordinates { 
(1, 327.907)
(2, 281.848)
(3, 239.725)
(4, 210.875)
(5, 183.919)
(6, 160.527)
(7, 138.759)
(8, 119.036)
(9, 101.752)
(10, 86.432)
(11, 73.0288)
(12, 61.3868)
(13, 51.5)
(14, 43.2695)
(15, 36.3206)
(16, 30.1736)
(17, 25.236)
(18, 21.3498)
(19, 18.3255)
(20, 15.8703)
(21, 14.0253)
(22, 12.3351)
(23, 10.5417)
(24, 8.84208)
(25, 7.32355)
(26, 6.14854)
(27, 5.25583)
(28, 4.50995)
(29, 3.90127)
(30, 3.42915)
(31, 2.99705)
(32, 2.63779)
(33, 2.37176)
(34, 2.1639)
(35, 2.02265)
(36, 1.94782)
(37, 1.91649)
(38, 1.93115)
(39, 1.99368)
(40, 2.08575)
(41, 2.18703)
(42, 2.2859)
(43, 2.37866)
(44, 2.46141)
(45, 2.51851)
(46, 2.53816)
(47, 2.52212)
(48, 2.47845)
(49, 2.41219)
(50, 2.32391)
(51, 2.22105)
(52, 2.10047)
(53, 1.96295)
(54, 1.80263)
(55, 1.63607)
(56, 1.48038)
(57, 1.34713)
(58, 1.21918)
(59, 1.09023)
(60, 0.959978)
(61, 0.836623)
(62, 0.726584)
(63, 0.628819)
(64, 0.553348)
(65, 0.515112)
(66, 0.511104)
(67, 0.523277)
(68, 0.540614)
(69, 0.554891)
(70, 0.564279)
(71, 0.568859)
(72, 0.569355)
(73, 0.565156)
(74, 0.555048)
(75, 0.536281)
(76, 0.506078)
(77, 0.466318)
(78, 0.419737)
(79, 0.373839)
(80, 0.332322)
(81, 0.298962)
(82, 0.27222)
(83, 0.24992)
(84, 0.230716)
(85, 0.214997)
(86, 0.20313)
(87, 0.192906)
(88, 0.183765)
(89, 0.176249)
(90, 0.168971)
(91, 0.160533)
(92, 0.152563)
(93, 0.145094)
(94, 0.136778)
(95, 0.128723)
(96, 0.12092)
(97, 0.113344)
(98, 0.10617)
(99, 0.0985129)
(100, 0.0904728)
(101, 0.0828538)
(102, 0.0752776)
(103, 0.0676207)
(104, 0.0619073)
(105, 0.0572312)
(106, 0.0525291)
(107, 0.0479789)
(108, 0.043782)
(109, 0.0401667)
(110, 0.0365555)
(111, 0.032903)
(112, 0.0300674)
(113, 0.028011)
(114, 0.0266361)
(115, 0.0258032)
(116, 0.0253681)
(117, 0.0251056)
(118, 0.0247293)
(119, 0.0241457)
(120, 0.0233773)
(121, 0.0224779)
(122, 0.0213358)
(123, 0.020056)
(124, 0.0185758)
(125, 0.0168547)
(126, 0.014988)
(127, 0.013001)
(128, 0.0113398)
(129, 0.0102205)
(130, 0.00957613)
(131, 0.00927836)
(132, 0.00914922)
(133, 0.00899951)
(134, 0.00881675)
(135, 0.00859034)
(136, 0.00833838)
(137, 0.00805807)
(138, 0.00766882)
(139, 0.00716893)
(140, 0.00659252)
(141, 0.00587447)
(142, 0.00510672)
(143, 0.00441467)
(144, 0.00382017)
(145, 0.00334617)
(146, 0.00294376)
(147, 0.00259019)
(148, 0.00228778)
(149, 0.00203876)
(150, 0.00185664)
(151, 0.00172755)
(152, 0.00164431)
(153, 0.00160181)
(154, 0.00156703)
(155, 0.00153479)
(156, 0.00149521)
(157, 0.00144085)
(158, 0.00139144)
(159, 0.00135879)
(160, 0.00132572)
(161, 0.00128287)
(162, 0.00122423)
(163, 0.00115828)
(164, 0.00109185)
(165, 0.00102069)
(166, 0.000961416)
(167, 0.000914818)
(168, 0.000871316)
(169, 0.000822469)
(170, 0.000772938)
(171, 0.000724341)
(172, 0.000676131)
(173, 0.000624146)
(174, 0.000572521)
(175, 0.000524961)
(176, 0.000479674)
(177, 0.000439665)
(178, 0.000404008)
(179, 0.000371126)
(180, 0.000344917)
(181, 0.000321417)
(182, 0.000301819)
(183, 0.000286737)
(184, 0.000273068)
(185, 0.000261255)
(186, 0.000249731)
(187, 0.000238158)
(188, 0.000225186)
(189, 0.000209903)
(190, 0.000191773)
(191, 0.00017256)
(192, 0.000157103)
(193, 0.000145128)
(194, 0.00013678)
(195, 0.000131841)
(196, 0.000131977)
(197, 0.000133024)
(198, 0.000133058)
(199, 0.000132255)
(200, 0.000130557)
(201, 0.000127894)
(202, 0.000123564)
(203, 0.000117636)
(204, 0.000110139)
(205, 0.000100902)
(206, 8.93023e-05)
(207, 7.83399e-05)
(208, 6.96897e-05)
(209, 6.32545e-05)
(210, 5.89366e-05)
(211, 5.6325e-05)
(212, 5.48527e-05)
(213, 5.35923e-05)
(214, 5.19912e-05)
(215, 4.9924e-05)
(216, 4.74624e-05)
(217, 4.48223e-05)
(218, 4.23795e-05)
(219, 3.9565e-05)
(220, 3.66341e-05)
(221, 3.34948e-05)
(222, 3.01999e-05)
(223, 2.77536e-05)
(224, 2.65918e-05)
(225, 2.61116e-05)
(226, 2.61285e-05)
(227, 2.6209e-05)
(228, 2.60682e-05)
(229, 2.56632e-05)
(230, 2.48757e-05)
(231, 2.37297e-05)
(232, 2.2357e-05)
(233, 2.08663e-05)
(234, 1.93236e-05)
(235, 1.768e-05)
(236, 1.58943e-05)
(237, 1.41666e-05)
(238, 1.26475e-05)
(239, 1.13168e-05)
(240, 1.02981e-05)
(241, 9.44054e-06)
(242, 8.69048e-06)
(243, 7.90413e-06)
(244, 7.08623e-06)
(245, 6.36348e-06)
(246, 5.71268e-06)
(247, 5.18136e-06)
(248, 4.80599e-06)
(249, 4.54614e-06)
(250, 4.32901e-06)
(251, 4.12118e-06)
(252, 3.94125e-06)
(253, 3.75765e-06)
(254, 3.55313e-06)
(255, 3.3044e-06)
(256, 3.00795e-06)
(257, 2.69814e-06)
(258, 2.44068e-06)
(259, 2.22923e-06)
(260, 2.04381e-06)
(261, 1.8909e-06)
(262, 1.76985e-06)
(263, 1.69633e-06)
(264, 1.65691e-06)
(265, 1.62585e-06)
(266, 1.57795e-06)
(267, 1.50635e-06)
(268, 1.42589e-06)
(269, 1.35657e-06)
(270, 1.29439e-06)
(271, 1.2369e-06)
(272, 1.17103e-06)
(273, 1.08617e-06)
(274, 1.00022e-06)
(275, 9.16704e-07)
(276, 8.43822e-07)
(277, 7.88107e-07)
(278, 7.42537e-07)
(279, 6.9891e-07)
(280, 6.52429e-07)
(281, 6.07548e-07)
(282, 5.64547e-07)
(283, 5.33674e-07)
(284, 5.13483e-07)
(285, 4.99049e-07)
(286, 4.85932e-07)
(287, 4.67688e-07)
(288, 4.44968e-07)
(289, 4.20924e-07)
(290, 4.00797e-07)
(291, 3.83374e-07)
(292, 3.67804e-07)
(293, 3.51247e-07)
(294, 3.32137e-07)
(295, 3.10877e-07)
(296, 2.89855e-07)
(297, 2.72495e-07)
(298, 2.59772e-07)
(299, 2.4723e-07)
(300, 2.3575e-07)
(301, 2.23566e-07)
(302, 2.10361e-07)
(303, 1.9801e-07)
(304, 1.86563e-07)
(305, 1.73648e-07)
(306, 1.62576e-07)
(307, 1.53564e-07)
(308, 1.45346e-07)
(309, 1.38304e-07)
(310, 1.31728e-07)
(311, 1.25418e-07)
(312, 1.19002e-07)
(313, 1.11151e-07)
(314, 1.01702e-07)
(315, 9.25054e-08)
(316, 8.4307e-08)
(317, 7.67618e-08)
(318, 7.04762e-08)
(319, 6.46142e-08)
(320, 5.93482e-08)
(321, 5.50948e-08)
(322, 5.15342e-08)
(323, 4.85833e-08)
(324, 4.57916e-08)
(325, 0)
(326, 0)
(327, 0)
(328, 0)
(329, 0)
(330, 0)
(331, 0)
(332, 0)
(333, 0)
(334, 0)
};
\addplot[ 
line width=1.5pt, color= mypurple,  solid, mark=o,  mark repeat= 33, mark phase = 3]
coordinates { 
(1, 197.588)
(2, 320.894)
(3, 260.444)
(4, 217.838)
(5, 167.86)
(6, 145.956)
(7, 136.947)
(8, 129.65)
(9, 110.992)
(10, 92.221)
(11, 72.8225)
(12, 63.9988)
(13, 52.8586)
(14, 47.8684)
(15, 40.3059)
(16, 33.1265)
(17, 27.6418)
(18, 21.5575)
(19, 17.6998)
(20, 15.7871)
(21, 13.3535)
(22, 13.0036)
(23, 11.7044)
(24, 10.8299)
(25, 8.64635)
(26, 6.94714)
(27, 5.54703)
(28, 4.99123)
(29, 4.07058)
(30, 3.60674)
(31, 3.16435)
(32, 2.56805)
(33, 2.13551)
(34, 1.92135)
(35, 1.69916)
(36, 1.62618)
(37, 1.54136)
(38, 1.43098)
(39, 1.31049)
(40, 1.1594)
(41, 1.05302)
(42, 1.05673)
(43, 1.07855)
(44, 1.15906)
(45, 1.23901)
(46, 1.2851)
(47, 1.35193)
(48, 1.45029)
(49, 1.50139)
(50, 1.55441)
(51, 1.53375)
(52, 1.49875)
(53, 1.47497)
(54, 1.45465)
(55, 1.36081)
(56, 1.22476)
(57, 1.13638)
(58, 1.14288)
(59, 1.12099)
(60, 1.0608)
(61, 0.942797)
(62, 0.838266)
(63, 0.766571)
(64, 0.654433)
(65, 0.541215)
(66, 0.460628)
(67, 0.443571)
(68, 0.416913)
(69, 0.407017)
(70, 0.384054)
(71, 0.363964)
(72, 0.359896)
(73, 0.378509)
(74, 0.409987)
(75, 0.456849)
(76, 0.505472)
(77, 0.517367)
(78, 0.512106)
(79, 0.478746)
(80, 0.44196)
(81, 0.396765)
(82, 0.367867)
(83, 0.342821)
(84, 0.324017)
(85, 0.294573)
(86, 0.276105)
(87, 0.263561)
(88, 0.253809)
(89, 0.249666)
(90, 0.25474)
(91, 0.255)
(92, 0.247161)
(93, 0.239672)
(94, 0.232985)
(95, 0.223611)
(96, 0.219424)
(97, 0.21597)
(98, 0.213391)
(99, 0.215277)
(100, 0.211922)
(101, 0.206151)
(102, 0.202353)
(103, 0.193597)
(104, 0.182634)
(105, 0.180467)
(106, 0.179754)
(107, 0.178222)
(108, 0.175676)
(109, 0.171964)
(110, 0.171457)
(111, 0.168035)
(112, 0.163724)
(113, 0.160482)
(114, 0.157989)
(115, 0.156238)
(116, 0.154207)
(117, 0.15332)
(118, 0.152252)
(119, 0.151779)
(120, 0.151111)
(121, 0.151481)
(122, 0.152075)
(123, 0.152634)
(124, 0.153205)
(125, 0.153373)
(126, 0.153286)
(127, 0.151885)
(128, 0.149535)
(129, 0.147275)
(130, 0.145752)
(131, 0.144624)
(132, 0.144082)
(133, 0.143819)
(134, 0.143471)
(135, 0.143378)
(136, 0.143085)
(137, 0.143177)
(138, 0.143589)
(139, 0.143703)
(140, 0.143756)
(141, 0.143944)
(142, 0.143337)
(143, 0.142812)
(144, 0.142125)
(145, 0.141652)
(146, 0.141371)
(147, 0.141159)
(148, 0.140968)
(149, 0.140812)
(150, 0.14067)
(151, 0.140579)
(152, 0.14051)
(153, 0.140456)
(154, 0.140459)
(155, 0.140448)
(156, 0.140452)
(157, 0.140461)
(158, 0.140424)
(159, 0.140422)
(160, 0.140416)
(161, 0.140443)
(162, 0.14046)
(163, 0.140444)
(164, 0.140461)
(165, 0.140437)
(166, 0.140421)
(167, 0.140404)
(168, 0.140409)
(169, 0.140404)
(170, 0.140404)
(171, 0.140398)
(172, 0.140394)
(173, 0.140395)
(174, 0.140385)
(175, 0.14038)
(176, 0.140377)
(177, 0.140369)
(178, 0.140368)
(179, 0.140362)
(180, 0.140359)
(181, 0.140355)
(182, 0.14035)
(183, 0.140349)
(184, 0.140345)
(185, 0.140345)
(186, 0.140342)
(187, 0.140343)
(188, 0.140341)
(189, 0.140341)
(190, 0.14034)
(191, 0.140339)
(192, 0.140337)
(193, 0.140336)
(194, 0.140335)
(195, 0.140334)
(196, 0.140334)
(197, 0.140333)
(198, 0.140334)
(199, 0.140333)
(200, 0.140333)
(201, 0.140333)
(202, 0.140334)
(203, 0.140334)
(204, 0.140334)
(205, 0.140335)
(206, 0.140335)
(207, 0.140335)
(208, 0.140335)
(209, 0.140335)
(210, 0.140334)
(211, 0.140334)
(212, 0.140334)
(213, 0.140334)
(214, 0.140334)
(215, 0.140334)
(216, 0.140334)
(217, 0.140333)
(218, 0.140334)
(219, 0.140333)
(220, 0.140333)
(221, 0.140333)
(222, 0.140332)
(223, 0.140332)
(224, 0.140332)
(225, 0.140332)
(226, 0.140332)
(227, 0.140332)
(228, 0.140332)
(229, 0.140332)
(230, 0.140332)
(231, 0.140332)
(232, 0.140332)
(233, 0.140332)
(234, 0.140332)
(235, 0.140332)
(236, 0.140333)
(237, 0.140333)
(238, 0.140333)
(239, 0.140333)
(240, 0.140333)
(241, 0.140333)
(242, 0.140333)
(243, 0.140333)
(244, 0.140333)
(245, 0.140333)
(246, 0.140333)
(247, 0.140333)
(248, 0.140333)
(249, 0.140333)
(250, 0.140333)
(251, 0.140333)
(252, 0.140333)
(253, 0.140333)
(254, 0.140333)
(255, 0.140333)
(256, 0.140333)
(257, 0.140333)
(258, 0.140333)
(259, 0.140333)
(260, 0.140333)
(261, 0.140333)
(262, 0.140333)
(263, 0.140333)
(264, 0.140333)
(265, 0.140333)
(266, 0.140333)
(267, 0.140333)
(268, 0.140333)
(269, 0.140333)
(270, 0.140333)
(271, 0.140333)
(272, 0.140333)
(273, 0.140333)
(274, 0.140333)
(275, 0.140333)
(276, 0.140333)
(277, 0.140333)
(278, 0.140333)
(279, 0.140333)
(280, 0.140333)
(281, 0.140333)
(282, 0.140333)
(283, 0.140333)
(284, 0.140333)
(285, 0.140333)
(286, 0.140333)
(287, 0.140333)
(288, 0.140333)
(289, 0.140333)
(290, 0.140333)
(291, 0.140333)
(292, 0.140333)
(293, 0.140333)
(294, 0.140333)
(295, 0.140333)
(296, 0.140333)
(297, 0.140333)
(298, 0.140333)
(299, 0.140333)
(300, 0.140333)
(301, 0.140333)
(302, 0.140333)
(303, 0.140333)
(304, 0.140333)
(305, 0.140333)
(306, 0.140333)
(307, 0.140333)
(308, 0.140333)
(309, 0.140333)
(310, 0.140333)
(311, 0.140333)
(312, 0.140333)
(313, 0.140333)
(314, 0.140333)
(315, 0.140333)
(316, 0.140333)
(317, 0.140333)
(318, 0.140333)
(319, 0.140333)
(320, 0.140333)
(321, 0.140333)
(322, 0.140333)
(323, 0.140333)
(324, 0.140333)
(325, 0.140333)
(326, 0.140333)
(327, 0.140333)
(328, 0.140333)
(329, 0.140333)
(330, 0.140333)
(331, 0.140333)
(332, 0.140333)
(333, 0.140333)
(334, 0.140333)
};
\addplot[ 
line width=1.5pt, color= myrose ,   solid, mark=x, mark repeat= 33, mark phase = 24]
coordinates { 
(1, 37.0293)
(2, 50.331)
(3, 38.6401)
(4, 30.4675)
(5, 20.4699)
(6, 16.1861)
(7, 15.2929)
(8, 14.3189)
(9, 12.5554)
(10, 11.2666)
(11, 9.31291)
(12, 8.16425)
(13, 6.97264)
(14, 6.5327)
(15, 5.42646)
(16, 4.34029)
(17, 3.5842)
(18, 2.84266)
(19, 2.38299)
(20, 2.09971)
(21, 1.77172)
(22, 1.7273)
(23, 1.53562)
(24, 1.37813)
(25, 1.05028)
(26, 0.802077)
(27, 0.620418)
(28, 0.562194)
(29, 0.473191)
(30, 0.434542)
(31, 0.387727)
(32, 0.313485)
(33, 0.258676)
(34, 0.23106)
(35, 0.195988)
(36, 0.180634)
(37, 0.167365)
(38, 0.151623)
(39, 0.136469)
(40, 0.120533)
(41, 0.111808)
(42, 0.117752)
(43, 0.123901)
(44, 0.135702)
(45, 0.148761)
(46, 0.156453)
(47, 0.164)
(48, 0.175036)
(49, 0.182091)
(50, 0.187018)
(51, 0.184377)
(52, 0.179397)
(53, 0.180817)
(54, 0.182205)
(55, 0.168741)
(56, 0.149169)
(57, 0.136133)
(58, 0.132684)
(59, 0.130907)
(60, 0.120655)
(61, 0.10688)
(62, 0.0936679)
(63, 0.0838817)
(64, 0.0705882)
(65, 0.059028)
(66, 0.0502373)
(67, 0.048264)
(68, 0.0453142)
(69, 0.0437405)
(70, 0.042035)
(71, 0.040125)
(72, 0.0402393)
(73, 0.0420764)
(74, 0.0455879)
(75, 0.0507382)
(76, 0.0555431)
(77, 0.0565013)
(78, 0.0550109)
(79, 0.0512376)
(80, 0.0471661)
(81, 0.0438284)
(82, 0.0418137)
(83, 0.0402301)
(84, 0.0386557)
(85, 0.0365345)
(86, 0.0352779)
(87, 0.0344908)
(88, 0.0338136)
(89, 0.0336963)
(90, 0.033867)
(91, 0.033899)
(92, 0.0334511)
(93, 0.0330415)
(94, 0.0328504)
(95, 0.0323749)
(96, 0.0323638)
(97, 0.032245)
(98, 0.0321652)
(99, 0.032339)
(100, 0.0322271)
(101, 0.0320639)
(102, 0.0320207)
(103, 0.0317505)
(104, 0.0314614)
(105, 0.0314355)
(106, 0.031389)
(107, 0.0313572)
(108, 0.0313444)
(109, 0.0312448)
(110, 0.031293)
(111, 0.0312549)
(112, 0.0311738)
(113, 0.03119)
(114, 0.0311188)
(115, 0.0311285)
(116, 0.0311205)
(117, 0.0310807)
(118, 0.0311106)
(119, 0.0310826)
(120, 0.0310736)
(121, 0.0311001)
(122, 0.0310644)
(123, 0.0310894)
(124, 0.031094)
(125, 0.0310606)
(126, 0.0310931)
(127, 0.0310511)
(128, 0.0310375)
(129, 0.0310375)
(130, 0.0310101)
(131, 0.0310174)
(132, 0.0310084)
(133, 0.031005)
(134, 0.0310066)
(135, 0.0309999)
(136, 0.0310002)
(137, 0.0309984)
(138, 0.0309963)
(139, 0.0309955)
(140, 0.0309932)
(141, 0.0309915)
(142, 0.030988)
(143, 0.0309847)
(144, 0.0309824)
(145, 0.0309798)
(146, 0.0309801)
(147, 0.0309787)
(148, 0.0309787)
(149, 0.0309793)
(150, 0.0309783)
(151, 0.0309794)
(152, 0.0309795)
(153, 0.0309795)
(154, 0.0309806)
(155, 0.0309807)
(156, 0.0309817)
(157, 0.0309826)
(158, 0.0309831)
(159, 0.0309842)
(160, 0.0309846)
(161, 0.0309859)
(162, 0.0309866)
(163, 0.0309875)
(164, 0.0309889)
(165, 0.0309894)
(166, 0.0309902)
(167, 0.0309908)
(168, 0.0309911)
(169, 0.0309915)
(170, 0.0309916)
(171, 0.0309918)
(172, 0.0309921)
(173, 0.0309917)
(174, 0.030992)
(175, 0.0309918)
(176, 0.0309914)
(177, 0.0309916)
(178, 0.0309912)
(179, 0.030991)
(180, 0.030991)
(181, 0.0309905)
(182, 0.0309906)
(183, 0.0309903)
(184, 0.0309901)
(185, 0.0309901)
(186, 0.0309898)
(187, 0.0309898)
(188, 0.0309897)
(189, 0.0309894)
(190, 0.0309894)
(191, 0.0309892)
(192, 0.0309891)
(193, 0.030989)
(194, 0.0309889)
(195, 0.0309888)
(196, 0.0309888)
(197, 0.0309888)
(198, 0.0309888)
(199, 0.0309888)
(200, 0.0309888)
(201, 0.0309888)
(202, 0.0309888)
(203, 0.0309888)
(204, 0.0309889)
(205, 0.0309889)
(206, 0.0309889)
(207, 0.030989)
(208, 0.030989)
(209, 0.030989)
(210, 0.030989)
(211, 0.030989)
(212, 0.030989)
(213, 0.030989)
(214, 0.030989)
(215, 0.030989)
(216, 0.0309889)
(217, 0.0309889)
(218, 0.0309889)
(219, 0.0309888)
(220, 0.0309888)
(221, 0.0309887)
(222, 0.0309887)
(223, 0.0309887)
(224, 0.0309886)
(225, 0.0309886)
(226, 0.0309886)
(227, 0.0309886)
(228, 0.0309886)
(229, 0.0309886)
(230, 0.0309886)
(231, 0.0309886)
(232, 0.0309886)
(233, 0.0309887)
(234, 0.0309887)
(235, 0.0309887)
(236, 0.0309887)
(237, 0.0309888)
(238, 0.0309888)
(239, 0.0309888)
(240, 0.0309888)
(241, 0.0309888)
(242, 0.0309888)
(243, 0.0309888)
(244, 0.0309888)
(245, 0.0309888)
(246, 0.0309888)
(247, 0.0309888)
(248, 0.0309888)
(249, 0.0309888)
(250, 0.0309888)
(251, 0.0309888)
(252, 0.0309888)
(253, 0.0309888)
(254, 0.0309888)
(255, 0.0309888)
(256, 0.0309888)
(257, 0.0309888)
(258, 0.0309888)
(259, 0.0309888)
(260, 0.0309888)
(261, 0.0309888)
(262, 0.0309888)
(263, 0.0309888)
(264, 0.0309888)
(265, 0.0309888)
(266, 0.0309888)
(267, 0.0309888)
(268, 0.0309888)
(269, 0.0309888)
(270, 0.0309888)
(271, 0.0309888)
(272, 0.0309888)
(273, 0.0309888)
(274, 0.0309888)
(275, 0.0309888)
(276, 0.0309888)
(277, 0.0309888)
(278, 0.0309888)
(279, 0.0309888)
(280, 0.0309888)
(281, 0.0309888)
(282, 0.0309888)
(283, 0.0309888)
(284, 0.0309888)
(285, 0.0309888)
(286, 0.0309888)
(287, 0.0309888)
(288, 0.0309888)
(289, 0.0309888)
(290, 0.0309888)
(291, 0.0309888)
(292, 0.0309888)
(293, 0.0309888)
(294, 0.0309888)
(295, 0.0309888)
(296, 0.0309888)
(297, 0.0309888)
(298, 0.0309888)
(299, 0.0309888)
(300, 0.0309888)
(301, 0.0309888)
(302, 0.0309888)
(303, 0.0309888)
(304, 0.0309888)
(305, 0.0309888)
(306, 0.0309888)
(307, 0.0309888)
(308, 0.0309888)
(309, 0.0309888)
(310, 0.0309888)
(311, 0.0309888)
(312, 0.0309888)
(313, 0.0309888)
(314, 0.0309888)
(315, 0.0309888)
(316, 0.0309888)
(317, 0.0309888)
(318, 0.0309888)
(319, 0.0309888)
(320, 0.0309888)
(321, 0.0309888)
(322, 0.0309888)
(323, 0.0309888)
(324, 0.0309888)
(325, 0.0309888)
(326, 0.0309888)
(327, 0.0309888)
(328, 0.0309888)
(329, 0.0309888)
(330, 0.0309888)
(331, 0.0309888)
(332, 0.0309888)
(333, 0.0309888)
(334, 0.0309888)
};

%% file: 4-4.tex
\subsection{{Test problem 4: r}evisiting the problem with highly variable coefficients}
\label{4-4}
The results of \cref{4-3-2} exhibit
significant oscillations of the residual norm,
as the linear system is very ill-conditioned. 
Although preconditioned CG minimizes 
$\|\Bx_k - \Bx\|_\BA$, the ratio
$\|\Br_k\|_2 / \|\Bx_k - \Bx\|_\BA$ can range
(in principle) from $\lambda_{\min}$ to 
$\lambda_{\max}$, which allows for substantial oscillations if the condition number of $\BA$ is
large. This effect can be mitigated by better preconditioning and/or deflated versions of
CG \cite{saad2000deflated,Nicolaides:1987:DeflatedCG,Dostal1988}.
Deflation, in particular, can remove the smallest eigenvalues (and largest if desired) and drastically improve the condition number, generally leading to convergence of the residual norm that is monotonic or nearly so. 
Consequently, deflation and better preconditioning allow for a small delay parameter $d$ and more reliable behavior of error estimators, improving the efficiency of stopping criteria. 

To demonstrate the benefits of deflation in this context, we also run \cref{4-3-2},
using the recycling conjugate gradients method (recycling CG) \cite{bolten2022krylov,soodhalter2020survey,MelloStu_2010}.
Recycling CG is appropriate for 
Poisson problems as they often occur
in a sequence of linear systems
arising as the pressure Poisson 
solve in incompressible Navier-Stokes problems \cite{AmrStuSwi_2015}
In this example, the 
recycle space 
is obtained from recycling CG 
by solving the Poisson equation with the source function $f = 10 + 50 \sin x$. The recycle
space basis 
has twenty orthonormal vectors that approximate eigenvectors corresponding to the first twenty smallest eigenvalues of the linear system, and the subspace is updated every twenty CG iterations.

\begin{table}[h!]
\centering
\caption{Numbers of iterations (iter.) and quality ratios (qual. \eqref{ev}) of stopping criteria to the solution in \cref{4-3-2} solved by the preconditioned recycling CG.}
\begin{tabular}{ l|c|c|c|c|c|c} 
\hline
 \multirow{ 2}{*}{Criterion} & \multicolumn{2}{c|}{$N=4$}  & \multicolumn{2}{c|}{$N=6$} & \multicolumn{2}{c}{$N=8$} \\ \cline{2-7}
 \multirow{ 2}{*}{} & iter & qual.  & iter & qual.  & iter & qual.   \\ 
\hline
$\eta_{\text{alg}} \leq \tau \eta_{\text{R}}$ & 44 & 1.03 &67 & 1.04 &106 & 1.06  \\ 
 $\eta_{\text{alg}} \leq \tau \eta_{\text{FC}}$ &54 & 1.00 &78 & 1.00 &129 & 1.00  \\ 
 $ \|\Br_k\|_{\Bw} \leq \tau \eta_{\text{RF}}^{\Bw}$ &43 & 1.04 &66 & 1.04 &107 & 1.06  \\ 
 $ \|\Br_k^p\|_{\Bw} \leq \tau \eta_{\text{RF}}^{\Bw,p} , \, \forall p$  &95 & 1.00 &155 & 1.00 &288 & 1.00  \\ 
  $\|\Br_k\|\leq 10^{-8}\|\Br_0\|$ &131 & 1.00 &191 & 1.00 &324 & 1.00  \\ 
\hline
\end{tabular}
\label{tb:test4}
\end{table}

\begin{figure}[htbp]
\begin{tikzpicture}
\begin{groupplot}[group style={
                      group name=myplot,
                      group size= 2 by 2,
                       horizontal sep = 50pt,
                        vertical sep= 2cm
                       },height=6cm,width=7.1cm]
 \nextgroupplot[  
    xlabel={iteration},
    ymode=log,  
    ymin=1e-7,
    ymax=1e1 ]
 \input{Data/rcgomega_121}   
\nextgroupplot[ 
    xlabel={iteration},
    ymode=log,
    ymin=1e-7,
    ymax=1e2,
]
 \input{Data/rcgomega_122}   
    \end{groupplot}
\path (myplot c1r1.south west|-current bounding box.south)--
      coordinate(legendpos)
      (myplot c2r1.south east|-current bounding box.south);
\matrix[
    matrix of nodes,
    anchor=south,
    draw,
    inner sep=0.2em,
    draw
  ]at([yshift=-9ex]legendpos)
{
    \ref{plot:totalerr}& total error&[5pt]
    \ref{plot:exactalg}& $\|\Bx_k - \Bx\|_{\BA}$ &[5pt]
    \ref{plot:etarf}& $\eta_{\text{RF}}^{\Bw}$ & [5pt]
    \ref{plot:res}& $ \| \Br_k\|_{\Bw}$  & [5pt]\\
    \ref{plot:etaalg}& $\eta_{\text{alg}}$ & [5pt] 
    \ref{plot:etar}& $\eta_{\text{R}}$ & [5pt]
    \ref{plot:bdm}& \added{$\eta_{\text{FC}}$}\deleted{$\eta_{\text{BDM}}$} \\ };
\end{tikzpicture}
\caption{Convergence history of the Poisson problem with a highly variable diffusion coefficient and a good recycle space in test problem 4, with the polynomial degree $N=6$ and the delay parameter $d=10$.}
\label{iteration32_rcg}
\end{figure}

\cref{tb:test4} displays the number of iterations and quality ratios of criteria. Criteria $\eta_{\text{R}}$,  $\eta_{\text{RF}}^{\Bw}$, \added{and} \added{$\eta_{\text{FC}}$}\deleted{$\eta_{\text{BDM}}$}\deleted{, and $\eta_{\underline{\text{BDM}}}$} have desired quality ratios. 
\deleted{The indicator $\eta_{\underline{\text{BDM}}}$, as a lower bound for  $\eta_{\text{BDM}}$, requires more iterations.}
The subdomain-based criterion \deleted{also} requires a great number of extra iterations. 
If the norm of the residual exhibits a roughly monotonic decrease, employing this criterion 
becomes unnecessary.
Consistent with  previous examples, the criterion based on relative norm of the residual expends a significant number of unnecessary iterations.

\cref{iteration32_rcg} illustrates the history of errors norm and indicators for this problem solved by the preconditioned 
recycling conjugate gradient algorithm. As the well-chosen recycle subspace lessens the impact of the the ill-conditioned linear system  and consequently, the residual in the iterative process tends to decrease monotonically. 
The algebraic error estimator and a posteriori error estimators capture the behavior of exact errors well. The separation  between $\|\Br_k\|_{\Bw}$ and $\eta_{\text{RF}}^{\Bw}$ almost coincides  with the separation between the algebraic error and the total error. This is the instance where employing suitable recycle subspace can be helpful in the efficient termination of iteration process.

%% file: Data/rcgomega_121.tex
\addplot[ 
line width=2.5pt, color=gray, solid ]
coordinates { 
(1, 1.692)
(2, 1.43788)
(3, 1.22019)
(4, 1.06635)
(5, 0.944568)
(6, 0.858706)
(7, 0.794706)
(8, 0.74039)
(9, 0.684211)
(10, 0.595181)
(11, 0.475938)
(12, 0.409327)
(13, 0.355027)
(14, 0.303336)
(15, 0.261247)
(16, 0.224979)
(17, 0.197048)
(18, 0.174931)
(19, 0.156495)
(20, 0.139115)
(21, 0.124017)
(22, 0.112254)
(23, 0.101967)
(24, 0.0926251)
(25, 0.0843759)
(26, 0.0767658)
(27, 0.0700164)
(28, 0.0640419)
(29, 0.0587866)
(30, 0.054228)
(31, 0.049825)
(32, 0.0457928)
(33, 0.0419408)
(34, 0.0383392)
(35, 0.0350395)
(36, 0.0318247)
(37, 0.0287)
(38, 0.0258412)
(39, 0.0233738)
(40, 0.0211263)
(41, 0.0193154)
(42, 0.0176513)
(43, 0.016103)
(44, 0.0147155)
(45, 0.0133962)
(46, 0.0121348)
(47, 0.0109402)
(48, 0.00981905)
(49, 0.0088196)
(50, 0.00798362)
(51, 0.00730152)
(52, 0.00675126)
(53, 0.00631962)
(54, 0.00595491)
(55, 0.00563102)
(56, 0.00534136)
(57, 0.00508814)
(58, 0.0048655)
(59, 0.00467122)
(60, 0.00451157)
(61, 0.00438084)
(62, 0.00426998)
(63, 0.00417923)
(64, 0.00410754)
(65, 0.00405262)
(66, 0.00401078)
(67, 0.00397707)
(68, 0.00394999)
(69, 0.00392845)
(70, 0.00391074)
(71, 0.00389611)
(72, 0.00388407)
(73, 0.00387399)
(74, 0.00386608)
(75, 0.0038601)
(76, 0.00385599)
(77, 0.00385299)
(78, 0.00385066)
(79, 0.00384871)
(80, 0.00384684)
(81, 0.00384539)
(82, 0.00384413)
(83, 0.00384307)
(84, 0.00384222)
(85, 0.00384156)
(86, 0.00384101)
(87, 0.00384059)
(88, 0.00384026)
(89, 0.00383998)
(90, 0.00383973)
(91, 0.00383953)
(92, 0.00383936)
(93, 0.00383923)
(94, 0.00383913)
(95, 0.00383904)
(96, 0.00383898)
(97, 0.00383893)
(98, 0.00383889)
(99, 0.00383887)
(100, 0.00383884)
(101, 0.00383882)
(102, 0.0038388)
(103, 0.00383879)
(104, 0.00383877)
(105, 0.00383876)
(106, 0.00383875)
(107, 0.00383874)
(108, 0.00383874)
(109, 0.00383874)
(110, 0.00383873)
(111, 0.00383873)
(112, 0.00383873)
(113, 0.00383872)
(114, 0.00383872)
(115, 0.00383872)
(116, 0.00383872)
(117, 0.00383872)
(118, 0.00383872)
(119, 0.00383872)
(120, 0.00383872)
(121, 0.00383872)
(122, 0.00383872)
(123, 0.00383872)
(124, 0.00383872)
(125, 0.00383872)
(126, 0.00383872)
(127, 0.00383872)
(128, 0.00383872)
(129, 0.00383872)
(130, 0.00383872)
(131, 0.00383872)
(132, 0.00383872)
(133, 0.00383872)
(134, 0.00383872)
(135, 0.00383872)
(136, 0.00383872)
(137, 0.00383872)
(138, 0.00383872)
(139, 0.00383872)
(140, 0.00383872)
(141, 0.00383872)
(142, 0.00383872)
(143, 0.00383872)
(144, 0.00383872)
(145, 0.00383872)
(146, 0.00383872)
(147, 0.00383872)
(148, 0.00383872)
(149, 0.00383872)
(150, 0.00383872)
(151, 0.00383872)
(152, 0.00383872)
(153, 0.00383872)
(154, 0.00383872)
(155, 0.00383872)
(156, 0.00383872)
(157, 0.00383872)
(158, 0.00383872)
(159, 0.00383872)
(160, 0.00383872)
(161, 0.00383872)
(162, 0.00383872)
(163, 0.00383872)
(164, 0.00383872)
(165, 0.00383872)
(166, 0.00383872)
(167, 0.00383872)
(168, 0.00383872)
(169, 0.00383872)
(170, 0.00383872)
(171, 0.00383872)
(172, 0.00383872)
(173, 0.00383872)
(174, 0.00383872)
(175, 0.00383872)
(176, 0.00383872)
(177, 0.00383872)
(178, 0.00383872)
(179, 0.00383872)
(180, 0.00383872)
(181, 0.00383872)
(182, 0.00383872)
(183, 0.00383872)
(184, 0.00383872)
(185, 0.00383872)
(186, 0.00383872)
(187, 0.00383872)
(188, 0.00383872)
(189, 0.00383872)
(190, 0.00383872)
(191, 0.00383872)
(192, 0.00383872)
(193, 0.00383872)
(194, 0.00383872)
(195, 0.00383872)
(196, 0.00383872)
(197, 0.00383872)
(198, 0.00383872)
(199, 0.00383872)
(200, 0.00383872)
(201, 0.00383872)
(202, 0.00383872)
(203, 0.00383872)
(204, 0.00383872)
(205, 0.00383872)
(206, 0.00383872)
(207, 0.00383872)
};
\addplot[ 
line width=2.5pt, color=mygreen, dashed]
coordinates { 
(1, 1.692)
(2, 1.43788)
(3, 1.22019)
(4, 1.06634)
(5, 0.94456)
(6, 0.858697)
(7, 0.794696)
(8, 0.74038)
(9, 0.6842)
(10, 0.595169)
(11, 0.475922)
(12, 0.409309)
(13, 0.355006)
(14, 0.303311)
(15, 0.261219)
(16, 0.224946)
(17, 0.197011)
(18, 0.174889)
(19, 0.156448)
(20, 0.139062)
(21, 0.123958)
(22, 0.112188)
(23, 0.101895)
(24, 0.0925455)
(25, 0.0842885)
(26, 0.0766698)
(27, 0.0699111)
(28, 0.0639268)
(29, 0.0586612)
(30, 0.0540919)
(31, 0.0496769)
(32, 0.0456316)
(33, 0.0417648)
(34, 0.0381466)
(35, 0.0348286)
(36, 0.0315923)
(37, 0.0284422)
(38, 0.0255545)
(39, 0.0230564)
(40, 0.0207746)
(41, 0.0189301)
(42, 0.0172288)
(43, 0.0156387)
(44, 0.014206)
(45, 0.0128344)
(46, 0.0115116)
(47, 0.0102446)
(48, 0.00903759)
(49, 0.00794037)
(50, 0.00700018)
(51, 0.006211)
(52, 0.00555372)
(53, 0.00502014)
(54, 0.0045525)
(55, 0.00411979)
(56, 0.00371408)
(57, 0.00333968)
(58, 0.00298955)
(59, 0.00266169)
(60, 0.00237035)
(61, 0.00211094)
(62, 0.00187002)
(63, 0.00165233)
(64, 0.00146156)
(65, 0.00129921)
(66, 0.00116217)
(67, 0.00103987)
(68, 0.000930959)
(69, 0.000834861)
(70, 0.00074709)
(71, 0.000666267)
(72, 0.000591841)
(73, 0.000521582)
(74, 0.000459128)
(75, 0.000405697)
(76, 0.000364616)
(77, 0.000331387)
(78, 0.000303059)
(79, 0.000277205)
(80, 0.000249887)
(81, 0.00022645)
(82, 0.000203853)
(83, 0.000182891)
(84, 0.000164099)
(85, 0.000147812)
(86, 0.000132839)
(87, 0.00011985)
(88, 0.000108747)
(89, 9.85949e-05)
(90, 8.80264e-05)
(91, 7.89522e-05)
(92, 7.05302e-05)
(93, 6.29515e-05)
(94, 5.61051e-05)
(95, 4.98383e-05)
(96, 4.46196e-05)
(97, 4.03471e-05)
(98, 3.6844e-05)
(99, 3.38761e-05)
(100, 3.07357e-05)
(101, 2.81672e-05)
(102, 2.57078e-05)
(103, 2.32791e-05)
(104, 2.08027e-05)
(105, 1.84258e-05)
(106, 1.62568e-05)
(107, 1.44391e-05)
(108, 1.30178e-05)
(109, 1.18385e-05)
(110, 1.02098e-05)
(111, 8.96244e-06)
(112, 7.87872e-06)
(113, 6.9665e-06)
(114, 6.2046e-06)
(115, 5.55257e-06)
(116, 4.97874e-06)
(117, 4.44595e-06)
(118, 3.86553e-06)
(119, 3.32453e-06)
(120, 2.91447e-06)
(121, 2.52586e-06)
(122, 2.19576e-06)
(123, 1.90341e-06)
(124, 1.66658e-06)
(125, 1.46971e-06)
(126, 1.31515e-06)
(127, 1.18788e-06)
(128, 1.06285e-06)
(129, 9.46916e-07)
(130, 8.62963e-07)
(131, 7.87297e-07)
(132, 7.23243e-07)
(133, 6.63171e-07)
(134, 6.06578e-07)
(135, 5.52987e-07)
(136, 5.04599e-07)
(137, 4.59936e-07)
(138, 4.19237e-07)
(139, 3.81189e-07)
(140, 3.47519e-07)
(141, 3.16063e-07)
(142, 2.87433e-07)
(143, 2.62005e-07)
(144, 2.39457e-07)
(145, 2.19613e-07)
(146, 2.01707e-07)
(147, 1.86261e-07)
(148, 1.72504e-07)
(149, 1.59465e-07)
(150, 1.48768e-07)
(151, 1.39848e-07)
(152, 1.32515e-07)
(153, 1.25937e-07)
(154, 1.19843e-07)
(155, 1.14033e-07)
(156, 1.08599e-07)
(157, 1.03997e-07)
(158, 1.00394e-07)
(159, 9.76543e-08)
(160, 9.56005e-08)
(161, 9.40886e-08)
(162, 9.29386e-08)
(163, 9.20067e-08)
(164, 9.12187e-08)
(165, 9.05734e-08)
(166, 9.00757e-08)
(167, 8.97068e-08)
(168, 8.94445e-08)
(169, 8.92498e-08)
(170, 8.90807e-08)
(171, 8.89597e-08)
(172, 8.88692e-08)
(173, 8.88054e-08)
(174, 8.87672e-08)
(175, 8.87532e-08)
(176, 8.87569e-08)
(177, 8.87723e-08)
(178, 8.87929e-08)
(179, 8.88109e-08)
(180, 8.88245e-08)
(181, 8.88383e-08)
(182, 8.88533e-08)
(183, 8.88666e-08)
(184, 8.88777e-08)
(185, 8.88852e-08)
(186, 8.88894e-08)
(187, 8.88906e-08)
(188, 8.88892e-08)
(189, 8.88863e-08)
(190, 8.88821e-08)
(191, 8.88782e-08)
(192, 8.88744e-08)
(193, 8.8871e-08)
(194, 8.88683e-08)
(195, 8.88665e-08)
(196, 8.88656e-08)
(197, 8.88656e-08)
(198, 8.88663e-08)
(199, 8.88672e-08)
(200, 8.88679e-08)
(201, 8.88686e-08)
(202, 8.88695e-08)
(203, 8.88704e-08)
(204, 8.88713e-08)
(205, 8.88722e-08)
(206, 8.88729e-08)
(207, 8.88735e-08)
};
\addplot[ 
line width=1.5pt, color= mycyan, densely dotted]
coordinates { 
(1, 1.38196)
(2, 1.18333)
(3, 0.850172)
(4, 0.663407)
(5, 0.531441)
(6, 0.425305)
(7, 0.367986)
(8, 0.352936)
(9, 0.388971)
(10, 0.565017)
(11, 0.219519)
(12, 0.185748)
(13, 0.196002)
(14, 0.185553)
(15, 0.162734)
(16, 0.136539)
(17, 0.113545)
(18, 0.0936048)
(19, 0.0836247)
(20, 0.0867382)
(21, 0.0650387)
(22, 0.059562)
(23, 0.0537626)
(24, 0.0476975)
(25, 0.0437119)
(26, 0.0392647)
(27, 0.0358443)
(28, 0.0325774)
(29, 0.0274707)
(30, 0.0252699)
(31, 0.0236056)
(32, 0.0228193)
(33, 0.0212168)
(34, 0.0194208)
(35, 0.017554)
(36, 0.0174118)
(37, 0.0157942)
(38, 0.0138436)
(39, 0.0118929)
(40, 0.0103348)
(41, 0.00983906)
(42, 0.00924594)
(43, 0.00827248)
(44, 0.00731466)
(45, 0.00682945)
(46, 0.0062226)
(47, 0.00580429)
(48, 0.00531459)
(49, 0.00481995)
(50, 0.00410199)
(51, 0.00368294)
(52, 0.00308422)
(53, 0.00265663)
(54, 0.00237574)
(55, 0.00215246)
(56, 0.00197358)
(57, 0.00178817)
(58, 0.00168178)
(59, 0.00151516)
(60, 0.00129841)
(61, 0.0011952)
(62, 0.00107415)
(63, 0.000962223)
(64, 0.000854577)
(65, 0.000749791)
(66, 0.000663145)
(67, 0.000598445)
(68, 0.000522621)
(69, 0.000466394)
(70, 0.000417369)
(71, 0.000381963)
(72, 0.000351669)
(73, 0.000313726)
(74, 0.00027645)
(75, 0.000229991)
(76, 0.000192203)
(77, 0.000168679)
(78, 0.000151385)
(79, 0.000142326)
(80, 0.000138996)
(81, 0.000129044)
(82, 0.000114867)
(83, 0.000101604)
(84, 8.66906e-05)
(85, 7.82704e-05)
(86, 7.05155e-05)
(87, 6.40933e-05)
(88, 5.76729e-05)
(89, 5.49618e-05)
(90, 5.32297e-05)
(91, 4.51802e-05)
(92, 4.03859e-05)
(93, 3.52947e-05)
(94, 3.24914e-05)
(95, 2.90456e-05)
(96, 2.46619e-05)
(97, 2.13457e-05)
(98, 1.84537e-05)
(99, 1.63438e-05)
(100, 1.70758e-05)
(101, 1.44878e-05)
(102, 1.31483e-05)
(103, 1.28188e-05)
(104, 1.21626e-05)
(105, 1.10805e-05)
(106, 9.52268e-06)
(107, 8.0849e-06)
(108, 6.64296e-06)
(109, 6.11586e-06)
(110, 7.17818e-06)
(111, 5.56306e-06)
(112, 4.76741e-06)
(113, 3.93949e-06)
(114, 3.42305e-06)
(115, 3.03845e-06)
(116, 2.7326e-06)
(117, 2.74137e-06)
(118, 2.73952e-06)
(119, 2.45089e-06)
(120, 1.68734e-06)
(121, 1.36699e-06)
(122, 1.34799e-06)
(123, 1.19296e-06)
(124, 1.05213e-06)
(125, 8.76525e-07)
(126, 7.26978e-07)
(127, 6.38112e-07)
(128, 5.97592e-07)
(129, 5.38843e-07)
(130, 4.15687e-07)
(131, 3.81356e-07)
(132, 3.57511e-07)
(133, 3.26793e-07)
(134, 3.00371e-07)
(135, 2.74383e-07)
(136, 2.53953e-07)
(137, 2.34042e-07)
(138, 2.09827e-07)
(139, 1.9153e-07)
(140, 1.81905e-07)
(141, 1.79198e-07)
(142, 1.57355e-07)
(143, 1.38322e-07)
(144, 1.20335e-07)
(145, 1.08863e-07)
(146, 9.82167e-08)
(147, 8.6999e-08)
(148, 8.09649e-08)
(149, 7.11066e-08)
(150, 6.65082e-08)
(151, 5.83027e-08)
(152, 5.12262e-08)
(153, 4.84897e-08)
(154, 4.60076e-08)
(155, 4.4003e-08)
(156, 4.05162e-08)
(157, 3.57849e-08)
(158, 3.19338e-08)
(159, 2.5676e-08)
(160, 2.18998e-08)
(161, 1.91031e-08)
(162, 1.67581e-08)
(163, 1.53957e-08)
(164, 1.40537e-08)
(165, 1.17256e-08)
(166, 1.00404e-08)
(167, 8.32731e-09)
(168, 7.24381e-09)
(169, 6.24079e-09)
(170, 5.6828e-09)
(171, 5.12953e-09)
(172, 4.64119e-09)
(173, 4.19094e-09)
(174, 3.75055e-09)
(175, 3.28951e-09)
(176, 2.88322e-09)
(177, 2.39589e-09)
(178, 2.01259e-09)
(179, 1.65232e-09)
(180, 1.48386e-09)
(181, 1.41686e-09)
(182, 1.2948e-09)
(183, 1.19641e-09)
(184, 1.12747e-09)
(185, 1.07502e-09)
(186, 1.01047e-09)
(187, 9.08192e-10)
(188, 8.18458e-10)
(189, 7.12983e-10)
(190, 6.09414e-10)
(191, 5.29012e-10)
(192, 4.5596e-10)
(193, 3.93748e-10)
(194, 3.41739e-10)
(195, 3.04991e-10)
(196, 2.88605e-10)
(197, 2.67007e-10)
(198, 2.51383e-10)
(199, 2.36194e-10)
(200, 2.02731e-10)
(201, 1.7971e-10)
(202, 1.65326e-10)
(203, 1.49517e-10)
(204, 1.34245e-10)
(205, 1.18322e-10)
(206, 1.04646e-10)
(207, 9.27037e-11)
};
\addplot[ 
line width=1.5pt, color= myteal,  solid, mark=+,  mark repeat= 21, mark phase = 4]
coordinates { 
(1, 2.47507)
(2, 2.18672)
(3, 1.55729)
(4, 1.18602)
(5, 0.937311)
(6, 0.749431)
(7, 0.65592)
(8, 0.638703)
(9, 0.710119)
(10, 1.01316)
(11, 0.345588)
(12, 0.276417)
(13, 0.301376)
(14, 0.296968)
(15, 0.266264)
(16, 0.225084)
(17, 0.186171)
(18, 0.152365)
(19, 0.136601)
(20, 0.142398)
(21, 0.107327)
(22, 0.0994652)
(23, 0.0900814)
(24, 0.0797436)
(25, 0.0729656)
(26, 0.0657107)
(27, 0.0606895)
(28, 0.0560212)
(29, 0.0482307)
(30, 0.0450353)
(31, 0.0425683)
(32, 0.0413318)
(33, 0.0384964)
(34, 0.0353849)
(35, 0.0323694)
(36, 0.0321496)
(37, 0.0296683)
(38, 0.0266844)
(39, 0.0238404)
(40, 0.0218964)
(41, 0.0214534)
(42, 0.0207871)
(43, 0.0195528)
(44, 0.018323)
(45, 0.0176923)
(46, 0.0170139)
(47, 0.0165736)
(48, 0.0161252)
(49, 0.0157203)
(50, 0.015147)
(51, 0.0148472)
(52, 0.0144592)
(53, 0.0142051)
(54, 0.0140584)
(55, 0.0139547)
(56, 0.0138738)
(57, 0.0138036)
(58, 0.0137625)
(59, 0.0137047)
(60, 0.0136456)
(61, 0.0136151)
(62, 0.0135816)
(63, 0.0135584)
(64, 0.0135387)
(65, 0.013527)
(66, 0.013518)
(67, 0.0135068)
(68, 0.0134975)
(69, 0.0134936)
(70, 0.0134892)
(71, 0.0134874)
(72, 0.0134863)
(73, 0.0134845)
(74, 0.0134831)
(75, 0.0134812)
(76, 0.0134798)
(77, 0.0134792)
(78, 0.0134783)
(79, 0.0134778)
(80, 0.0134776)
(81, 0.0134771)
(82, 0.0134765)
(83, 0.0134763)
(84, 0.013476)
(85, 0.0134759)
(86, 0.0134758)
(87, 0.0134756)
(88, 0.0134755)
(89, 0.0134755)
(90, 0.0134754)
(91, 0.0134754)
(92, 0.0134755)
(93, 0.0134754)
(94, 0.0134753)
(95, 0.0134752)
(96, 0.0134751)
(97, 0.013475)
(98, 0.0134749)
(99, 0.0134749)
(100, 0.0134748)
(101, 0.0134746)
(102, 0.0134746)
(103, 0.0134745)
(104, 0.0134744)
(105, 0.0134743)
(106, 0.0134742)
(107, 0.0134741)
(108, 0.0134741)
(109, 0.0134741)
(110, 0.0134741)
(111, 0.013474)
(112, 0.013474)
(113, 0.0134739)
(114, 0.0134739)
(115, 0.0134739)
(116, 0.0134738)
(117, 0.0134738)
(118, 0.0134738)
(119, 0.0134738)
(120, 0.0134738)
(121, 0.0134738)
(122, 0.0134738)
(123, 0.0134738)
(124, 0.0134738)
(125, 0.0134738)
(126, 0.0134738)
(127, 0.0134738)
(128, 0.0134738)
(129, 0.0134738)
(130, 0.0134738)
(131, 0.0134738)
(132, 0.0134738)
(133, 0.0134738)
(134, 0.0134738)
(135, 0.0134738)
(136, 0.0134738)
(137, 0.0134738)
(138, 0.0134738)
(139, 0.0134738)
(140, 0.0134738)
(141, 0.0134738)
(142, 0.0134738)
(143, 0.0134738)
(144, 0.0134738)
(145, 0.0134738)
(146, 0.0134738)
(147, 0.0134738)
(148, 0.0134738)
(149, 0.0134738)
(150, 0.0134738)
(151, 0.0134738)
(152, 0.0134738)
(153, 0.0134738)
(154, 0.0134738)
(155, 0.0134738)
(156, 0.0134738)
(157, 0.0134738)
(158, 0.0134738)
(159, 0.0134738)
(160, 0.0134738)
(161, 0.0134738)
(162, 0.0134738)
(163, 0.0134738)
(164, 0.0134738)
(165, 0.0134738)
(166, 0.0134738)
(167, 0.0134738)
(168, 0.0134738)
(169, 0.0134738)
(170, 0.0134738)
(171, 0.0134738)
(172, 0.0134738)
(173, 0.0134738)
(174, 0.0134738)
(175, 0.0134738)
(176, 0.0134738)
(177, 0.0134738)
(178, 0.0134738)
(179, 0.0134738)
(180, 0.0134738)
(181, 0.0134738)
(182, 0.0134738)
(183, 0.0134738)
(184, 0.0134738)
(185, 0.0134738)
(186, 0.0134738)
(187, 0.0134738)
(188, 0.0134738)
(189, 0.0134738)
(190, 0.0134738)
(191, 0.0134738)
(192, 0.0134738)
(193, 0.0134738)
(194, 0.0134738)
(195, 0.0134738)
(196, 0.0134738)
(197, 0.0134738)
(198, 0.0134738)
(199, 0.0134738)
(200, 0.0134738)
(201, 0.0134738)
(202, 0.0134738)
(203, 0.0134738)
(204, 0.0134738)
(205, 0.0134738)
(206, 0.0134738)
(207, 0.0134738)
};

%% file: Data/rcgomega_122.tex
\addplot[ 
line width=2.5pt, color=gray, solid ]
coordinates { 
(1, 1.692)
(2, 1.43788)
(3, 1.22019)
(4, 1.06635)
(5, 0.944568)
(6, 0.858706)
(7, 0.794706)
(8, 0.74039)
(9, 0.684211)
(10, 0.595181)
(11, 0.475938)
(12, 0.409327)
(13, 0.355027)
(14, 0.303336)
(15, 0.261247)
(16, 0.224979)
(17, 0.197048)
(18, 0.174931)
(19, 0.156495)
(20, 0.139115)
(21, 0.124017)
(22, 0.112254)
(23, 0.101967)
(24, 0.0926251)
(25, 0.0843759)
(26, 0.0767658)
(27, 0.0700164)
(28, 0.0640419)
(29, 0.0587866)
(30, 0.054228)
(31, 0.049825)
(32, 0.0457928)
(33, 0.0419408)
(34, 0.0383392)
(35, 0.0350395)
(36, 0.0318247)
(37, 0.0287)
(38, 0.0258412)
(39, 0.0233738)
(40, 0.0211263)
(41, 0.0193154)
(42, 0.0176513)
(43, 0.016103)
(44, 0.0147155)
(45, 0.0133962)
(46, 0.0121348)
(47, 0.0109402)
(48, 0.00981905)
(49, 0.0088196)
(50, 0.00798362)
(51, 0.00730152)
(52, 0.00675126)
(53, 0.00631962)
(54, 0.00595491)
(55, 0.00563102)
(56, 0.00534136)
(57, 0.00508814)
(58, 0.0048655)
(59, 0.00467122)
(60, 0.00451157)
(61, 0.00438084)
(62, 0.00426998)
(63, 0.00417923)
(64, 0.00410754)
(65, 0.00405262)
(66, 0.00401078)
(67, 0.00397707)
(68, 0.00394999)
(69, 0.00392845)
(70, 0.00391074)
(71, 0.00389611)
(72, 0.00388407)
(73, 0.00387399)
(74, 0.00386608)
(75, 0.0038601)
(76, 0.00385599)
(77, 0.00385299)
(78, 0.00385066)
(79, 0.00384871)
(80, 0.00384684)
(81, 0.00384539)
(82, 0.00384413)
(83, 0.00384307)
(84, 0.00384222)
(85, 0.00384156)
(86, 0.00384101)
(87, 0.00384059)
(88, 0.00384026)
(89, 0.00383998)
(90, 0.00383973)
(91, 0.00383953)
(92, 0.00383936)
(93, 0.00383923)
(94, 0.00383913)
(95, 0.00383904)
(96, 0.00383898)
(97, 0.00383893)
(98, 0.00383889)
(99, 0.00383887)
(100, 0.00383884)
(101, 0.00383882)
(102, 0.0038388)
(103, 0.00383879)
(104, 0.00383877)
(105, 0.00383876)
(106, 0.00383875)
(107, 0.00383874)
(108, 0.00383874)
(109, 0.00383874)
(110, 0.00383873)
(111, 0.00383873)
(112, 0.00383873)
(113, 0.00383872)
(114, 0.00383872)
(115, 0.00383872)
(116, 0.00383872)
(117, 0.00383872)
(118, 0.00383872)
(119, 0.00383872)
(120, 0.00383872)
(121, 0.00383872)
(122, 0.00383872)
(123, 0.00383872)
(124, 0.00383872)
(125, 0.00383872)
(126, 0.00383872)
(127, 0.00383872)
(128, 0.00383872)
(129, 0.00383872)
(130, 0.00383872)
(131, 0.00383872)
(132, 0.00383872)
(133, 0.00383872)
(134, 0.00383872)
(135, 0.00383872)
(136, 0.00383872)
(137, 0.00383872)
(138, 0.00383872)
(139, 0.00383872)
(140, 0.00383872)
(141, 0.00383872)
(142, 0.00383872)
(143, 0.00383872)
(144, 0.00383872)
(145, 0.00383872)
(146, 0.00383872)
(147, 0.00383872)
(148, 0.00383872)
(149, 0.00383872)
(150, 0.00383872)
(151, 0.00383872)
(152, 0.00383872)
(153, 0.00383872)
(154, 0.00383872)
(155, 0.00383872)
(156, 0.00383872)
(157, 0.00383872)
(158, 0.00383872)
(159, 0.00383872)
(160, 0.00383872)
(161, 0.00383872)
(162, 0.00383872)
(163, 0.00383872)
(164, 0.00383872)
(165, 0.00383872)
(166, 0.00383872)
(167, 0.00383872)
(168, 0.00383872)
(169, 0.00383872)
(170, 0.00383872)
(171, 0.00383872)
(172, 0.00383872)
(173, 0.00383872)
(174, 0.00383872)
(175, 0.00383872)
(176, 0.00383872)
(177, 0.00383872)
(178, 0.00383872)
(179, 0.00383872)
(180, 0.00383872)
(181, 0.00383872)
(182, 0.00383872)
(183, 0.00383872)
(184, 0.00383872)
(185, 0.00383872)
(186, 0.00383872)
(187, 0.00383872)
(188, 0.00383872)
(189, 0.00383872)
(190, 0.00383872)
(191, 0.00383872)
(192, 0.00383872)
(193, 0.00383872)
(194, 0.00383872)
(195, 0.00383872)
(196, 0.00383872)
(197, 0.00383872)
(198, 0.00383872)
(199, 0.00383872)
(200, 0.00383872)
(201, 0.00383872)
(202, 0.00383872)
(203, 0.00383872)
(204, 0.00383872)
(205, 0.00383872)
(206, 0.00383872)
(207, 0.00383872)
};
\addplot[ 
line width=2.5pt, color=mygreen, dashed]
coordinates { 
(1, 1.692)
(2, 1.43788)
(3, 1.22019)
(4, 1.06634)
(5, 0.94456)
(6, 0.858697)
(7, 0.794696)
(8, 0.74038)
(9, 0.6842)
(10, 0.595169)
(11, 0.475922)
(12, 0.409309)
(13, 0.355006)
(14, 0.303311)
(15, 0.261219)
(16, 0.224946)
(17, 0.197011)
(18, 0.174889)
(19, 0.156448)
(20, 0.139062)
(21, 0.123958)
(22, 0.112188)
(23, 0.101895)
(24, 0.0925455)
(25, 0.0842885)
(26, 0.0766698)
(27, 0.0699111)
(28, 0.0639268)
(29, 0.0586612)
(30, 0.0540919)
(31, 0.0496769)
(32, 0.0456316)
(33, 0.0417648)
(34, 0.0381466)
(35, 0.0348286)
(36, 0.0315923)
(37, 0.0284422)
(38, 0.0255545)
(39, 0.0230564)
(40, 0.0207746)
(41, 0.0189301)
(42, 0.0172288)
(43, 0.0156387)
(44, 0.014206)
(45, 0.0128344)
(46, 0.0115116)
(47, 0.0102446)
(48, 0.00903759)
(49, 0.00794037)
(50, 0.00700018)
(51, 0.006211)
(52, 0.00555372)
(53, 0.00502014)
(54, 0.0045525)
(55, 0.00411979)
(56, 0.00371408)
(57, 0.00333968)
(58, 0.00298955)
(59, 0.00266169)
(60, 0.00237035)
(61, 0.00211094)
(62, 0.00187002)
(63, 0.00165233)
(64, 0.00146156)
(65, 0.00129921)
(66, 0.00116217)
(67, 0.00103987)
(68, 0.000930959)
(69, 0.000834861)
(70, 0.00074709)
(71, 0.000666267)
(72, 0.000591841)
(73, 0.000521582)
(74, 0.000459128)
(75, 0.000405697)
(76, 0.000364616)
(77, 0.000331387)
(78, 0.000303059)
(79, 0.000277205)
(80, 0.000249887)
(81, 0.00022645)
(82, 0.000203853)
(83, 0.000182891)
(84, 0.000164099)
(85, 0.000147812)
(86, 0.000132839)
(87, 0.00011985)
(88, 0.000108747)
(89, 9.85949e-05)
(90, 8.80264e-05)
(91, 7.89522e-05)
(92, 7.05302e-05)
(93, 6.29515e-05)
(94, 5.61051e-05)
(95, 4.98383e-05)
(96, 4.46196e-05)
(97, 4.03471e-05)
(98, 3.6844e-05)
(99, 3.38761e-05)
(100, 3.07357e-05)
(101, 2.81672e-05)
(102, 2.57078e-05)
(103, 2.32791e-05)
(104, 2.08027e-05)
(105, 1.84258e-05)
(106, 1.62568e-05)
(107, 1.44391e-05)
(108, 1.30178e-05)
(109, 1.18385e-05)
(110, 1.02098e-05)
(111, 8.96244e-06)
(112, 7.87872e-06)
(113, 6.9665e-06)
(114, 6.2046e-06)
(115, 5.55257e-06)
(116, 4.97874e-06)
(117, 4.44595e-06)
(118, 3.86553e-06)
(119, 3.32453e-06)
(120, 2.91447e-06)
(121, 2.52586e-06)
(122, 2.19576e-06)
(123, 1.90341e-06)
(124, 1.66658e-06)
(125, 1.46971e-06)
(126, 1.31515e-06)
(127, 1.18788e-06)
(128, 1.06285e-06)
(129, 9.46916e-07)
(130, 8.62963e-07)
(131, 7.87297e-07)
(132, 7.23243e-07)
(133, 6.63171e-07)
(134, 6.06578e-07)
(135, 5.52987e-07)
(136, 5.04599e-07)
(137, 4.59936e-07)
(138, 4.19237e-07)
(139, 3.81189e-07)
(140, 3.47519e-07)
(141, 3.16063e-07)
(142, 2.87433e-07)
(143, 2.62005e-07)
(144, 2.39457e-07)
(145, 2.19613e-07)
(146, 2.01707e-07)
(147, 1.86261e-07)
(148, 1.72504e-07)
(149, 1.59465e-07)
(150, 1.48768e-07)
(151, 1.39848e-07)
(152, 1.32515e-07)
(153, 1.25937e-07)
(154, 1.19843e-07)
(155, 1.14033e-07)
(156, 1.08599e-07)
(157, 1.03997e-07)
(158, 1.00394e-07)
(159, 9.76543e-08)
(160, 9.56005e-08)
(161, 9.40886e-08)
(162, 9.29386e-08)
(163, 9.20067e-08)
(164, 9.12187e-08)
(165, 9.05734e-08)
(166, 9.00757e-08)
(167, 8.97068e-08)
(168, 8.94445e-08)
(169, 8.92498e-08)
(170, 8.90807e-08)
(171, 8.89597e-08)
(172, 8.88692e-08)
(173, 8.88054e-08)
(174, 8.87672e-08)
(175, 8.87532e-08)
(176, 8.87569e-08)
(177, 8.87723e-08)
(178, 8.87929e-08)
(179, 8.88109e-08)
(180, 8.88245e-08)
(181, 8.88383e-08)
(182, 8.88533e-08)
(183, 8.88666e-08)
(184, 8.88777e-08)
(185, 8.88852e-08)
(186, 8.88894e-08)
(187, 8.88906e-08)
(188, 8.88892e-08)
(189, 8.88863e-08)
(190, 8.88821e-08)
(191, 8.88782e-08)
(192, 8.88744e-08)
(193, 8.8871e-08)
(194, 8.88683e-08)
(195, 8.88665e-08)
(196, 8.88656e-08)
(197, 8.88656e-08)
(198, 8.88663e-08)
(199, 8.88672e-08)
(200, 8.88679e-08)
(201, 8.88686e-08)
(202, 8.88695e-08)
(203, 8.88704e-08)
(204, 8.88713e-08)
(205, 8.88722e-08)
(206, 8.88729e-08)
(207, 8.88735e-08)
};
\addplot[ 
line width=1.5pt, color= myolive,  solid, mark=square, mark repeat= 21, mark phase = 7]
coordinates { 
(1, 1.62368)
(2, 1.37839)
(3, 1.1674)
(4, 1.0223)
(5, 0.907721)
(6, 0.82871)
(7, 0.769889)
(8, 0.719428)
(9, 0.666074)
(10, 0.578695)
(11, 0.459496)
(12, 0.393634)
(13, 0.340069)
(14, 0.288848)
(15, 0.247246)
(16, 0.211477)
(17, 0.18419)
(18, 0.162787)
(19, 0.145034)
(20, 0.128111)
(21, 0.113568)
(22, 0.102488)
(23, 0.0929426)
(24, 0.0843179)
(25, 0.0767563)
(26, 0.0698583)
(27, 0.063864)
(28, 0.0585969)
(29, 0.0539401)
(30, 0.0499435)
(31, 0.0459287)
(32, 0.0422541)
(33, 0.0387263)
(34, 0.0354027)
(35, 0.0323776)
(36, 0.0294204)
(37, 0.0265331)
(38, 0.023903)
(39, 0.021646)
(40, 0.0195597)
(41, 0.0178822)
(42, 0.0163091)
(43, 0.0148111)
(44, 0.0134568)
(45, 0.0121553)
(46, 0.010896)
(47, 0.00968495)
(48, 0.00852881)
(49, 0.00748097)
(50, 0.00658665)
(51, 0.00584127)
(52, 0.00522942)
(53, 0.00474043)
(54, 0.00431151)
(55, 0.00390956)
(56, 0.00352757)
(57, 0.00317366)
(58, 0.0028409)
(59, 0.00252737)
(60, 0.00224953)
(61, 0.00200303)
(62, 0.00177389)
(63, 0.00156784)
(64, 0.00138757)
(65, 0.00123425)
(66, 0.00110349)
(67, 0.000985656)
(68, 0.000880249)
(69, 0.000787496)
(70, 0.000704059)
(71, 0.000626604)
(72, 0.000555626)
(73, 0.000488466)
(74, 0.000428802)
(75, 0.000377812)
(76, 0.000339557)
(77, 0.000308955)
(78, 0.000282876)
(79, 0.000259078)
(80, 0.000233867)
(81, 0.000212239)
(82, 0.00019126)
(83, 0.000171713)
(84, 0.000154206)
(85, 0.000139153)
(86, 0.000125118)
(87, 0.000112851)
(88, 0.000102312)
(89, 9.25892e-05)
(90, 8.24833e-05)
(91, 7.37542e-05)
(92, 6.5676e-05)
(93, 5.84876e-05)
(94, 5.21049e-05)
(95, 4.63065e-05)
(96, 4.15526e-05)
(97, 3.76752e-05)
(98, 3.44682e-05)
(99, 3.1741e-05)
(100, 2.89915e-05)
(101, 2.67046e-05)
(102, 2.44722e-05)
(103, 2.22137e-05)
(104, 1.98573e-05)
(105, 1.75706e-05)
(106, 1.54768e-05)
(107, 1.37386e-05)
(108, 1.24314e-05)
(109, 1.13627e-05)
(110, 9.78538e-06)
(111, 8.5993e-06)
(112, 7.56653e-06)
(113, 6.70124e-06)
(114, 5.97627e-06)
(115, 5.35414e-06)
(116, 4.80151e-06)
(117, 4.284e-06)
(118, 3.71635e-06)
(119, 3.18679e-06)
(120, 2.78396e-06)
(121, 2.40048e-06)
(122, 2.07397e-06)
(123, 1.78519e-06)
(124, 1.55363e-06)
(125, 1.36338e-06)
(126, 1.21648e-06)
(127, 1.09752e-06)
(128, 9.79283e-07)
(129, 8.69719e-07)
(130, 7.92999e-07)
(131, 7.24293e-07)
(132, 6.66939e-07)
(133, 6.1243e-07)
(134, 5.60328e-07)
(135, 5.10196e-07)
(136, 4.64789e-07)
(137, 4.22282e-07)
(138, 3.83284e-07)
(139, 3.46866e-07)
(140, 3.14332e-07)
(141, 2.8338e-07)
(142, 2.54717e-07)
(143, 2.29181e-07)
(144, 2.06647e-07)
(145, 1.87108e-07)
(146, 1.69662e-07)
(147, 1.54604e-07)
(148, 1.40853e-07)
(149, 1.27192e-07)
(150, 1.15595e-07)
(151, 1.05501e-07)
(152, 9.68489e-08)
(153, 8.86075e-08)
(154, 8.04064e-08)
(155, 7.1828e-08)
(156, 6.2887e-08)
(157, 5.43369e-08)
(158, 4.6708e-08)
(159, 4.01774e-08)
(160, 3.4886e-08)
(161, 3.05648e-08)
(162, 2.68814e-08)
(163, 2.3584e-08)
(164, 2.05018e-08)
(165, 1.77339e-08)
(166, 1.53718e-08)
(167, 1.3418e-08)
(168, 1.18698e-08)
(169, 1.05789e-08)
(170, 9.31191e-09)
(171, 8.25989e-09)
(172, 7.31856e-09)
(173, 6.42594e-09)
(174, 5.59705e-09)
(175, 4.85871e-09)
(176, 4.2306e-09)
(177, 3.72482e-09)
(178, 3.33411e-09)
(179, 3.05495e-09)
(180, 2.79423e-09)
(181, 2.56747e-09)
(182, 2.35878e-09)
(183, 2.16285e-09)
(184, 1.97005e-09)
(185, 1.77731e-09)
(186, 1.58115e-09)
(187, 1.39791e-09)
(188, 1.22066e-09)
(189, 1.06668e-09)
(190, 9.25939e-10)
(191, 8.24027e-10)
(192, 7.41923e-10)
(193, 6.75199e-10)
(194, 6.17973e-10)
(195, 5.68939e-10)
(196, 5.23331e-10)
(197, 4.79147e-10)
(198, 0)
(199, 0)
(200, 0)
(201, 0)
(202, 0)
(203, 0)
(204, 0)
(205, 0)
(206, 0)
(207, 0)
};
\addplot[ 
line width=1.5pt, color= mypurple,  solid, mark=o,  mark repeat= 21, mark phase = 2]
coordinates { 
(1, 2.36155)
(2, 2.07009)
(3, 1.44004)
(4, 1.09926)
(5, 0.870976)
(6, 0.697549)
(7, 0.609562)
(8, 0.593634)
(9, 0.671447)
(10, 0.971934)
(11, 0.337619)
(12, 0.278683)
(13, 0.303327)
(14, 0.295715)
(15, 0.261438)
(16, 0.221503)
(17, 0.184023)
(18, 0.149829)
(19, 0.134563)
(20, 0.142361)
(21, 0.106226)
(22, 0.0981226)
(23, 0.0896184)
(24, 0.0793777)
(25, 0.0729657)
(26, 0.0666874)
(27, 0.0621369)
(28, 0.0573742)
(29, 0.049553)
(30, 0.04629)
(31, 0.0438652)
(32, 0.0427316)
(33, 0.0403614)
(34, 0.0379097)
(35, 0.0354082)
(36, 0.0352586)
(37, 0.0329858)
(38, 0.0303628)
(39, 0.0280369)
(40, 0.0265485)
(41, 0.0261624)
(42, 0.0255343)
(43, 0.0245368)
(44, 0.0235817)
(45, 0.0231127)
(46, 0.0226335)
(47, 0.0222851)
(48, 0.0219349)
(49, 0.0216533)
(50, 0.0212384)
(51, 0.0210271)
(52, 0.020772)
(53, 0.0206122)
(54, 0.0205239)
(55, 0.0204666)
(56, 0.0204166)
(57, 0.0203745)
(58, 0.0203454)
(59, 0.0203063)
(60, 0.0202721)
(61, 0.0202487)
(62, 0.0202214)
(63, 0.0202059)
(64, 0.0201936)
(65, 0.0201913)
(66, 0.0201898)
(67, 0.0201795)
(68, 0.0201735)
(69, 0.020174)
(70, 0.0201716)
(71, 0.020171)
(72, 0.0201717)
(73, 0.0201713)
(74, 0.0201714)
(75, 0.0201707)
(76, 0.0201704)
(77, 0.0201704)
(78, 0.0201698)
(79, 0.0201692)
(80, 0.0201693)
(81, 0.0201689)
(82, 0.0201684)
(83, 0.0201687)
(84, 0.0201688)
(85, 0.0201689)
(86, 0.020169)
(87, 0.0201687)
(88, 0.0201688)
(89, 0.0201688)
(90, 0.0201688)
(91, 0.0201689)
(92, 0.0201691)
(93, 0.0201691)
(94, 0.020169)
(95, 0.0201688)
(96, 0.0201687)
(97, 0.0201687)
(98, 0.0201686)
(99, 0.0201686)
(100, 0.0201685)
(101, 0.0201682)
(102, 0.0201682)
(103, 0.0201681)
(104, 0.0201679)
(105, 0.0201678)
(106, 0.0201677)
(107, 0.0201676)
(108, 0.0201676)
(109, 0.0201676)
(110, 0.0201676)
(111, 0.0201675)
(112, 0.0201675)
(113, 0.0201674)
(114, 0.0201674)
(115, 0.0201674)
(116, 0.0201673)
(117, 0.0201673)
(118, 0.0201673)
(119, 0.0201673)
(120, 0.0201673)
(121, 0.0201673)
(122, 0.0201673)
(123, 0.0201673)
(124, 0.0201673)
(125, 0.0201673)
(126, 0.0201673)
(127, 0.0201673)
(128, 0.0201673)
(129, 0.0201673)
(130, 0.0201673)
(131, 0.0201673)
(132, 0.0201673)
(133, 0.0201673)
(134, 0.0201673)
(135, 0.0201673)
(136, 0.0201673)
(137, 0.0201673)
(138, 0.0201673)
(139, 0.0201673)
(140, 0.0201673)
(141, 0.0201673)
(142, 0.0201673)
(143, 0.0201673)
(144, 0.0201673)
(145, 0.0201673)
(146, 0.0201673)
(147, 0.0201673)
(148, 0.0201673)
(149, 0.0201673)
(150, 0.0201673)
(151, 0.0201673)
(152, 0.0201673)
(153, 0.0201673)
(154, 0.0201673)
(155, 0.0201673)
(156, 0.0201673)
(157, 0.0201673)
(158, 0.0201673)
(159, 0.0201673)
(160, 0.0201673)
(161, 0.0201673)
(162, 0.0201673)
(163, 0.0201673)
(164, 0.0201673)
(165, 0.0201673)
(166, 0.0201673)
(167, 0.0201673)
(168, 0.0201673)
(169, 0.0201673)
(170, 0.0201673)
(171, 0.0201673)
(172, 0.0201673)
(173, 0.0201673)
(174, 0.0201673)
(175, 0.0201673)
(176, 0.0201673)
(177, 0.0201673)
(178, 0.0201673)
(179, 0.0201673)
(180, 0.0201673)
(181, 0.0201673)
(182, 0.0201673)
(183, 0.0201673)
(184, 0.0201673)
(185, 0.0201673)
(186, 0.0201673)
(187, 0.0201673)
(188, 0.0201673)
(189, 0.0201673)
(190, 0.0201673)
(191, 0.0201673)
(192, 0.0201673)
(193, 0.0201673)
(194, 0.0201673)
(195, 0.0201673)
(196, 0.0201673)
(197, 0.0201673)
(198, 0.0201673)
(199, 0.0201673)
(200, 0.0201673)
(201, 0.0201673)
(202, 0.0201673)
(203, 0.0201673)
(204, 0.0201673)
(205, 0.0201673)
(206, 0.0201673)
(207, 0.0201673)
};
\addplot[ 
line width=1.5pt, color= myrose ,   solid, mark=x, mark repeat= 21, mark phase = 16]
coordinates { 
(1, 0.371914)
(2, 0.311356)
(3, 0.205601)
(4, 0.148841)
(5, 0.116161)
(6, 0.0952938)
(7, 0.0835557)
(8, 0.086533)
(9, 0.102955)
(10, 0.155187)
(11, 0.0524461)
(12, 0.042462)
(13, 0.0479079)
(14, 0.0481426)
(15, 0.0433385)
(16, 0.0354836)
(17, 0.0286522)
(18, 0.0239469)
(19, 0.0216416)
(20, 0.0211919)
(21, 0.0168957)
(22, 0.0158789)
(23, 0.0144463)
(24, 0.0129546)
(25, 0.0119392)
(26, 0.0112169)
(27, 0.0107594)
(28, 0.0103687)
(29, 0.00953055)
(30, 0.00925311)
(31, 0.00892346)
(32, 0.00879465)
(33, 0.00844732)
(34, 0.0080048)
(35, 0.00769974)
(36, 0.00752781)
(37, 0.00728185)
(38, 0.00693937)
(39, 0.0066765)
(40, 0.00651563)
(41, 0.00644751)
(42, 0.00637829)
(43, 0.00626865)
(44, 0.00618778)
(45, 0.00613282)
(46, 0.00609038)
(47, 0.00604315)
(48, 0.00599402)
(49, 0.00596167)
(50, 0.0059276)
(51, 0.00591096)
(52, 0.00590116)
(53, 0.00588915)
(54, 0.00588514)
(55, 0.00588351)
(56, 0.00587869)
(57, 0.00587684)
(58, 0.00587322)
(59, 0.00586973)
(60, 0.00586873)
(61, 0.00586731)
(62, 0.00586265)
(63, 0.00586068)
(64, 0.00585991)
(65, 0.00586132)
(66, 0.00586241)
(67, 0.00586149)
(68, 0.00586048)
(69, 0.00586136)
(70, 0.00586134)
(71, 0.00586182)
(72, 0.00586186)
(73, 0.00586213)
(74, 0.00586223)
(75, 0.00586246)
(76, 0.00586262)
(77, 0.00586269)
(78, 0.00586261)
(79, 0.00586259)
(80, 0.00586254)
(81, 0.00586265)
(82, 0.00586266)
(83, 0.00586278)
(84, 0.00586293)
(85, 0.00586308)
(86, 0.00586313)
(87, 0.00586316)
(88, 0.00586318)
(89, 0.00586323)
(90, 0.00586332)
(91, 0.00586337)
(92, 0.00586343)
(93, 0.00586349)
(94, 0.00586343)
(95, 0.00586341)
(96, 0.00586342)
(97, 0.00586339)
(98, 0.00586338)
(99, 0.00586338)
(100, 0.00586332)
(101, 0.00586329)
(102, 0.00586327)
(103, 0.00586323)
(104, 0.0058632)
(105, 0.00586316)
(106, 0.00586313)
(107, 0.00586311)
(108, 0.0058631)
(109, 0.0058631)
(110, 0.00586309)
(111, 0.00586308)
(112, 0.00586307)
(113, 0.00586305)
(114, 0.00586304)
(115, 0.00586303)
(116, 0.00586303)
(117, 0.00586303)
(118, 0.00586303)
(119, 0.00586303)
(120, 0.00586302)
(121, 0.00586302)
(122, 0.00586303)
(123, 0.00586302)
(124, 0.00586302)
(125, 0.00586302)
(126, 0.00586302)
(127, 0.00586303)
(128, 0.00586303)
(129, 0.00586303)
(130, 0.00586303)
(131, 0.00586303)
(132, 0.00586303)
(133, 0.00586303)
(134, 0.00586303)
(135, 0.00586303)
(136, 0.00586303)
(137, 0.00586303)
(138, 0.00586303)
(139, 0.00586303)
(140, 0.00586303)
(141, 0.00586303)
(142, 0.00586303)
(143, 0.00586303)
(144, 0.00586303)
(145, 0.00586303)
(146, 0.00586303)
(147, 0.00586303)
(148, 0.00586303)
(149, 0.00586303)
(150, 0.00586303)
(151, 0.00586303)
(152, 0.00586303)
(153, 0.00586303)
(154, 0.00586303)
(155, 0.00586303)
(156, 0.00586303)
(157, 0.00586303)
(158, 0.00586303)
(159, 0.00586303)
(160, 0.00586303)
(161, 0.00586303)
(162, 0.00586303)
(163, 0.00586303)
(164, 0.00586303)
(165, 0.00586303)
(166, 0.00586303)
(167, 0.00586303)
(168, 0.00586303)
(169, 0.00586303)
(170, 0.00586303)
(171, 0.00586303)
(172, 0.00586303)
(173, 0.00586303)
(174, 0.00586303)
(175, 0.00586303)
(176, 0.00586303)
(177, 0.00586303)
(178, 0.00586303)
(179, 0.00586303)
(180, 0.00586303)
(181, 0.00586303)
(182, 0.00586303)
(183, 0.00586303)
(184, 0.00586303)
(185, 0.00586303)
(186, 0.00586303)
(187, 0.00586303)
(188, 0.00586303)
(189, 0.00586303)
(190, 0.00586303)
(191, 0.00586303)
(192, 0.00586303)
(193, 0.00586303)
(194, 0.00586303)
(195, 0.00586303)
(196, 0.00586303)
(197, 0.00586303)
(198, 0.00586303)
(199, 0.00586303)
(200, 0.00586303)
(201, 0.00586303)
(202, 0.00586303)
(203, 0.00586303)
(204, 0.00586303)
(205, 0.00586303)
(206, 0.00586303)
(207, 0.00586303)
};
\addplot[ 
line width=1.5pt, color= mysand,  solid, mark=diamond, mark repeat= 21, mark phase = 16]
coordinates { 
(1, 0.236913)
(2, 0.203147)
(3, 0.136481)
(4, 0.0985551)
(5, 0.0760475)
(6, 0.061863)
(7, 0.0545029)
(8, 0.0563331)
(9, 0.0662766)
(10, 0.0984884)
(11, 0.0324684)
(12, 0.0256182)
(13, 0.0286222)
(14, 0.0288516)
(15, 0.0260752)
(16, 0.0215035)
(17, 0.0173588)
(18, 0.0144752)
(19, 0.0130279)
(20, 0.0127973)
(21, 0.010053)
(22, 0.00938988)
(23, 0.00849485)
(24, 0.00755646)
(25, 0.00697087)
(26, 0.00648833)
(27, 0.00615741)
(28, 0.00586725)
(29, 0.00528045)
(30, 0.00506924)
(31, 0.00485386)
(32, 0.00480659)
(33, 0.00458209)
(34, 0.00429159)
(35, 0.00406492)
(36, 0.00395115)
(37, 0.0037719)
(38, 0.00352888)
(39, 0.00333712)
(40, 0.00321135)
(41, 0.0031623)
(42, 0.00311491)
(43, 0.00303044)
(44, 0.00297132)
(45, 0.00293044)
(46, 0.00290086)
(47, 0.0028677)
(48, 0.00283636)
(49, 0.00281101)
(50, 0.00278374)
(51, 0.00276902)
(52, 0.00275754)
(53, 0.00274649)
(54, 0.00274206)
(55, 0.00273928)
(56, 0.00273513)
(57, 0.00273312)
(58, 0.00273026)
(59, 0.00272767)
(60, 0.00272563)
(61, 0.0027248)
(62, 0.00272159)
(63, 0.00272052)
(64, 0.00271954)
(65, 0.00272001)
(66, 0.00272006)
(67, 0.00271968)
(68, 0.00271875)
(69, 0.00271915)
(70, 0.00271886)
(71, 0.00271908)
(72, 0.00271892)
(73, 0.00271904)
(74, 0.0027189)
(75, 0.002719)
(76, 0.00271895)
(77, 0.00271898)
(78, 0.00271888)
(79, 0.00271887)
(80, 0.0027188)
(81, 0.00271886)
(82, 0.00271882)
(83, 0.00271887)
(84, 0.00271891)
(85, 0.00271898)
(86, 0.00271899)
(87, 0.002719)
(88, 0.002719)
(89, 0.00271902)
(90, 0.00271906)
(91, 0.00271907)
(92, 0.0027191)
(93, 0.00271912)
(94, 0.00271909)
(95, 0.00271908)
(96, 0.00271908)
(97, 0.00271906)
(98, 0.00271906)
(99, 0.00271906)
(100, 0.00271903)
(101, 0.00271902)
(102, 0.00271901)
(103, 0.00271899)
(104, 0.00271897)
(105, 0.00271895)
(106, 0.00271894)
(107, 0.00271893)
(108, 0.00271892)
(109, 0.00271892)
(110, 0.00271892)
(111, 0.00271891)
(112, 0.00271891)
(113, 0.0027189)
(114, 0.00271889)
(115, 0.00271889)
(116, 0.00271889)
(117, 0.00271889)
(118, 0.00271889)
(119, 0.00271889)
(120, 0.00271889)
(121, 0.00271889)
(122, 0.00271889)
(123, 0.00271889)
(124, 0.00271889)
(125, 0.00271889)
(126, 0.00271889)
(127, 0.00271889)
(128, 0.00271889)
(129, 0.00271889)
(130, 0.00271889)
(131, 0.00271889)
(132, 0.00271889)
(133, 0.00271889)
(134, 0.00271889)
(135, 0.00271889)
(136, 0.00271889)
(137, 0.00271889)
(138, 0.00271889)
(139, 0.00271889)
(140, 0.00271889)
(141, 0.00271889)
(142, 0.00271889)
(143, 0.00271889)
(144, 0.00271889)
(145, 0.00271889)
(146, 0.00271889)
(147, 0.00271889)
(148, 0.00271889)
(149, 0.00271889)
(150, 0.00271889)
(151, 0.00271889)
(152, 0.00271889)
(153, 0.00271889)
(154, 0.00271889)
(155, 0.00271889)
(156, 0.00271889)
(157, 0.00271889)
(158, 0.00271889)
(159, 0.00271889)
(160, 0.00271889)
(161, 0.00271889)
(162, 0.00271889)
(163, 0.00271889)
(164, 0.00271889)
(165, 0.00271889)
(166, 0.00271889)
(167, 0.00271889)
(168, 0.00271889)
(169, 0.00271889)
(170, 0.00271889)
(171, 0.00271889)
(172, 0.00271889)
(173, 0.00271889)
(174, 0.00271889)
(175, 0.00271889)
(176, 0.00271889)
(177, 0.00271889)
(178, 0.00271889)
(179, 0.00271889)
(180, 0.00271889)
(181, 0.00271889)
(182, 0.00271889)
(183, 0.00271889)
(184, 0.00271889)
(185, 0.00271889)
(186, 0.00271889)
(187, 0.00271889)
(188, 0.00271889)
(189, 0.00271889)
(190, 0.00271889)
(191, 0.00271889)
(192, 0.00271889)
(193, 0.00271889)
(194, 0.00271889)
(195, 0.00271889)
(196, 0.00271889)
(197, 0.00271889)
(198, 0.00271889)
(199, 0.00271889)
(200, 0.00271889)
(201, 0.00271889)
(202, 0.00271889)
(203, 0.00271889)
(204, 0.00271889)
(205, 0.00271889)
(206, 0.00271889)
(207, 0.00271889)
};
\addplot[ 
line width=2.5pt, color=gray, solid ]
coordinates { 
(1, 1.692)
(2, 1.43788)
(3, 1.22019)
(4, 1.06635)
(5, 0.944568)
(6, 0.858706)
(7, 0.794706)
(8, 0.74039)
(9, 0.684211)
(10, 0.595181)
(11, 0.475938)
(12, 0.409327)
(13, 0.355027)
(14, 0.303336)
(15, 0.261247)
(16, 0.224979)
(17, 0.197048)
(18, 0.174931)
(19, 0.156495)
(20, 0.139115)
(21, 0.124017)
(22, 0.112254)
(23, 0.101967)
(24, 0.0926251)
(25, 0.0843759)
(26, 0.0767658)
(27, 0.0700164)
(28, 0.0640419)
(29, 0.0587866)
(30, 0.054228)
(31, 0.049825)
(32, 0.0457928)
(33, 0.0419408)
(34, 0.0383392)
(35, 0.0350395)
(36, 0.0318247)
(37, 0.0287)
(38, 0.0258412)
(39, 0.0233738)
(40, 0.0211263)
(41, 0.0193154)
(42, 0.0176513)
(43, 0.016103)
(44, 0.0147155)
(45, 0.0133962)
(46, 0.0121348)
(47, 0.0109402)
(48, 0.00981905)
(49, 0.0088196)
(50, 0.00798362)
(51, 0.00730152)
(52, 0.00675126)
(53, 0.00631962)
(54, 0.00595491)
(55, 0.00563102)
(56, 0.00534136)
(57, 0.00508814)
(58, 0.0048655)
(59, 0.00467122)
(60, 0.00451157)
(61, 0.00438084)
(62, 0.00426998)
(63, 0.00417923)
(64, 0.00410754)
(65, 0.00405262)
(66, 0.00401078)
(67, 0.00397707)
(68, 0.00394999)
(69, 0.00392845)
(70, 0.00391074)
(71, 0.00389611)
(72, 0.00388407)
(73, 0.00387399)
(74, 0.00386608)
(75, 0.0038601)
(76, 0.00385599)
(77, 0.00385299)
(78, 0.00385066)
(79, 0.00384871)
(80, 0.00384684)
(81, 0.00384539)
(82, 0.00384413)
(83, 0.00384307)
(84, 0.00384222)
(85, 0.00384156)
(86, 0.00384101)
(87, 0.00384059)
(88, 0.00384026)
(89, 0.00383998)
(90, 0.00383973)
(91, 0.00383953)
(92, 0.00383936)
(93, 0.00383923)
(94, 0.00383913)
(95, 0.00383904)
(96, 0.00383898)
(97, 0.00383893)
(98, 0.00383889)
(99, 0.00383887)
(100, 0.00383884)
(101, 0.00383882)
(102, 0.0038388)
(103, 0.00383879)
(104, 0.00383877)
(105, 0.00383876)
(106, 0.00383875)
(107, 0.00383874)
(108, 0.00383874)
(109, 0.00383874)
(110, 0.00383873)
(111, 0.00383873)
(112, 0.00383873)
(113, 0.00383872)
(114, 0.00383872)
(115, 0.00383872)
(116, 0.00383872)
(117, 0.00383872)
(118, 0.00383872)
(119, 0.00383872)
(120, 0.00383872)
(121, 0.00383872)
(122, 0.00383872)
(123, 0.00383872)
(124, 0.00383872)
(125, 0.00383872)
(126, 0.00383872)
(127, 0.00383872)
(128, 0.00383872)
(129, 0.00383872)
(130, 0.00383872)
(131, 0.00383872)
(132, 0.00383872)
(133, 0.00383872)
(134, 0.00383872)
(135, 0.00383872)
(136, 0.00383872)
(137, 0.00383872)
(138, 0.00383872)
(139, 0.00383872)
(140, 0.00383872)
(141, 0.00383872)
(142, 0.00383872)
(143, 0.00383872)
(144, 0.00383872)
(145, 0.00383872)
(146, 0.00383872)
(147, 0.00383872)
(148, 0.00383872)
(149, 0.00383872)
(150, 0.00383872)
(151, 0.00383872)
(152, 0.00383872)
(153, 0.00383872)
(154, 0.00383872)
(155, 0.00383872)
(156, 0.00383872)
(157, 0.00383872)
(158, 0.00383872)
(159, 0.00383872)
(160, 0.00383872)
(161, 0.00383872)
(162, 0.00383872)
(163, 0.00383872)
(164, 0.00383872)
(165, 0.00383872)
(166, 0.00383872)
(167, 0.00383872)
(168, 0.00383872)
(169, 0.00383872)
(170, 0.00383872)
(171, 0.00383872)
(172, 0.00383872)
(173, 0.00383872)
(174, 0.00383872)
(175, 0.00383872)
(176, 0.00383872)
(177, 0.00383872)
(178, 0.00383872)
(179, 0.00383872)
(180, 0.00383872)
(181, 0.00383872)
(182, 0.00383872)
(183, 0.00383872)
(184, 0.00383872)
(185, 0.00383872)
(186, 0.00383872)
(187, 0.00383872)
(188, 0.00383872)
(189, 0.00383872)
(190, 0.00383872)
(191, 0.00383872)
(192, 0.00383872)
(193, 0.00383872)
(194, 0.00383872)
(195, 0.00383872)
(196, 0.00383872)
(197, 0.00383872)
(198, 0.00383872)
(199, 0.00383872)
(200, 0.00383872)
(201, 0.00383872)
(202, 0.00383872)
(203, 0.00383872)
(204, 0.00383872)
(205, 0.00383872)
(206, 0.00383872)
(207, 0.00383872)
};
\addplot[ 
line width=2.5pt, color=mygreen, dashed]
coordinates { 
(1, 1.692)
(2, 1.43788)
(3, 1.22019)
(4, 1.06634)
(5, 0.94456)
(6, 0.858697)
(7, 0.794696)
(8, 0.74038)
(9, 0.6842)
(10, 0.595169)
(11, 0.475922)
(12, 0.409309)
(13, 0.355006)
(14, 0.303311)
(15, 0.261219)
(16, 0.224946)
(17, 0.197011)
(18, 0.174889)
(19, 0.156448)
(20, 0.139062)
(21, 0.123958)
(22, 0.112188)
(23, 0.101895)
(24, 0.0925455)
(25, 0.0842885)
(26, 0.0766698)
(27, 0.0699111)
(28, 0.0639268)
(29, 0.0586612)
(30, 0.0540919)
(31, 0.0496769)
(32, 0.0456316)
(33, 0.0417648)
(34, 0.0381466)
(35, 0.0348286)
(36, 0.0315923)
(37, 0.0284422)
(38, 0.0255545)
(39, 0.0230564)
(40, 0.0207746)
(41, 0.0189301)
(42, 0.0172288)
(43, 0.0156387)
(44, 0.014206)
(45, 0.0128344)
(46, 0.0115116)
(47, 0.0102446)
(48, 0.00903759)
(49, 0.00794037)
(50, 0.00700018)
(51, 0.006211)
(52, 0.00555372)
(53, 0.00502014)
(54, 0.0045525)
(55, 0.00411979)
(56, 0.00371408)
(57, 0.00333968)
(58, 0.00298955)
(59, 0.00266169)
(60, 0.00237035)
(61, 0.00211094)
(62, 0.00187002)
(63, 0.00165233)
(64, 0.00146156)
(65, 0.00129921)
(66, 0.00116217)
(67, 0.00103987)
(68, 0.000930959)
(69, 0.000834861)
(70, 0.00074709)
(71, 0.000666267)
(72, 0.000591841)
(73, 0.000521582)
(74, 0.000459128)
(75, 0.000405697)
(76, 0.000364616)
(77, 0.000331387)
(78, 0.000303059)
(79, 0.000277205)
(80, 0.000249887)
(81, 0.00022645)
(82, 0.000203853)
(83, 0.000182891)
(84, 0.000164099)
(85, 0.000147812)
(86, 0.000132839)
(87, 0.00011985)
(88, 0.000108747)
(89, 9.85949e-05)
(90, 8.80264e-05)
(91, 7.89522e-05)
(92, 7.05302e-05)
(93, 6.29515e-05)
(94, 5.61051e-05)
(95, 4.98383e-05)
(96, 4.46196e-05)
(97, 4.03471e-05)
(98, 3.6844e-05)
(99, 3.38761e-05)
(100, 3.07357e-05)
(101, 2.81672e-05)
(102, 2.57078e-05)
(103, 2.32791e-05)
(104, 2.08027e-05)
(105, 1.84258e-05)
(106, 1.62568e-05)
(107, 1.44391e-05)
(108, 1.30178e-05)
(109, 1.18385e-05)
(110, 1.02098e-05)
(111, 8.96244e-06)
(112, 7.87872e-06)
(113, 6.9665e-06)
(114, 6.2046e-06)
(115, 5.55257e-06)
(116, 4.97874e-06)
(117, 4.44595e-06)
(118, 3.86553e-06)
(119, 3.32453e-06)
(120, 2.91447e-06)
(121, 2.52586e-06)
(122, 2.19576e-06)
(123, 1.90341e-06)
(124, 1.66658e-06)
(125, 1.46971e-06)
(126, 1.31515e-06)
(127, 1.18788e-06)
(128, 1.06285e-06)
(129, 9.46916e-07)
(130, 8.62963e-07)
(131, 7.87297e-07)
(132, 7.23243e-07)
(133, 6.63171e-07)
(134, 6.06578e-07)
(135, 5.52987e-07)
(136, 5.04599e-07)
(137, 4.59936e-07)
(138, 4.19237e-07)
(139, 3.81189e-07)
(140, 3.47519e-07)
(141, 3.16063e-07)
(142, 2.87433e-07)
(143, 2.62005e-07)
(144, 2.39457e-07)
(145, 2.19613e-07)
(146, 2.01707e-07)
(147, 1.86261e-07)
(148, 1.72504e-07)
(149, 1.59465e-07)
(150, 1.48768e-07)
(151, 1.39848e-07)
(152, 1.32515e-07)
(153, 1.25937e-07)
(154, 1.19843e-07)
(155, 1.14033e-07)
(156, 1.08599e-07)
(157, 1.03997e-07)
(158, 1.00394e-07)
(159, 9.76543e-08)
(160, 9.56005e-08)
(161, 9.40886e-08)
(162, 9.29386e-08)
(163, 9.20067e-08)
(164, 9.12187e-08)
(165, 9.05734e-08)
(166, 9.00757e-08)
(167, 8.97068e-08)
(168, 8.94445e-08)
(169, 8.92498e-08)
(170, 8.90807e-08)
(171, 8.89597e-08)
(172, 8.88692e-08)
(173, 8.88054e-08)
(174, 8.87672e-08)
(175, 8.87532e-08)
(176, 8.87569e-08)
(177, 8.87723e-08)
(178, 8.87929e-08)
(179, 8.88109e-08)
(180, 8.88245e-08)
(181, 8.88383e-08)
(182, 8.88533e-08)
(183, 8.88666e-08)
(184, 8.88777e-08)
(185, 8.88852e-08)
(186, 8.88894e-08)
(187, 8.88906e-08)
(188, 8.88892e-08)
(189, 8.88863e-08)
(190, 8.88821e-08)
(191, 8.88782e-08)
(192, 8.88744e-08)
(193, 8.8871e-08)
(194, 8.88683e-08)
(195, 8.88665e-08)
(196, 8.88656e-08)
(197, 8.88656e-08)
(198, 8.88663e-08)
(199, 8.88672e-08)
(200, 8.88679e-08)
(201, 8.88686e-08)
(202, 8.88695e-08)
(203, 8.88704e-08)
(204, 8.88713e-08)
(205, 8.88722e-08)
(206, 8.88729e-08)
(207, 8.88735e-08)
};
\addplot[ 
line width=1.5pt, color= myolive,  solid, mark=square, mark repeat= 21, mark phase = 7]
coordinates { 
(1, 1.62368)
(2, 1.37839)
(3, 1.1674)
(4, 1.0223)
(5, 0.907721)
(6, 0.82871)
(7, 0.769889)
(8, 0.719428)
(9, 0.666074)
(10, 0.578695)
(11, 0.459496)
(12, 0.393634)
(13, 0.340069)
(14, 0.288848)
(15, 0.247246)
(16, 0.211477)
(17, 0.18419)
(18, 0.162787)
(19, 0.145034)
(20, 0.128111)
(21, 0.113568)
(22, 0.102488)
(23, 0.0929426)
(24, 0.0843179)
(25, 0.0767563)
(26, 0.0698583)
(27, 0.063864)
(28, 0.0585969)
(29, 0.0539401)
(30, 0.0499435)
(31, 0.0459287)
(32, 0.0422541)
(33, 0.0387263)
(34, 0.0354027)
(35, 0.0323776)
(36, 0.0294204)
(37, 0.0265331)
(38, 0.023903)
(39, 0.021646)
(40, 0.0195597)
(41, 0.0178822)
(42, 0.0163091)
(43, 0.0148111)
(44, 0.0134568)
(45, 0.0121553)
(46, 0.010896)
(47, 0.00968495)
(48, 0.00852881)
(49, 0.00748097)
(50, 0.00658665)
(51, 0.00584127)
(52, 0.00522942)
(53, 0.00474043)
(54, 0.00431151)
(55, 0.00390956)
(56, 0.00352757)
(57, 0.00317366)
(58, 0.0028409)
(59, 0.00252737)
(60, 0.00224953)
(61, 0.00200303)
(62, 0.00177389)
(63, 0.00156784)
(64, 0.00138757)
(65, 0.00123425)
(66, 0.00110349)
(67, 0.000985656)
(68, 0.000880249)
(69, 0.000787496)
(70, 0.000704059)
(71, 0.000626604)
(72, 0.000555626)
(73, 0.000488466)
(74, 0.000428802)
(75, 0.000377812)
(76, 0.000339557)
(77, 0.000308955)
(78, 0.000282876)
(79, 0.000259078)
(80, 0.000233867)
(81, 0.000212239)
(82, 0.00019126)
(83, 0.000171713)
(84, 0.000154206)
(85, 0.000139153)
(86, 0.000125118)
(87, 0.000112851)
(88, 0.000102312)
(89, 9.25892e-05)
(90, 8.24833e-05)
(91, 7.37542e-05)
(92, 6.5676e-05)
(93, 5.84876e-05)
(94, 5.21049e-05)
(95, 4.63065e-05)
(96, 4.15526e-05)
(97, 3.76752e-05)
(98, 3.44682e-05)
(99, 3.1741e-05)
(100, 2.89915e-05)
(101, 2.67046e-05)
(102, 2.44722e-05)
(103, 2.22137e-05)
(104, 1.98573e-05)
(105, 1.75706e-05)
(106, 1.54768e-05)
(107, 1.37386e-05)
(108, 1.24314e-05)
(109, 1.13627e-05)
(110, 9.78538e-06)
(111, 8.5993e-06)
(112, 7.56653e-06)
(113, 6.70124e-06)
(114, 5.97627e-06)
(115, 5.35414e-06)
(116, 4.80151e-06)
(117, 4.284e-06)
(118, 3.71635e-06)
(119, 3.18679e-06)
(120, 2.78396e-06)
(121, 2.40048e-06)
(122, 2.07397e-06)
(123, 1.78519e-06)
(124, 1.55363e-06)
(125, 1.36338e-06)
(126, 1.21648e-06)
(127, 1.09752e-06)
(128, 9.79283e-07)
(129, 8.69719e-07)
(130, 7.92999e-07)
(131, 7.24293e-07)
(132, 6.66939e-07)
(133, 6.1243e-07)
(134, 5.60328e-07)
(135, 5.10196e-07)
(136, 4.64789e-07)
(137, 4.22282e-07)
(138, 3.83284e-07)
(139, 3.46866e-07)
(140, 3.14332e-07)
(141, 2.8338e-07)
(142, 2.54717e-07)
(143, 2.29181e-07)
(144, 2.06647e-07)
(145, 1.87108e-07)
(146, 1.69662e-07)
(147, 1.54604e-07)
(148, 1.40853e-07)
(149, 1.27192e-07)
(150, 1.15595e-07)
(151, 1.05501e-07)
(152, 9.68489e-08)
(153, 8.86075e-08)
(154, 8.04064e-08)
(155, 7.1828e-08)
(156, 6.2887e-08)
(157, 5.43369e-08)
(158, 4.6708e-08)
(159, 4.01774e-08)
(160, 3.4886e-08)
(161, 3.05648e-08)
(162, 2.68814e-08)
(163, 2.3584e-08)
(164, 2.05018e-08)
(165, 1.77339e-08)
(166, 1.53718e-08)
(167, 1.3418e-08)
(168, 1.18698e-08)
(169, 1.05789e-08)
(170, 9.31191e-09)
(171, 8.25989e-09)
(172, 7.31856e-09)
(173, 6.42594e-09)
(174, 5.59705e-09)
(175, 4.85871e-09)
(176, 4.2306e-09)
(177, 3.72482e-09)
(178, 3.33411e-09)
(179, 3.05495e-09)
(180, 2.79423e-09)
(181, 2.56747e-09)
(182, 2.35878e-09)
(183, 2.16285e-09)
(184, 1.97005e-09)
(185, 1.77731e-09)
(186, 1.58115e-09)
(187, 1.39791e-09)
(188, 1.22066e-09)
(189, 1.06668e-09)
(190, 9.25939e-10)
(191, 8.24027e-10)
(192, 7.41923e-10)
(193, 6.75199e-10)
(194, 6.17973e-10)
(195, 5.68939e-10)
(196, 5.23331e-10)
(197, 4.79147e-10)
(198, 0)
(199, 0)
(200, 0)
(201, 0)
(202, 0)
(203, 0)
(204, 0)
(205, 0)
(206, 0)
(207, 0)
};
\addplot[ 
line width=1.5pt, color= mypurple,  solid, mark=o,  mark repeat= 21, mark phase = 2]
coordinates { 
(1, 2.36155)
(2, 2.07009)
(3, 1.44004)
(4, 1.09926)
(5, 0.870976)
(6, 0.697549)
(7, 0.609562)
(8, 0.593634)
(9, 0.671447)
(10, 0.971934)
(11, 0.337619)
(12, 0.278683)
(13, 0.303327)
(14, 0.295715)
(15, 0.261438)
(16, 0.221503)
(17, 0.184023)
(18, 0.149829)
(19, 0.134563)
(20, 0.142361)
(21, 0.106226)
(22, 0.0981226)
(23, 0.0896184)
(24, 0.0793777)
(25, 0.0729657)
(26, 0.0666874)
(27, 0.0621369)
(28, 0.0573742)
(29, 0.049553)
(30, 0.04629)
(31, 0.0438652)
(32, 0.0427316)
(33, 0.0403614)
(34, 0.0379097)
(35, 0.0354082)
(36, 0.0352586)
(37, 0.0329858)
(38, 0.0303628)
(39, 0.0280369)
(40, 0.0265485)
(41, 0.0261624)
(42, 0.0255343)
(43, 0.0245368)
(44, 0.0235817)
(45, 0.0231127)
(46, 0.0226335)
(47, 0.0222851)
(48, 0.0219349)
(49, 0.0216533)
(50, 0.0212384)
(51, 0.0210271)
(52, 0.020772)
(53, 0.0206122)
(54, 0.0205239)
(55, 0.0204666)
(56, 0.0204166)
(57, 0.0203745)
(58, 0.0203454)
(59, 0.0203063)
(60, 0.0202721)
(61, 0.0202487)
(62, 0.0202214)
(63, 0.0202059)
(64, 0.0201936)
(65, 0.0201913)
(66, 0.0201898)
(67, 0.0201795)
(68, 0.0201735)
(69, 0.020174)
(70, 0.0201716)
(71, 0.020171)
(72, 0.0201717)
(73, 0.0201713)
(74, 0.0201714)
(75, 0.0201707)
(76, 0.0201704)
(77, 0.0201704)
(78, 0.0201698)
(79, 0.0201692)
(80, 0.0201693)
(81, 0.0201689)
(82, 0.0201684)
(83, 0.0201687)
(84, 0.0201688)
(85, 0.0201689)
(86, 0.020169)
(87, 0.0201687)
(88, 0.0201688)
(89, 0.0201688)
(90, 0.0201688)
(91, 0.0201689)
(92, 0.0201691)
(93, 0.0201691)
(94, 0.020169)
(95, 0.0201688)
(96, 0.0201687)
(97, 0.0201687)
(98, 0.0201686)
(99, 0.0201686)
(100, 0.0201685)
(101, 0.0201682)
(102, 0.0201682)
(103, 0.0201681)
(104, 0.0201679)
(105, 0.0201678)
(106, 0.0201677)
(107, 0.0201676)
(108, 0.0201676)
(109, 0.0201676)
(110, 0.0201676)
(111, 0.0201675)
(112, 0.0201675)
(113, 0.0201674)
(114, 0.0201674)
(115, 0.0201674)
(116, 0.0201673)
(117, 0.0201673)
(118, 0.0201673)
(119, 0.0201673)
(120, 0.0201673)
(121, 0.0201673)
(122, 0.0201673)
(123, 0.0201673)
(124, 0.0201673)
(125, 0.0201673)
(126, 0.0201673)
(127, 0.0201673)
(128, 0.0201673)
(129, 0.0201673)
(130, 0.0201673)
(131, 0.0201673)
(132, 0.0201673)
(133, 0.0201673)
(134, 0.0201673)
(135, 0.0201673)
(136, 0.0201673)
(137, 0.0201673)
(138, 0.0201673)
(139, 0.0201673)
(140, 0.0201673)
(141, 0.0201673)
(142, 0.0201673)
(143, 0.0201673)
(144, 0.0201673)
(145, 0.0201673)
(146, 0.0201673)
(147, 0.0201673)
(148, 0.0201673)
(149, 0.0201673)
(150, 0.0201673)
(151, 0.0201673)
(152, 0.0201673)
(153, 0.0201673)
(154, 0.0201673)
(155, 0.0201673)
(156, 0.0201673)
(157, 0.0201673)
(158, 0.0201673)
(159, 0.0201673)
(160, 0.0201673)
(161, 0.0201673)
(162, 0.0201673)
(163, 0.0201673)
(164, 0.0201673)
(165, 0.0201673)
(166, 0.0201673)
(167, 0.0201673)
(168, 0.0201673)
(169, 0.0201673)
(170, 0.0201673)
(171, 0.0201673)
(172, 0.0201673)
(173, 0.0201673)
(174, 0.0201673)
(175, 0.0201673)
(176, 0.0201673)
(177, 0.0201673)
(178, 0.0201673)
(179, 0.0201673)
(180, 0.0201673)
(181, 0.0201673)
(182, 0.0201673)
(183, 0.0201673)
(184, 0.0201673)
(185, 0.0201673)
(186, 0.0201673)
(187, 0.0201673)
(188, 0.0201673)
(189, 0.0201673)
(190, 0.0201673)
(191, 0.0201673)
(192, 0.0201673)
(193, 0.0201673)
(194, 0.0201673)
(195, 0.0201673)
(196, 0.0201673)
(197, 0.0201673)
(198, 0.0201673)
(199, 0.0201673)
(200, 0.0201673)
(201, 0.0201673)
(202, 0.0201673)
(203, 0.0201673)
(204, 0.0201673)
(205, 0.0201673)
(206, 0.0201673)
(207, 0.0201673)
};
\addplot[ 
line width=1.5pt, color= myrose ,   solid, mark=x, mark repeat= 21, mark phase = 16]
coordinates { 
(1, 0.371914)
(2, 0.311356)
(3, 0.205601)
(4, 0.148841)
(5, 0.116161)
(6, 0.0952938)
(7, 0.0835557)
(8, 0.086533)
(9, 0.102955)
(10, 0.155187)
(11, 0.0524461)
(12, 0.042462)
(13, 0.0479079)
(14, 0.0481426)
(15, 0.0433385)
(16, 0.0354836)
(17, 0.0286522)
(18, 0.0239469)
(19, 0.0216416)
(20, 0.0211919)
(21, 0.0168957)
(22, 0.0158789)
(23, 0.0144463)
(24, 0.0129546)
(25, 0.0119392)
(26, 0.0112169)
(27, 0.0107594)
(28, 0.0103687)
(29, 0.00953055)
(30, 0.00925311)
(31, 0.00892346)
(32, 0.00879465)
(33, 0.00844732)
(34, 0.0080048)
(35, 0.00769974)
(36, 0.00752781)
(37, 0.00728185)
(38, 0.00693937)
(39, 0.0066765)
(40, 0.00651563)
(41, 0.00644751)
(42, 0.00637829)
(43, 0.00626865)
(44, 0.00618778)
(45, 0.00613282)
(46, 0.00609038)
(47, 0.00604315)
(48, 0.00599402)
(49, 0.00596167)
(50, 0.0059276)
(51, 0.00591096)
(52, 0.00590116)
(53, 0.00588915)
(54, 0.00588514)
(55, 0.00588351)
(56, 0.00587869)
(57, 0.00587684)
(58, 0.00587322)
(59, 0.00586973)
(60, 0.00586873)
(61, 0.00586731)
(62, 0.00586265)
(63, 0.00586068)
(64, 0.00585991)
(65, 0.00586132)
(66, 0.00586241)
(67, 0.00586149)
(68, 0.00586048)
(69, 0.00586136)
(70, 0.00586134)
(71, 0.00586182)
(72, 0.00586186)
(73, 0.00586213)
(74, 0.00586223)
(75, 0.00586246)
(76, 0.00586262)
(77, 0.00586269)
(78, 0.00586261)
(79, 0.00586259)
(80, 0.00586254)
(81, 0.00586265)
(82, 0.00586266)
(83, 0.00586278)
(84, 0.00586293)
(85, 0.00586308)
(86, 0.00586313)
(87, 0.00586316)
(88, 0.00586318)
(89, 0.00586323)
(90, 0.00586332)
(91, 0.00586337)
(92, 0.00586343)
(93, 0.00586349)
(94, 0.00586343)
(95, 0.00586341)
(96, 0.00586342)
(97, 0.00586339)
(98, 0.00586338)
(99, 0.00586338)
(100, 0.00586332)
(101, 0.00586329)
(102, 0.00586327)
(103, 0.00586323)
(104, 0.0058632)
(105, 0.00586316)
(106, 0.00586313)
(107, 0.00586311)
(108, 0.0058631)
(109, 0.0058631)
(110, 0.00586309)
(111, 0.00586308)
(112, 0.00586307)
(113, 0.00586305)
(114, 0.00586304)
(115, 0.00586303)
(116, 0.00586303)
(117, 0.00586303)
(118, 0.00586303)
(119, 0.00586303)
(120, 0.00586302)
(121, 0.00586302)
(122, 0.00586303)
(123, 0.00586302)
(124, 0.00586302)
(125, 0.00586302)
(126, 0.00586302)
(127, 0.00586303)
(128, 0.00586303)
(129, 0.00586303)
(130, 0.00586303)
(131, 0.00586303)
(132, 0.00586303)
(133, 0.00586303)
(134, 0.00586303)
(135, 0.00586303)
(136, 0.00586303)
(137, 0.00586303)
(138, 0.00586303)
(139, 0.00586303)
(140, 0.00586303)
(141, 0.00586303)
(142, 0.00586303)
(143, 0.00586303)
(144, 0.00586303)
(145, 0.00586303)
(146, 0.00586303)
(147, 0.00586303)
(148, 0.00586303)
(149, 0.00586303)
(150, 0.00586303)
(151, 0.00586303)
(152, 0.00586303)
(153, 0.00586303)
(154, 0.00586303)
(155, 0.00586303)
(156, 0.00586303)
(157, 0.00586303)
(158, 0.00586303)
(159, 0.00586303)
(160, 0.00586303)
(161, 0.00586303)
(162, 0.00586303)
(163, 0.00586303)
(164, 0.00586303)
(165, 0.00586303)
(166, 0.00586303)
(167, 0.00586303)
(168, 0.00586303)
(169, 0.00586303)
(170, 0.00586303)
(171, 0.00586303)
(172, 0.00586303)
(173, 0.00586303)
(174, 0.00586303)
(175, 0.00586303)
(176, 0.00586303)
(177, 0.00586303)
(178, 0.00586303)
(179, 0.00586303)
(180, 0.00586303)
(181, 0.00586303)
(182, 0.00586303)
(183, 0.00586303)
(184, 0.00586303)
(185, 0.00586303)
(186, 0.00586303)
(187, 0.00586303)
(188, 0.00586303)
(189, 0.00586303)
(190, 0.00586303)
(191, 0.00586303)
(192, 0.00586303)
(193, 0.00586303)
(194, 0.00586303)
(195, 0.00586303)
(196, 0.00586303)
(197, 0.00586303)
(198, 0.00586303)
(199, 0.00586303)
(200, 0.00586303)
(201, 0.00586303)
(202, 0.00586303)
(203, 0.00586303)
(204, 0.00586303)
(205, 0.00586303)
(206, 0.00586303)
(207, 0.00586303)
};
\addplot[ 
line width=1.5pt, color= mysand,  solid, mark=diamond, mark repeat= 21, mark phase = 16]
coordinates { 
(1, 0.236913)
(2, 0.203147)
(3, 0.136481)
(4, 0.0985551)
(5, 0.0760475)
(6, 0.061863)
(7, 0.0545029)
(8, 0.0563331)
(9, 0.0662766)
(10, 0.0984884)
(11, 0.0324684)
(12, 0.0256182)
(13, 0.0286222)
(14, 0.0288516)
(15, 0.0260752)
(16, 0.0215035)
(17, 0.0173588)
(18, 0.0144752)
(19, 0.0130279)
(20, 0.0127973)
(21, 0.010053)
(22, 0.00938988)
(23, 0.00849485)
(24, 0.00755646)
(25, 0.00697087)
(26, 0.00648833)
(27, 0.00615741)
(28, 0.00586725)
(29, 0.00528045)
(30, 0.00506924)
(31, 0.00485386)
(32, 0.00480659)
(33, 0.00458209)
(34, 0.00429159)
(35, 0.00406492)
(36, 0.00395115)
(37, 0.0037719)
(38, 0.00352888)
(39, 0.00333712)
(40, 0.00321135)
(41, 0.0031623)
(42, 0.00311491)
(43, 0.00303044)
(44, 0.00297132)
(45, 0.00293044)
(46, 0.00290086)
(47, 0.0028677)
(48, 0.00283636)
(49, 0.00281101)
(50, 0.00278374)
(51, 0.00276902)
(52, 0.00275754)
(53, 0.00274649)
(54, 0.00274206)
(55, 0.00273928)
(56, 0.00273513)
(57, 0.00273312)
(58, 0.00273026)
(59, 0.00272767)
(60, 0.00272563)
(61, 0.0027248)
(62, 0.00272159)
(63, 0.00272052)
(64, 0.00271954)
(65, 0.00272001)
(66, 0.00272006)
(67, 0.00271968)
(68, 0.00271875)
(69, 0.00271915)
(70, 0.00271886)
(71, 0.00271908)
(72, 0.00271892)
(73, 0.00271904)
(74, 0.0027189)
(75, 0.002719)
(76, 0.00271895)
(77, 0.00271898)
(78, 0.00271888)
(79, 0.00271887)
(80, 0.0027188)
(81, 0.00271886)
(82, 0.00271882)
(83, 0.00271887)
(84, 0.00271891)
(85, 0.00271898)
(86, 0.00271899)
(87, 0.002719)
(88, 0.002719)
(89, 0.00271902)
(90, 0.00271906)
(91, 0.00271907)
(92, 0.0027191)
(93, 0.00271912)
(94, 0.00271909)
(95, 0.00271908)
(96, 0.00271908)
(97, 0.00271906)
(98, 0.00271906)
(99, 0.00271906)
(100, 0.00271903)
(101, 0.00271902)
(102, 0.00271901)
(103, 0.00271899)
(104, 0.00271897)
(105, 0.00271895)
(106, 0.00271894)
(107, 0.00271893)
(108, 0.00271892)
(109, 0.00271892)
(110, 0.00271892)
(111, 0.00271891)
(112, 0.00271891)
(113, 0.0027189)
(114, 0.00271889)
(115, 0.00271889)
(116, 0.00271889)
(117, 0.00271889)
(118, 0.00271889)
(119, 0.00271889)
(120, 0.00271889)
(121, 0.00271889)
(122, 0.00271889)
(123, 0.00271889)
(124, 0.00271889)
(125, 0.00271889)
(126, 0.00271889)
(127, 0.00271889)
(128, 0.00271889)
(129, 0.00271889)
(130, 0.00271889)
(131, 0.00271889)
(132, 0.00271889)
(133, 0.00271889)
(134, 0.00271889)
(135, 0.00271889)
(136, 0.00271889)
(137, 0.00271889)
(138, 0.00271889)
(139, 0.00271889)
(140, 0.00271889)
(141, 0.00271889)
(142, 0.00271889)
(143, 0.00271889)
(144, 0.00271889)
(145, 0.00271889)
(146, 0.00271889)
(147, 0.00271889)
(148, 0.00271889)
(149, 0.00271889)
(150, 0.00271889)
(151, 0.00271889)
(152, 0.00271889)
(153, 0.00271889)
(154, 0.00271889)
(155, 0.00271889)
(156, 0.00271889)
(157, 0.00271889)
(158, 0.00271889)
(159, 0.00271889)
(160, 0.00271889)
(161, 0.00271889)
(162, 0.00271889)
(163, 0.00271889)
(164, 0.00271889)
(165, 0.00271889)
(166, 0.00271889)
(167, 0.00271889)
(168, 0.00271889)
(169, 0.00271889)
(170, 0.00271889)
(171, 0.00271889)
(172, 0.00271889)
(173, 0.00271889)
(174, 0.00271889)
(175, 0.00271889)
(176, 0.00271889)
(177, 0.00271889)
(178, 0.00271889)
(179, 0.00271889)
(180, 0.00271889)
(181, 0.00271889)
(182, 0.00271889)
(183, 0.00271889)
(184, 0.00271889)
(185, 0.00271889)
(186, 0.00271889)
(187, 0.00271889)
(188, 0.00271889)
(189, 0.00271889)
(190, 0.00271889)
(191, 0.00271889)
(192, 0.00271889)
(193, 0.00271889)
(194, 0.00271889)
(195, 0.00271889)
(196, 0.00271889)
(197, 0.00271889)
(198, 0.00271889)
(199, 0.00271889)
(200, 0.00271889)
(201, 0.00271889)
(202, 0.00271889)
(203, 0.00271889)
(204, 0.00271889)
(205, 0.00271889)
(206, 0.00271889)
(207, 0.00271889)
};

%% file: 4-5.tex
\subsection{\added{Runtime Comparison}}
\label{4-5}
\added{
To quantitatively demonstrate the effectiveness of stopping criterion  \eqref{eq:criterion}, we consider the Poisson equation on $[-0.5,0.5]^3$ with Dirichlet boundary conditions and forcing function 
$f=\sin(2\pi{}x)\sin(2\pi{}y)\sin(2\pi{}z)$, using hexahedral spectral elements on Kershaw mesh ($\varepsilon=0.3$) \cite{ kershaw1981differencing}, which is used as the basis of benchmark problems by the Center for Efficient Exascale Discretization \cite{kolev2021high}.  The number of elements is chosen to be multiple of 6 in each axis to align with the Kershaw regions and such that there are a total of approximately 4M degrees of freedom(DoFs).
Solving the linear system by the conjugate gradient algorithm and setting $\sigma=0.1$, we run the calculations on an NVIDIA H100 SXM GPU using polynomial degree $N=3, \dots, 8$.
The experiments were performed using libParanumal \cite{ChalmersKarakusAustinSwirydowiczWarburton2020} and employed a highly optimized matrix-free  \cite{swirydowicz2019acceleration} preconditioned conjugate gradient algorithm with FP64 outer iteration precision and FP32 preconditioner precision. The matrix-free operations are implemented using the OCCA API and OKL kernel language \cite{occa} and executed using the CUDA backend.
Following the approach of \cite{karakus2019gpu}, the preconditioner applies a hybrid multigrid preconditioner with pMG reducing polynomial degrees to approximately halve the number of DoFs per level and a second order Chebyshev smoother, paired with an AMG hybrid multigrid hierarchy for the coarse grid solve.
}

\added{
In Table \ref{tb:4}, the results of applying the relative residual norm criterion and the $\eta_{\text{RF}}^{\Bw}$ stopping criterion every second iteration are presented. 
As evidenced by the iteration counts and total errors, a fixed residual relative convergence tolerance can result in significant over-iteration and it can be particularly acute for low order discretizations when comparing in time-to-solution for calculations with the same number of DoFs.
The last two columns display the iteration count ratio and runtime ratio for applying the $\eta_{\text{RF}}^{\Bw}$  stopping criterion compared to the relative residual criterion.
The results indicate that the additional overhead of evaluating  the indicator $\eta_{\text{RF}}^{\Bw}$ every second iteration is minimal, with the difference between the iteration count ratio and time ratio being less than 10\%.
}
\added{
\setlength{\tabcolsep}{4pt}
\renewcommand{\arraystretch}{1.1} 
\begin{table*}[htbp]
\centering
\caption{\added{Relative performance of the error estimate based stopping criterion \eqref{eq:criterion} applied every second iteration versus a more standard relative residual based stopping criterion. 
}}
\begin{tabular}{ c|c|c c c c |c c c c | c | c} 
\hline
\multirow{2}{*}{N} & \multirow{2}{*}{DoFs} &  \multicolumn{4}{c|}{ $\|\Br_k\|_{\Bw} \leq 10^{-1} \eta_{\text{RF}}^{\Bw}$ } & \multicolumn{4}{c|}{$\|\Br_k\|_{l^2}\leq 10^{-10} \|\Br_0\|_{l^2}$}  & Iter & Time \\ \cline{3-10}
\multirow{2}{*}{}  &  \multirow{2}{*}{}  & Iter & $\|u-u^k_h\|_{E}$ & $\|\Br_k\|_{l^2}$ & Time & Iter & $\|u-u^k_h\|_{E}$ & $\|\Br_k\|_{l^2}$  & Time & ratio & ratio  \\ 
\hline
    3 & 4.2M &    24 & 7.6e-03 & 1.6e-04 & 1.3e-01 &    78 & 7.3e-03 & 8.4e-11 & 3.9e-01 & 3.2 &    3.0 \\
    4 & 4.7M &    38 & 1.1e-03 & 1.6e-05 & 1.6e-01 &    96 & 1.0e-03 & 8.4e-11 & 3.8e-01 & 2.5 &    2.3 \\
    5 & 5.7M &    46 & 1.4e-04 & 2.4e-06 & 2.6e-01 &    93 & 1.4e-04 & 8.1e-11 & 5.0e-01 & 2.0 &    1.9 \\
    6 & 5.7M &    60 & 2.2e-05 & 3.7e-07 & 3.2e-01 &   104 & 2.2e-05 & 9.0e-11 & 5.2e-01 & 1.7 &    1.6 \\
    7 & 4.7M &    64 & 6.8e-06 & 9.8e-08 & 3.2e-01 &    99 & 6.7e-06 & 9.2e-11 & 4.6e-01 & 1.5 &    1.5 \\
    8 & 7.0M &    84 & 4.2e-07 & 7.8e-09 & 5.8e-01 &   106 & 4.1e-07 & 9.8e-11 & 6.9e-01 & 1.3 &    1.2 \\
\hline
\end{tabular}
\label{tb:4}
\end{table*}
}

%% file: 4-6ResultsSummary.tex
\subsection{Results summary}

\deleted{
We present a summary of the numerical experiments in this section as a score sheet in \cref{comparison}.
The methods are scored by whether they lead to premature termination with a quality ratio greater than 2 (0 points), whether they suggest stopping slightly early or slightly late (1 point), or if the quality ratio is less than 1.5 and the iteration stops at roughly the optimal point, ensuring adequate accuracy and fewer iterations (2 points). 
}
\deleted{
We summarize the performance of tested stopping criteria  as follows:}
\added{The estimator $\eta_{\text{FC}}$ provides the most accurate estimate for the total error; however, it is computationally expensive. The criterion $ \| \Br_k\|_{\Bw} \leq \tau \eta_{\text{RF}}^{\Bw}$ offers a competitive option as long as $\eta_{\text{RF}}^{\Bw}$ closely tracks $ \|\Br_k\|_{\Bw}$, which is usually the case except for \cref{4-3-2}. The subdomain-based  criterion is the only one that provides reliable termination for \cref{4-3-2}. Moreover, as presented in \cref{4-5}, $\eta_{\text{RF}}^{\Bw}$ is inexpensive to compute.
}
All criteria depending on the algebraic error estimator $\eta_{\text{alg}}$ fail  when the algebraic error remains almost constant for a relatively large number of iterations. In such cases, with a small delay parameter $d$, the estimator $\eta_{\text{alg}}$ is not accurate, and selecting an appropriate $d$ can be challenging. Additionally, in practice, the criteria relying on $\eta_{\text{alg}}$ 
include an additional $d$ iterations (with $d=10$ in all experiments) required to compute $\eta_{\text{alg}}$. 

%% file: 5-conclusion.tex
\section{Conclusion}
\label{sec:conclusions}

In this study, we have presented \deleted{three}\added{two} new stopping criteria and compared the proposed criteria with several existing stopping criteria for the conjugate gradient algorithm within the context of high-order finite elements for solving the Poisson equation.\added{
We have established reliability and efficiency theorems to ensure that the criterion prevents both over-solving and under-solving.}

\deleted{Criterion \ref{c4}  compares the error indicator $\eta_{\text{\underline{BDM}}}$ to the estimate of the algebraic error. 
The indicator $\eta_{\text{\underline{BDM}}}$ is computationally less expensive than $\eta_\text{{BDM}}$, and the associated stopping criterion is more robust than $\eta_\text{{BDM}}$ criterion.}

Criterion \ref{c5} compares error indicator $\eta_{\text{RF}}^{\Bw}$  to the weighted norm of the residual. The indicator is a natural upper bound for the weighted norm of the residual without involving any unknown constants. This criterion, which closely relies on the residual, offers advantages over criteria based on algebraic error estimation and a posteriori error estimation. It eliminates the difficulty of selecting an appropriate delay parameter in algebraic error estimation and has a more favorable computational cost. Furthermore, it is robust with respect to the mesh size, the polynomial degree, and the shape regularity of the mesh.

\added{Moreover}\deleted{Thirdly}, we proposed a subdomain-based criterion \ref{c6} for solving the Poisson problem with highly variable piecewise constant coefficient.
This stopping criterion terminates when the criterion is individually satisfied for each subdomain. It is the only tested criterion that ensures reliable termination for \cref{4-3-2} with highly variable coefficients in the absence of a good preconditioner
or deflation. 

 For problems with highly variable piecewise constant coefficients, criteria \ref{c1}, \ref{c3}, 
 and \ref{c5} recommend  termination at a reasonable iteration for \cref{4-3-1}, but they are not reliable for \cref{4-3-2}. In such cases, criterion \ref{c6} is used instead. However, it remains unclear when it is necessary to switch to the subdomain-based criterion. Further investigation is planned for future work.
Additionally, we plan to extend these criteria to  more general problems, such as problems with continuous variable coefficients and mixed problems. Furthermore, it is natural to consider applying these criteria to nonconforming finite element methods.

%% file: 6-appendix.tex
\appendix

\input{6-1appendixProof}

\section{Error indicators in subdomains}
\label{appen}

An element is defined as an overlap element if at least one of its edges lies on the interface of $\Omega_1$, $\Omega_2$, or $\Omega_3$. An element is an interior element if the element and all its edges are located in the interior of  $\Omega_1$, $\Omega_2$, or $\Omega_3$. An element is an exterior element if it is neither an overlap nor an interior element.
The nodes that are present in the overlap elements are referred to as overlap nodes, and the set of all overlap nodes  is represented by $\mathcal{S}_{\text{o}}$. Conversely, $\mathcal{S}_{\text{i}}$ represents the set of nodes that belong to the interior elements but not the overlap elements. Similarly, $\mathcal{S}_{\text{e}}$ denotes the set of nodes that belong to the exterior elements but not the overlap elements.

Since $\|u - u^k_h\|_E$ is the sum of  errors on all elements, 
we can define 
\added{$\|u - u^k_h\|_{E,\text{i}} $, $\|u - u^k_h\|_{E,\text{o}}$, and $\|u - u^k_h\|_{E,\text{e}}$}
\deleted{$\|u - u^k_h\|_{a,\text{i}} $, $\|u - u^k_h\|_{a,\text{o}}$, and $\|u - u^k_h\|_{a,\text{e}}$}
as the sum of errors on all elements in the interior subdomain, overlap subdomain, and exterior subdomain, respectively. We define the subdomain algebraic errors in a similar manner.

On the other hand, $\eta_{\text{RF}}$ and the linear system residual are based on nodes, rather than elements.  We define a diagonal matrix  $\BM_{\text{o}} \in \mathbb{R}^{N_s\times N_s}$ to represent the mask of overlap nodes $\mathcal{S}_{\text{o}}$ where  diagonal entries are defined as
\[
\left(\BM_{ \text{o}} \right)_{ii}= 
\left\{
\begin{aligned}
    1, & \quad  x_i\in \mathcal{S}_{\text{o}} \\
    0, & \quad\text{elsewhere}.
\end{aligned}
\right.
\]
Likewise, we define matrices  $\BM_{\text{i}}$ and $\BM_{\text{e}}$ for the interior subdomain and the exterior subdomain, respectively. We denote the restriction of  $\eta_{\text{RF}}$ in the overlap subdomain by $\eta_{\text{RF}}^{ \text{o}}$ 
\begin{equation*}
    \eta_{\text{RF}}^{ \Bw,\text{o}} := \| \BM_{ \text{o}} \BR_k   \|_{\Bw}  + \| \BM_{ \text{o}} \BF_k \|_{\Bw}.
\end{equation*}
The residual in the overlap subdomain is defined  as 
$$
\Br_k^{ \text{o}} = \BM_{ \text{o}} \Br_k.
$$
Analogously, we define $\eta_{\text{RF}}^{\text{i}}$ and $\Br_k^{ \text{i}}$ for the interior subdomain, $\eta_{\text{RF}}^{\text{e}}$ and   $\Br_k^{ \text{e}}$ for the exterior subdomain.
The subdomain-based stopping criterion is:
\[
    \|\Br_k^{ \text{i}}\|_{\Bw} \leq \tau \eta_{\text{RF}}^{\Bw,\text{i}}, \quad  \|\Br_k^{ \text{e}}\|_{\Bw} \leq \tau \eta_{\text{RF}}^{\Bw,\text{e}}, \text{ and }   \|\Br_k^{ \text{o}} \|_{\Bw} \leq \tau \eta_{\text{RF}}^{\Bw,\text{o}}.
\]

%% file: 6-1appendixProof.tex
\section{\added{Supporting Lemmas for the Proofs of \cref{thm:etaRFtoerror} and \cref{thm:errortoetaRF}}}
\label{sec:6-1}
\added{
In this section, we state several lemmas to support proofs of \cref{thm:etaRFtoerror} and \cref{thm:errortoetaRF}.
All lemmas share the same notations and assumption as in \cref{thm:etaRFtoerror}.
}
\added{
\subsection{Extremal eigenvalues of matrices}
In this subsection, we provide bounds for the largest and smallest eigenvalues of three matrices.
Define the matrix $\BM_{\ell}\in \mathbb{R}^{N_s\times N_s}$, whose $(i,j)$-th entry is
\begin{equation}
    \left(\BM_{\ell}\right)_{i,j} = \sum\limits_{\ell\in\mathcal{E}} \int_{\ell} \phi_j(s)\phi_i(s)ds.
    \label{eq:massl}
\end{equation}
We summarize the bounds on the smallest and largest eigenvalues of the stiffness matrix $\BA$, the mass matrix $\BM$, and the matrix $\BM_{\ell}$ \cite{canuto2007spectralsingle, canuto2007spectral, melenk2002condition}.
\begin{lemma}
    Let $\BA $ and $\BM$ be the stiffness matrix and the mass matrix, respectively, where $\BA_{i,j} = \left(\nabla\phi_j, \nabla\phi_i\right)$ and $\BM_{i,j} = \left(\phi_j, \phi_i\right)$. The matrix $\BM_{\ell}$ is defined in \cref{eq:massl}.
    Let  $\lammina$, $\lammaxa$, $\lammina$, $\lammaxa$, $\lambda_{\min}\left(\BM_{\ell}\right)$, and $\lambda_{\max}\left(\BM_{\ell}\right)$ denote the smallest and largest eigen-\\values of $\BM$, $\BA$ and $\BM_{\ell}$, respectively.
    Under \cref{assump}, there exist positive constants  $\underC_{\BA}$, $\overC_{\BA}$, $\underC_{\BM}$, $\overC_{\BM}$, $\underline{C}_{\ell}$, and $\overline{C}_{\ell}$ independent of polynomial degree $N$ and mesh size $h$, such that,
    \begin{eqnarray*}
        &\underC_{\BA} \frac{h^2}{N^2} \leq \lammina,\quad &\lammaxa \leq \overC_{\BA} N,\\   
        &\underC_{\BM} \frac{h^2}{N^4} \leq \lamminm,\quad &\lammaxm \leq \overC_{\BM}\frac{h^2}{N^2},\\
        &\underline{C}_{\ell} \frac{h}{N^2} \leq \lambda_{\min}\left(\BM_{\ell}\right),\quad &\lambda_{\max}\left(\BM_{\ell}\right) \leq \overline{C}_{\ell}\frac{h}{N}.
    \end{eqnarray*}
\label{lem:lambda}
\end{lemma}
\begin{remark}
    Similar bounds on $\lamminm$ and $\lammina$ can be verified for triangle elements with Warp \& Blend nodes numerically \cite{warburton2006explicit}. Numerical eigenvalues of matrices on the reference triangle are shown in \cref{fig:eigtri}. 
\end{remark}
\input{Data/figEigAM}

}
\added{
\subsection{Norm of residual and $\BA$-norm error}
Let $\Bw\in \mathbb{R}^{N_s}$ denote the vector corresponding to $w_h\in \Vhn$. The error is defined as $\Be_w = \Bx - \Bw$.
Due to $\Br(w_h) = \BA\Be_w$ and $\|\nabla(u_h - w_h)\|^2 = \Be_w^T\BA\Be_w$, we obtain the bound on the norm of the residual:
\begin{equation}
\lammina \|\nabla(u_h - w_h)\|^2 \leq \|\Br(w_h)\|_{l^2}^2 \leq \lammaxa \|\nabla(u_h - w_h)\|^2.
\label{res and error}
\end{equation}
}
\added{
\subsubsection{$\eta_{\text{RF}}(u_h^k)$ and $\eta_{\text{R}}(u_h^k)$}
Let $\Bv, \Bz \in \mathbb{R}^{N_s}$, where $\Bv_i = f_{h,N}(x_i) + \Delta w_h(x_i)$ and $ \Bz_i =  \Bn\cdot [ \nabla w_h](x_i)$. We find 
\[ 
\begin{aligned}
    \eta_{\text{RF}}^2(w_h) & = \|\BM\Bv\|^2_{l^2} +  \|\BM_{\ell}\Bz\|^2_{l^2}  =  \Bv^T\BM \BM\Bv +  \Bz^T\BM_{\ell} \BM_{\ell}\Bz.
\end{aligned}
\] 
Due to the quasi-uniform property of the triangulation, there exists constants $\underline{C}, \overline{C}>0$ such that
\[
\underC \left( \frac{h^2}{N^2}\Bv^T\BM\Bv +  \frac{h}{N}\Bz^T\BM_{\ell}\Bz \right)
\leq \eta_{\text{R}}^2(w_h) \leq \overC \left(\frac{h^2}{N^2}\Bv^T\BM\Bv +  \frac{h}{N}\Bz^T\BM_{\ell}\Bz\right).
\]
Therefore,
\begin{equation}
\begin{aligned}    
    \frac{1}{\overC}\min\left\{\lamminm \frac{N^2}{h^2}, \lamminml \frac{N}{h}\right\} \eta_R^2(w_h) & \leq 
    \eta_{\text{RF}}^2(w_h) \\
    & \leq  \frac{1}{\underC}\max\left\{\lammaxm \frac{N^2}{h^2}, \lammaxml \frac{N}{h}\right\} \eta_R^2(w_h).
    \label{etaRFetaR}
\end{aligned}
\end{equation}
Based on \cref{lem:lambda}, \cref{etaRFetaR} is equivalent to
\begin{equation}
\begin{aligned}    
    \lamminm \frac{N^2}{\overC h^2} \eta_R^2(w_h) \leq 
    \eta_{\text{RF}}^2(w_h) 
    \leq \lammaxm \frac{N^2}{\underC h^2} \eta_R^2(w_h).
    \label{etaRFetaR2}
\end{aligned}
\end{equation}
}
\added{
\subsubsection{ $\eta_{\text{R}}(w_h)$ and  $\|\nabla(u - w_h)\|^2$}
To obtain the relationship between the indicator $\eta_{\text{R}}(w_h)$ and the   total error $\|\nabla(u - w_h)\|^2$, we need three lemmas.
The first lemma states the connection between the indicator $\eta_{\text{R}}(u_h)$ and the discretization error $\|\nabla(u - u_h)\|^2$.
\begin{lemma}[Theorem 3.6 in \cite{melenk2001residual}]
Given $\varepsilon>0$, there exists $C_1, C_2>0$ independent of $h$ and $N$ such that
\begin{equation*}
\begin{aligned}
    & C_1 \|\nabla(u - u_h)\|^2 \leq \eta_R^2(u_h) +  \sum\limits_{K\in\mathcal{T}_h} \frac{h_K^2}{N^2} \|f_{h,K} - f\|_K^2  ,\\
    & \eta_R^2(u_h) \leq  C_2(\varepsilon) \left\{N^{2+2\varepsilon} \|\nabla(u - u_h)\|^2 + \sum\limits_{K\in\mathcal{T}_h} \frac{h_K^2}{N^{1-4\varepsilon}} \|f_{h,K} - f\|^2 \right\}.
\end{aligned}
\end{equation*}
\label{etaRbound}
\end{lemma}
The second lemma presents the basic inverse inequality  and the trace inequality from \cite{schwab1998p}.
\begin{lemma}[Theorem 4.76 in \cite{schwab1998p}]
There exists a constant $C>0$ independent of $N$, such that the following inequalities hold for any positive integer $N$ and for any polynomial $\phi_N \in \mathbb{Q}_N(\hatk)$: 
\[
\begin{aligned}
    \|\nabla\phi_N\|_{\hatk} \leq C N^2 \|\phi_N\|_{\hatk}, \quad
    \|\phi_N\|_{\ell} \leq C N \|\phi_N\|_{\hatk},
\end{aligned}
\]
where $\hatk = [-1,1]^2$ and $\ell$ is any one of edges of $\hatk$.
\label{inverse}
\end{lemma}
Moreover, we extend the result \cite[Lemma 3.1]{arioli2013stopping} to  clarify the $N$ dependence using \cref{inverse}.
\begin{lemma}
Let $v,w \in \Vhn$. Then there exists a constant $C>0$ independent of $h$ and $N$, such that
    \[
        \eta_R^2(v) \leq \eta_R^2(w) + CN^2 \|\nabla(v-w)\|^2.
    \]
    \label{lemma2}
\end{lemma}
\begin{proof}
Using the triangle inequality and  inequalities in \cref{inverse},
    \[
    \begin{aligned}
        \eta_R^2(v) & = \sum\limits_{K\in\mathcal{T}_h} \left(\frac{h_K^2}{N^2}\| f_h + \Delta w +  \Delta(v-w)\|^2_{K} + \sum\limits_{\ell\subset K\cap\left(\mathcal{E}_{int}\cup\mathcal{E}_{bd}^N\right)} \frac{h_K}{2N}\| [ \Bn\cdot\nabla w] + [ \Bn\cdot\nabla (v-w)] \|^2_{\ell}\right)  \\
        & \leq \eta_R^2(w) + 
         \sum\limits_{K\in\mathcal{T}_h} \left(\frac{h_K^2}{N^2} \| \Delta(v-w)\|^2_{K} + 
         \sum\limits_{\ell\subset K\cap\left(\mathcal{E}_{int}\cup\mathcal{E}_{bd}^N\right)}\frac{h_{\ell}}{2N}\| [ \Bn\cdot\nabla (v-w)] \|^2_{\ell} \right) \\
        & \leq  \eta_R^2(w) + \sum\limits_{K\in\mathcal{T}_h} \left( C_1\frac{h_K^2}{N^2} \frac{N^4}{h_K^2}\| \nabla(v-w)\|^2_{K} + 
        \sum\limits_{\ell\subset K\cap\left(\mathcal{E}_{int}\cup\mathcal{E}_{bd}^N\right)}
        \frac{h_{\ell}}{N}\| \nabla (v-w) \|^2_{\ell} \right) \\
        & \leq \eta_R^2(w) +  \sum\limits_{K\in\mathcal{T}_h} \left(C_1 N^2 \| \nabla(v-w)\|^2_{K} + C_2 \frac{h_{\ell}}{N} \frac{N^2}{h} \| \nabla (v-w) \|^2_{K}\right)  \\
        & \leq  \eta_R^2(w) + CN^2 \|\nabla(v-w)\|_{\Omega}^2, 
    \end{aligned}
    \]
    where $C_1$ and $C_2$ are positive constants derived from \cref{inverse} and the quasi-uniform property of the triangulation \cref{eq:quasi uniform}.
\end{proof}
Using \cref{etaRbound}, \cref{lemma2}, along with Galerkin orthogonality \cref{eq:galerkin orth}, we derive the relation-\\ship between the indicator $\eta_{\text{R}}(w_h)$ and the  total error $\|\nabla(u - w_h)\|^2$,
\begin{equation}
\begin{aligned}
\eta_R^2(w_h) & \leq \eta_R^2(u_h) + C_3N^2 \|\nabla(u_h -w_h)\|^2  \\
& \leq  C_4(\varepsilon) \left(N^{2+2\varepsilon} \|\nabla(u - u_h)\|^2 + \sum\limits_{K\in\mathcal{T}_h} \frac{h_K^2}{N^{1-4\varepsilon}} \|f_{h,K} - f\|_K^2 \right) +  C_3N^2\|\nabla(u_h -w_h)\|^2  \\
& \leq  \Tilde{C}_4(\varepsilon) N^{2+2\varepsilon} \left( \|\nabla(u - u_h)\|^2 + \|\nabla(u_h - w_h)\|^2 \right) +  C_4(\varepsilon) \sum\limits_{K\in\mathcal{T}_h} \frac{h_K^2}{N^{1-4\varepsilon}} \|f_{h,K} - f\|_K^2  \\
& \leq   \Tilde{C}_4(\varepsilon) N^{2+2\varepsilon} \|\nabla(u - w_h)\|^2  +  C_4(\varepsilon) \sum\limits_{K\in\mathcal{T}_h} \frac{h_K^2}{N^{1-4\varepsilon}} \|f_{h,K} - f\|_K^2. 
\end{aligned}
\label{etaRerror}
\end{equation}
Here $C_3$ is derived from \cref{lemma2}, $C_4(\varepsilon) $ is from \cref{etaRbound}, and $\Tilde{C}_4(\varepsilon) = \max(C_3, C_4(\varepsilon))$.
}

%% file: Data/figEigAM.tex
\begin{figure}
\centering
\begin{tikzpicture}
\begin{axis}[
    xlabel={Polynomial degree},
    ylabel={eigenvalue},
    xtick=data,
    xticklabels={2,3,4,5,6,8,, ,11,,,14,},
    xmode=log,
    ymode=log,
    legend pos= outer north east,
    ymajorgrids=true,
    grid style=dashed,
]
\addplot[
    color=blue,
    mark=star,
    ]
    coordinates {
(1, 1.5) 
(2, 4.74208) 
(3, 10.7221) 
(4, 19.8354) 
(5, 37.2156) 
(6, 75.9954) 
(7, 170.359) 
(8, 418.229) 
(9, 1107.09) 
(10, 3114) 
(11, 9167.63) 
(12, 27961.9) 
(13, 87631.2) 
    };
    \addlegendentry{$\lambda_{\max}(\BA)$}
\addplot[
    color=red,
    mark=square,
    ]
    coordinates {
(1, 0.5) 
(2, 0.311305) 
(3, 0.219056) 
(4, 0.176351) 
(5, 0.141426) 
(6, 0.117046) 
(7, 0.0970199) 
(8, 0.0818189) 
(9, 0.069476) 
(10, 0.0596922) 
(11, 0.0517123) 
(12, 0.0451905) 
(13, 0.0397986) 
    };
    \addlegendentry{$\lambda_{\min}(\BA)$}

\addplot[
    color=red,
    dashed,
]
    coordinates {

(1, 1.01018) 
(2, 0.505088) 
(3, 0.303053) 
(4, 0.202035) 
(5, 0.144311) 
(6, 0.108233) 
(7, 0.0841814) 
(8, 0.0673451) 
(9, 0.0551005) 
(10, 0.0459171) 
(11, 0.0388529) 
(12, 0.0333025) 
(13, 0.0288622) 
    };
    \addlegendentry{$N_p$}
\addplot[
    color=cyan,
    mark=triangle,
    ]
    coordinates {
(1, 0.0416667) 
(2, 0.0446288) 
(3, 0.0379003) 
(4, 0.0280971) 
(5, 0.0275884) 
(6, 0.0316924) 
(7, 0.0458138) 
(8, 0.0743847) 
(9, 0.144546) 
(10, 0.300661) 
(11, 0.702936) 
(12, 1.75406) 
(13, 4.40352)  
    };
    \addlegendentry{$\lambda_{\max}(\BM)$}
\addplot[
    color=black,
    mark=diamond,
    ]
    coordinates {
(1, 0.0104167) 
(2, 0.00259341) 
(3, 0.00111572) 
(4, 0.000486579) 
(5, 0.000266901) 
(6, 0.000147498) 
(7, 9.26337e-05) 
(8, 5.86049e-05) 
(9, 4.02201e-05) 
(10, 2.77106e-05) 
(11, 2.01998e-05) 
(12, 1.47389e-05) 
(13, 1.12122e-05) 
    };
    \addlegendentry{$\lambda_{\min}(\BM)$}

\addplot[
    color=black,
    densely dotted,
]
    coordinates {
(1, 0.0103169) 
(2, 0.00257922) 
(3, 0.000928518) 
(4, 0.000412675) 
(5, 0.000210548) 
(6, 0.000118433) 
(7, 7.16449e-05) 
(8, 4.58528e-05) 
(9, 3.06948e-05) 
(10, 2.13158e-05) 
(11, 1.52616e-05) 
(12, 1.12126e-05) 
(13, 8.42193e-06) 
    };
    \addlegendentry{$N_p^{-2}$}
\end{axis}
\end{tikzpicture}
\caption{The largest and smallest eigenvalues of the stiffness matrix ${\BA}$ and the mass matrix $\BM$ on the reference triangle element using Warp \& Blend nodes \cite{warburton2006explicit}. $N_p=(N+1)(N+2)/2$.}
\label{fig:eigtri}
\end{figure}